\RequirePackage{fix-cm}

\documentclass[msom,nonblindrev]{informs4}

\RequirePackage{bm}
\RequirePackage{bbm}
\RequirePackage{endnotes}

\OneAndAHalfSpacedXI

\makeatletter
\def\theARTICLEABSTRACT{\HOOKb
  \vspace*{18pt}
  \begin{minipage}[t]{\textwidth}\parindent1em
    \ABSfont
    \noindent\theABSTRACT\endgraf
    \vskip5pt
    \theFUNDING
    \theKEYWORDS
    \theSUBJECTCLASS
    \theAREAOFREVIEW
    \theMSCCLASS
    \theORMSCLASS
    \if@BLINDREV\else\theHISTORY\fi
    \noindent\hrulefill
  \end{minipage}\vspace*{0pt}}
\makeatother

\TheoremsNumberedThrough
\ECRepeatTheorems
\EquationsNumberedThrough

\usepackage{graphicx}
\usepackage{float}
\usepackage{placeins}
\usepackage{subcaption}

\allowdisplaybreaks[4]

\usepackage{longtable}
\usepackage{booktabs}
\usepackage{tabularx}
\usepackage{multirow}
\usepackage{threeparttable}
\usepackage{makecell}
\usepackage{rotating}

\usepackage{algorithm}
\usepackage[noend]{algpseudocode}

\usepackage{natbib}

\usepackage{hyperref}

\bibpunct[, ]{(}{)}{,}{a}{}{,}

\definecolor{msomlinkred}{RGB}{180,0,0}
\definecolor{msomciteblue}{RGB}{0,0,170}
\definecolor{msomurlblue}{RGB}{0,0,204}
\definecolor{msomecblue}{RGB}{0,0,166}
\hypersetup{
    colorlinks=true,
    linkcolor=msomlinkred,
    citecolor=msomciteblue,
    urlcolor=msomurlblue
}
\newcommand{\ECsection}[1]{\textcolor{msomecblue}{Section~#1 of the e-companion}}

\makeatletter

\@ifundefined{theorem2}
  {}
  {}

\@ifundefined{observation}
  {}
  {}

\makeatother

\newcommand{\bx}{\boldsymbol{x}}
\newcommand{\bw}{\boldsymbol{w}}
\newcommand{\bb}{\boldsymbol{b}}
\newcommand{\bt}{\boldsymbol{t}}
\newcommand{\by}{\boldsymbol{y}}
\newcommand{\bn}{\boldsymbol{n}}
\newcommand{\br}{\boldsymbol{r}}
\newcommand{\bs}{\boldsymbol{s}}

\newcommand{\bh}{\boldsymbol{h}}
\newcommand{\bg}{\boldsymbol{g}}
\newcommand{\bq}{\boldsymbol{q}}
\newcommand{\bomega}{\boldsymbol{\omega}}
\newcommand{\bN}{\boldsymbol{N}}
\newcommand{\bbeta}{\bm{\beta}}

\newcommand{\bmu}{\boldsymbol{\mu}}

\newcommand{\bpi}{\boldsymbol{\pi}}

\newcommand{\bD}{\boldsymbol{D}}
\newcommand{\bP}{\boldsymbol{P}}

\newcommand{\bG}{\boldsymbol{G}}

\algdef{SE}[DOWHILE]{Do}{doWhile}
  {\algorithmicdo}
  [1]{\algorithmicwhile\ #1}

\makeatletter
\newcommand*{\centerfloat}{\parindent\z@
    \leftskip\z@\@plus 1fil\@minus\marginparwidth
    \rightskip\leftskip
    \parfillskip\z@skip
}
\makeatother

\definecolor{upgreen}{RGB}{0,128,0}
\definecolor{downred}{RGB}{180,0,0}

\newcommand{\gain}[1]{\textcolor{upgreen}{\(\uparrow\)\hspace{0.20em}#1\%
    }}

\newcommand{\compacttablerows}{\renewcommand{\arraystretch}{1.15}\setlength{\extrarowheight}{0pt}\setlength{\aboverulesep}{0.10ex}\setlength{\belowrulesep}{0.10ex}}
\AtBeginEnvironment{gather}{\compactalignspacing}

\newcolumntype{R}{>{\centering\arraybackslash}p{0.140\textwidth}}
\newcolumntype{K}{>{\centering\arraybackslash}p{0.040\textwidth}}
\newcolumntype{O}{>{\centering\arraybackslash}p{0.043\textwidth}}
\newcolumntype{T}{>{\centering\arraybackslash}p{0.048\textwidth}}
\newcolumntype{I}{>{\centering\arraybackslash}p{0.050\textwidth}}
\newcolumntype{G}{>{\centering\arraybackslash}p{0.070\textwidth}}

\newcommand{\compacttableend}{\vspace{-3pt}}
\newcommand{\compactfigureend}{\vspace{-8pt}}

\usepackage{etoolbox}
\providecommand{\best}[1]{\textbf{#1}}
\newcommand{\compactalignspacing}{\setlength{\abovedisplayskip}{3pt plus 1pt minus 1pt}\setlength{\belowdisplayskip}{3pt plus 1pt minus 1pt}\setlength{\abovedisplayshortskip}{2pt plus 1pt}\setlength{\belowdisplayshortskip}{3pt plus 1pt minus 1pt}}

\newcommand{\compactequationspacing}{\setlength{\abovedisplayskip}{5pt plus 1pt minus 1pt}\setlength{\belowdisplayskip}{5pt plus 1pt minus 1pt}\setlength{\abovedisplayshortskip}{3pt plus 1pt}\setlength{\belowdisplayshortskip}{5pt plus 1pt minus 1pt}}

\AtBeginEnvironment{equation}{\compactequationspacing}
\AtBeginEnvironment{equation*}{\compactequationspacing}
\AtBeginEnvironment{align}{\compactalignspacing}
\AtBeginEnvironment{align*}{\compactalignspacing}

\begin{document}

\RUNAUTHOR{Xu et al.}

\RUNTITLE{Driver Surge Maps with Downstream Matching}

\TITLE{Driver Surge Maps for Supply Rebalancing with Downstream Matching}

\ARTICLEAUTHORS{\AUTHOR{Yida Xu\textsuperscript{1,\(\dagger\)},
Baichuan Mo\textsuperscript{1,\(\dagger\),*},
Zhengzhong Ricky You\textsuperscript{1,\(\dagger\)},
Yuan Zhou\textsuperscript{2},
Hao Yi Ong\textsuperscript{3},
Zuo-Jun Max Shen\textsuperscript{4}}
\AFF{\textsuperscript{1}Department of Civil Engineering, Tsinghua University, Beijing 100084, China}
\AFF{\textsuperscript{2}Yau Mathematical Sciences Center, Tsinghua University, Beijing 100084, China}
\AFF{\textsuperscript{3}Lyft Inc., San Francisco, CA 94107, USA}
\AFF{\textsuperscript{4}Department of Industrial Engineering and Operations Research, University of California, Berkeley, Berkeley, CA 94720, USA}
\AFF{\textsuperscript{\(\dagger\)}These authors contributed equally to this work. \textsuperscript{*}Corresponding author: \EMAIL{bmo@tsinghua.edu.cn}.}
}

\ABSTRACT{\textbf{Problem definition:}
Ride-hailing platforms use driver-facing surge maps to rebalance spatial supply by offering zone-level bonuses.
These bonuses can reshape supply by inducing nearby drivers to relocate and encouraging latent drivers to become dispatchable, but both responses are voluntary.
Within short decision windows, platforms risk offering bonuses that reshape driver supply without generating sufficient downstream matching value to justify the associated payments.
We therefore study real-time driver surge map design when voluntary driver responses and downstream matching jointly determine the value of zone-level bonuses.
\textbf{Methodology/results:}
We formulate the driver surge map design problem (DSMDP) to optimize zone-level bonuses while anticipating voluntary driver responses and downstream matching.
To assess how bonus adjustments affect net operating profit, we conduct a marginal analysis that identifies when small bonuses improve on the no-bonus map and how spatial smoothness links bonus decisions across zones.
Using these marginal values, we develop an iterative solution approach (ISA) to sample dual sensitivities of the matching linear program across nearby matching regimes and guide real-time bonus updates.
Theoretical analysis provides certificates of approximate stationarity and establishes conditions for stationary limits.
Computational experiments show that ISA improves net operating profit over the best-performing feasible benchmark by up to 19.37\% on the generated instances and by up to 8.43\% in the New York City study.
\textbf{Managerial implications:}
Local shortage alone is not a sufficient basis for bonus allocation.
Platforms should target responsive supply whose downstream matching value justifies the associated bonus payments.
Latent-driver activation can support higher-value matches even when all demand is served, but the resulting gains must also cover additional payments to drivers who would already be matched without incentives.
Results indicate diminishing returns to bonus spending and show that moderate spatial smoothness can preserve most of the gains from spatial targeting.
}

\maketitle

\section{Introduction}
Ride-hailing platforms such as Lyft and Uber perform millions of real-time rider--driver matches every day \citep{garg2022driver}.
Yet spatial imbalances between rider demand and available drivers remain a persistent operational problem, resulting in long waits, inefficient pickups, and lost service opportunities.
To improve matching success rates and mitigate spatial supply-demand imbalances, ride-hailing platforms adopt operational strategies on both the demand and supply sides. 
On the demand side, platforms
use dynamic rider pricing to regulate realized demand, and prior research has examined its implications
\citep{banerjee2015pricing,cachon2017role,chen2024real}.
On the supply side, platforms use spatial incentives to rebalance driver supply. Existing work studies how spatial prices and incentives reshape drivers' relocation and participation decisions and the resulting supply distribution \citep{afeche2023ride,besbes2021surge,bimpikis2019spatial}. 
In these models, the direction of influence runs primarily from incentives to driver responses and downstream matching outcomes, while the matching value of the induced supply does not in turn inform the incentive decision.
This feedback naturally raises the question of how the platform should construct the spatial signal shown to drivers.

\begin{figure}[!htbp]
        \centering
        \includegraphics[width=0.7\linewidth]{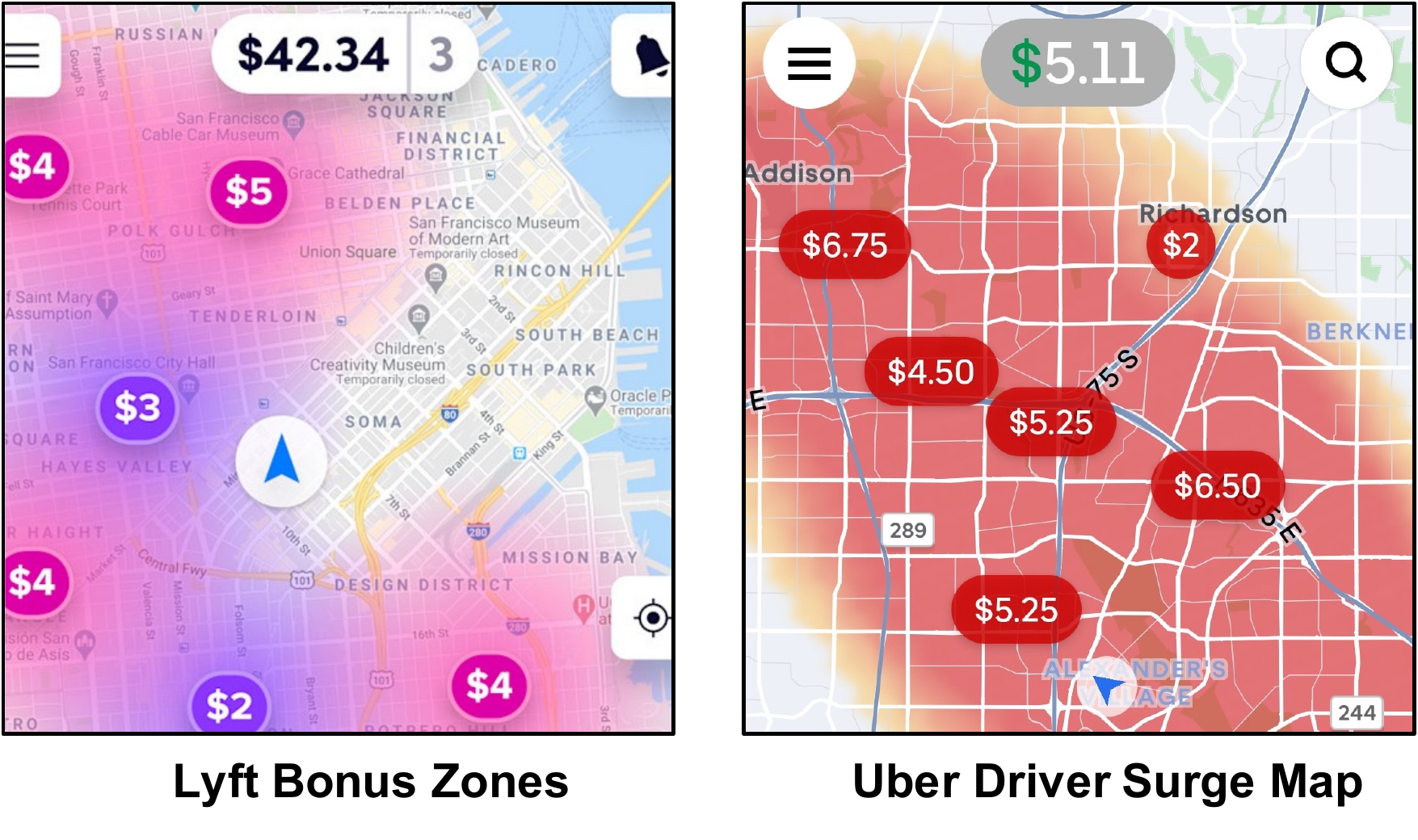}
    \caption{Examples of supply-side spatial incentive maps on ride-hailing platforms.}
\label{fig:driver_spatial_incentive_products}
\vspace{-10pt}
\end{figure}

To influence drivers' relocation and participation decisions, one can provide spatial information, such as demand heatmaps or repositioning recommendations, to guide drivers toward particular areas without explicitly specifying a monetary reward \citep{haferkamp2024heatmap, alnaggar2025heatmap}. 
Another approach is to offer zone-level monetary incentives through a driver-facing surge map, as considered in this study. 
Specifically, the platform may assign zone-level bonuses based on the current supply--demand imbalances, and matched drivers who move to these zones receive the corresponding bonus.
Such monetary incentive maps are widely used by ride-hailing platforms.
For example, Uber displays driver surge maps in which areas with higher rider demand and limited driver supply are highlighted by stronger colors and associated with higher potential earnings \citep{uber_driver_surge}.
Similarly, Lyft provides Bonus Zones, where drivers can unlock additional earnings by entering designated high-demand areas \citep{lyft_bonus_zones}.
Figure~\ref{fig:driver_spatial_incentive_products} shows examples of driver-facing maps.
Surge maps can reshape supply through two voluntary-response channels: (1) attracting nearby online drivers to voluntarily reposition toward undersupplied areas, and (2) encouraging nearby offline drivers, or drivers from competing platforms, to become active on the platform.

These voluntary-response channels have been studied through spatial pricing and compensation mechanisms \citep{bimpikis2019spatial,besbes2021surge}, incentive-compatible driver payments \citep{garg2022driver}, and real-time positioning incentives with spatial budgets \citep{ong2021driver}.
Explicitly incorporating voluntary driver responses into the model can be challenging, since drivers may choose among a large set of potential relocation destinations.
In real-time operations, incentive decisions are updated frequently; for example, Lyft's deployed positioning system generates fresh driver incentives every 2 minutes \citep{ong2021driver}. 
This short decision horizon makes it natural to restrict each driver's relocation alternatives to nearby reachable locations. 
Consequently, driver responses over the resulting finite choice set can then be represented using discrete choice models.
However, embedding such choice probabilities into the model introduces nonlinear terms, making the problem difficult to solve within a real-time operational window.


Motivated by these operational features, we study how a platform should design a zone-level driver surge map to maximize net operating profit within a short operational window.
We formulate the driver surge map design problem (DSMDP), in which the platform directly optimizes zone-level bonuses while anticipating voluntary driver responses and optimal downstream matching, so that bonus allocation reflects the matching value of induced supply and the additional payments on existing matches.
We use \textit{latent drivers} to include drivers who are offline or serving competing platforms.
To capture voluntary driver responses, we use a multinomial logit (MNL) model for dispatchable drivers' choices among nearby zones and an unavailable state, and a binary logit model for latent drivers to become available for dispatch or not within their current zones.
These responses determine the expected supply for matching, and bonuses are paid to matched drivers according to their post-response zones.
To determine how downstream matching should inform bonus decisions, we analyze the matching linear program (LP) obtained by fixing the surge map.

\begin{figure}[H]
		\centering
\includegraphics[width=0.88\linewidth]{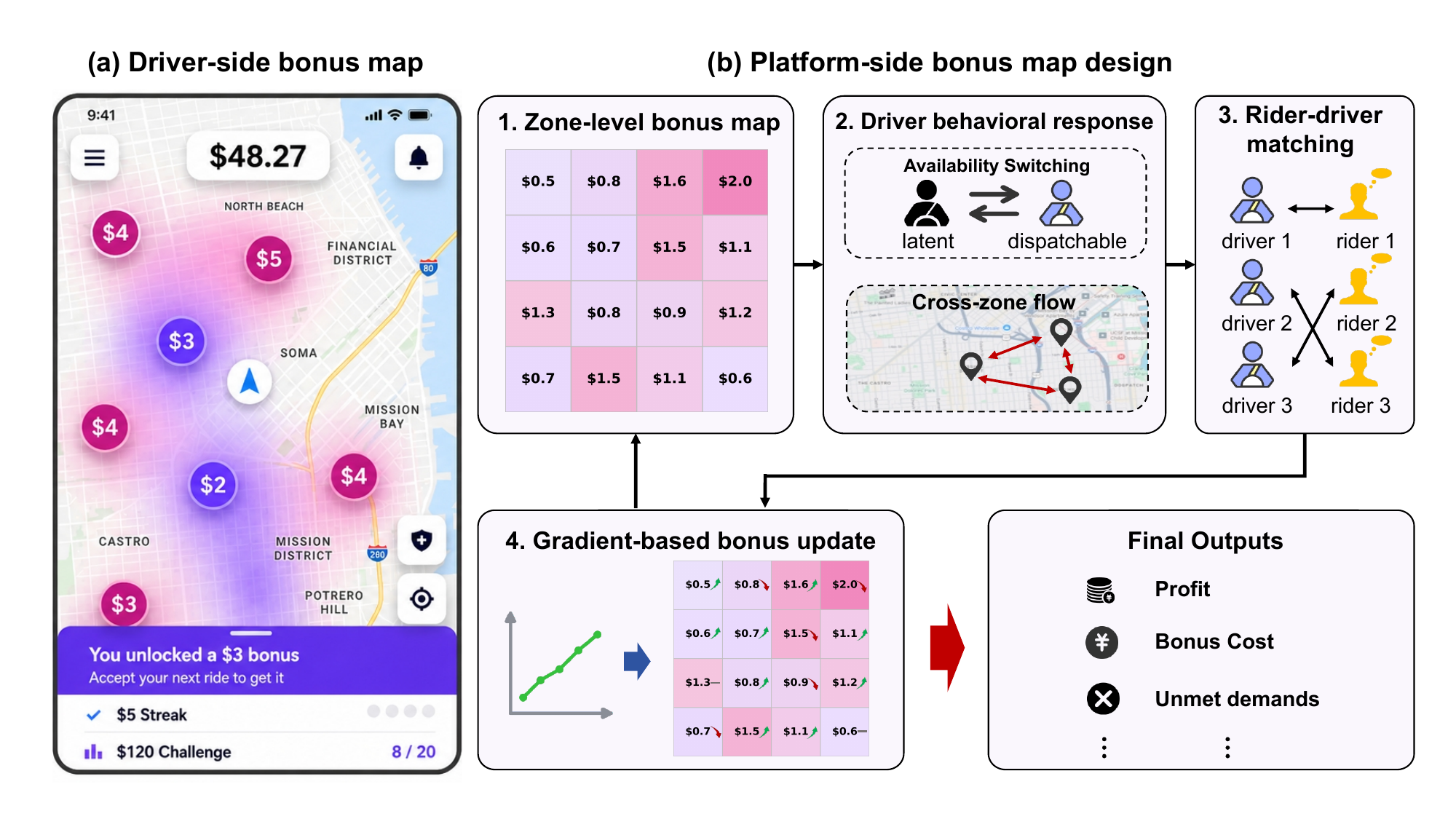}
\caption{Framework of the driver surge map design problem.}
\label{fig_framework}
\compactfigureend
\end{figure}

To reflect the interaction between bonus decisions, driver responses, and downstream matching, we develop an iterative solution approach (ISA) that uses dual information from the matching problem to guide bonus adjustments.
At each iteration, ISA fixes the zone-level bonuses, evaluates driver responses using the MNL and binary logit models, and solves the matching problem with the updated supply.
Dual information from the matching problem then guides bonus adjustments for the next iteration.
This feedback allows us to examine how marginal changes in bonuses affect driver responses and, in turn, net operating profit.
Figure~\ref{fig_framework} summarizes the operational sequence.
Our analysis provides finite-scale stationarity certificates
and sufficient conditions for constrained Clarke stationarity
of sampled limits.

We then examine the operational value of driver surge maps and show that offering bonuses improves net operating profit relative to the no-bonus map across the tested settings.
The gains are greatest when available driver supply is scarce, while returns diminish as bonus spending increases.
Furthermore, stronger latent-driver activation can increase net operating profit even when all demand is already served, as the additional supply supports higher-value matches.
This finding highlights the role of downstream matching opportunities in determining the value of attracting drivers.
Finally, moderate limits on bonus differences between neighboring zones preserve most of the gains from spatial targeting, suggesting that effective bonus allocation does not require sharp differences across adjacent zones.

We make three contributions.

\begin{itemize}

\item \emph{Model formulation and structural analysis.}
We formulate the DSMDP to optimize zone-level bonuses while anticipating voluntary driver responses and downstream matching.
The marginal analysis reveals how the matching value of induced supply, payments on existing matches, and budget opportunity cost shape bonus allocation, explaining why zones with similar local shortages may warrant different bonuses.
We derive conditions under which introducing bonuses improves net operating profit and establish how spatial smoothness couples marginal values across zones.

\item \emph{Real-time solution and theoretical analysis.}
We develop ISA to translate the marginal-value characterization into bonus updates within the real-time decision window.
The method combines matching-LP dual information with sensitivities from nearby matching regimes to address nonsmoothness in bonus optimization.
We establish finite-scale stationarity certificates and sufficient conditions for constrained Clarke stationarity of sampled limits.
A conic outer relaxation provides an upper bound on the optimal net operating profit.

\item \emph{Computational performance and operational insights.}
ISA achieves the highest net operating profit among the evaluated methods, outperforming the best-performing feasible benchmark by up to 19.37\% on the generated instances and by 2.38\%--8.43\% in the New York City case study.
Comparisons with the no-bonus map quantify the value of intervention, while sensitivity analyses examine how value varies with budget, driver responsiveness, behavioral-parameter estimation error, and spatial smoothness.

\end{itemize}

After reviewing the literature in Section~\ref{lit}, we present the formulation of DSMDP in Section~\ref{problem}. Section~\ref{sec:matching_value} shows how downstream matching determines marginal bonus value, and Section~\ref{theories} presents ISA and its theoretical analysis.
Section~\ref{exps} then reports computational performance and provides operational insights.
Finally, Section~\ref{sec_conclusion} concludes the paper.

\section{Literature Review} \label{lit}

Our study is closely related to three streams of literature: surge pricing and driver bonus design, driver behavioral response and spatial relocation, and real-time matching and optimization.

\textbf{Surge pricing and driver bonus design.}
Surge pricing and driver bonuses are central operational tools for managing supply--demand imbalance in ride-hailing platforms. A stream of studies has investigated how pricing, matching, and bonuses balance rider demand and driver supply. \citet{banerjee2015pricing} analyzed pricing effects on revenue and throughput, while \citet{cachon2017role} showed that state-dependent prices and payments improve provider utilization and market performance. From a system-level perspective, 
\citet{yan2020dynamic} reviewed dynamic pricing and matching, while \citet{castillo2025matching} showed that ride-hailing markets may suffer from ``wild goose chase'' failures, in which high demand, few idle drivers, and long pickup times reinforce each other and lead to substantial efficiency losses.

In ride-hailing markets, supply--demand imbalance is often spatial: some zones face excess rider demand while nearby zones have idle drivers. This spatial mismatch motivates a closely related stream on spatial pricing and driver bonus design. Rider-side pricing primarily regulates demand, while bonuses change drivers' earnings and influence labor supply, relocation, and trip acceptance \citep{chen2016dynamic}. \citet{bimpikis2019spatial} analyzed driver participation and location choices, while \citet{ma2022spatio} developed spatio-temporal pricing mechanisms. \citet{besbes2021surge} showed that local demand shocks can induce global relocation, while \citet{guda2019your} demonstrated that surge pricing and forecast communication guide movement across locations. These studies show how spatial pricing links local conditions to driver responses across the network. Beyond spatial relocation, bonus rules also shape driver earnings and strategic participation. \citet{garg2022driver} showed that additive surge can improve incentive compatibility over multiplicative surge, while \citet{bai2023hiding} showed that providers may limit online availability to induce surge pricing.

Lyft's Personal Power Zone (PPZ) system provides a close operational analogue to our setting. Specifically, it allocates real-time location-based positioning incentives under spatial budgets and computes driver-positioning decisions within seconds \citep{ong2021driver}. Whereas their approach computes bonuses after determining driver-positioning allocations, our DSMDP directly optimizes a driver surge map while anticipating voluntary driver relocation and changes in dispatch availability. As these responses reshape driver supply, downstream matching reveals where additional supply creates sufficient value to justify bonuses after accounting for payments on existing matches and budget opportunity cost, thereby guiding their allocation across zones.

\textbf{Driver behavioral response and spatial relocation.}
A key feature of ride-hailing systems is that driver supply responds behaviorally to bonuses. Drivers may decide whether to go online, remain online, reposition, or accept trips. \citet{lu2018surge} found that surge heatmap visibility affected Uber drivers' relocation and earnings. Beyond relocation, \citet{chen2016dynamic} examined labor supply responses, while \citet{miao2023effects} found heterogeneity by employment type, experience, and preferences. Driver repositioning has also been studied from behavioral and operational perspectives. \citet{afeche2023ride} examined platform controls with strategic drivers, while \citet{wang2025strategic} characterized when spontaneous relocation mitigates imbalance and when subsidies are needed. At the individual level, \citet{urata2021learning} modeled zonal search through discrete choice, while \citet{liu2026street} studied outside opportunities and platform participation. These studies account for how monetary incentives, location, and outside earning opportunities affect driver availability and relocation.

Another line of work studies heatmaps and repositioning recommendations that influence relocation without monetary rewards. \citet{alnaggar2025heatmap} and \citet{haferkamp2024heatmap} showed that heatmaps can improve spatial supply balance. In our setting, monetary zone-level bonuses affect relocation and dispatch availability. To represent these voluntary responses, we use a multinomial logit (MNL) model for dispatchable drivers' choices among reachable zones and an unavailable state, and a binary logit model for latent-driver activation within the current zone. These responses determine expected post-response supply for downstream matching.

\textbf{Real-time matching, dispatching, and optimization methods.}
Real-time rider--driver matching and order dispatching form another major stream of ride-hailing research. Matching determines rider--driver assignments and affects waiting time, utilization, revenue, and future supply. \citet{ozkan2020dynamic} studied dynamic matching using linear programs with spatial heterogeneity and waiting constraints. \citet{yan2020dynamic} reviewed matching and pricing under real-time requirements. Several studies developed scalable dispatching algorithms. \citet{xu2018large} proposed learning-and-planning framework incorporating future value, while \citet{qin2020ride} described sequential dispatching at DiDi. These studies show that real-time dispatching requires scalable algorithms accounting for immediate and future performance.

Another related stream considered vehicle repositioning and rebalancing. \citet{braverman2019empty} studied empty-car routing, while \citet{tuncel2023integrated} and \citet{gao2024online} integrated repositioning with matching. In DSMDP, bonuses influence voluntary movements and dispatch availability before supply enters matching. Our study examines how downstream matching should inform bonus decisions under voluntary supply response. For a fixed driver surge map, the matching LP evaluates post-response supply while accounting for bonus payments and the available budget. Combined with the supply response model, its dual information characterizes the marginal value of zone-level bonuses and guides feasible updates in ISA.

\section{Problem Description and Formulation} \label{problem}
This section introduces the DSMDP and develops its mathematical formulation.
We first describe the operational setting in Section~\ref{subsec:problem_description}, then model driver responses to posted bonuses and incorporate the resulting supply into the platform's optimization problem in Section~\ref{subsec:mathematical_formulation}.

\subsection{Problem Description}\label{subsec:problem_description}
We consider a ride-hailing system operating over discrete real-time decision periods.
At the beginning of each period, the platform observes rider demand and driver availability and posts a driver surge map, to which drivers respond through voluntary relocation and changes in dispatch availability.
The downstream matching system then assigns available drivers to rider requests, and bonuses are paid to matched drivers based on their post-response zones.
In practice, this sequence is repeated frequently, e.g., every 2 minutes.
For simplicity, we suppress the period index in the remainder of this section.

The service region is partitioned into a finite set of zones denoted by $\mathcal{Z}$, with cardinality $|\mathcal{Z}|$.
Each zone serves as the basic spatial unit for displaying driver-facing bonuses and conducting real-time matching operations.
In practice, platforms can define these units using a fine-grained geospatial grid, such as the cells used in Lyft's driver-positioning system \citep{ong2021driver}.

The current system state is characterized by zone-level rider demand and driver supply.
Specifically, $d_j$ denotes the number of ride requests awaiting service in zone $j\in\mathcal{Z}$, while $s_i^{\mathrm{On}}$ and $s_i^{\mathrm{Off}}$ denote the dispatchable and latent driver supplies in zone $i\in\mathcal{Z}$, respectively.
Dispatchable drivers are immediately available for matching, whereas latent drivers may view the surge map and become available in response to incentives.
We assume that their approximate locations are available from recent application activity or location signals.
Since driver bonuses can induce short-term relocation, let $\mathcal{Z}_i\subseteq\mathcal{Z}$ denote the set of zones that a driver initially located in zone $i$ can reach within one decision period.
We assume $i\in\mathcal{Z}_i$ and a symmetric reachability relation, under which $j\in\mathcal{Z}_i$ if and only if $i\in\mathcal{Z}_j$, as satisfied by the grid neighborhoods used in the experiments.
This neighborhood captures the spatial scope over which the platform can reshape driver supply in real time.

Given the observed system state, the platform may regulate the spatial distribution of driver supply by posting zone-level monetary incentives.
To this end, it posts a nonnegative driver-side \emph{surge bonus} $b_i\ge 0$ for each zone $i\in\mathcal{Z}$.
Following the \emph{additive surge} structure described by \citet{garg2022driver}, the bonus is attached to the service zone and is independent of the travel distance or travel time of a particular trip.
We collect these zone-level bonuses in the vector $\bb=\{b_i\}_{i\in\mathcal{Z}}$, which represents the driver surge map displayed by the platform.

Designing the driver surge map requires the platform to anticipate how the offered bonuses will reshape driver supply.
Specifically, dispatchable drivers may relocate across nearby zones or become unavailable for dispatch (e.g., switch to a different platform), whereas latent drivers may become available for service.
These behavioral responses alter the spatial supply distribution available to the matching system and therefore the operational value of the driver surge map.
We represent the downstream matching outcome by the matching vector $\bx=\{x_{ij}\}_{i,j\in\mathcal{Z}}$, where $x_{ij}$ denotes the number of available drivers in zone $i$ assigned to ride requests in zone $j$ during the current period.

The platform selects the driver surge map $\bb$ to maximize net operating profit while anticipating the optimal downstream matching outcome $\bx$.
For a match between a driver in zone $i$ and a rider in zone $j$, let $p_{ij}$ denote the pre-bonus matching value, which accounts for rider-side revenue, baseline driver compensation, and other non-bonus costs.
Thus, the net contribution of this match is $p_{ij}-b_i$, since bonus eligibility is determined by the driver's post-response zone at the time of dispatch rather than the rider's request zone.
Although ride-hailing operations are inherently dynamic, the DSMDP can be formulated as a static optimization problem over a single decision period.
When forecasts of future system conditions are available, their effects can be incorporated through opportunity-cost adjustments to the matching values $p_{ij}$ \citep{qin2020ride}, thereby preserving the static structure of the formulation.

\subsection{Mathematical Formulation}
\label{subsec:mathematical_formulation}
To formulate the DSMDP, we characterize the post-response driver supply available for matching under a given driver surge map.
Surge bonuses affect driver supply through two behavioral channels.
First, a dispatchable driver who observes the map may remain in the current zone, relocate to a reachable zone, or become unavailable for dispatch during the current period.
For any $i\in\mathcal{Z}$ and $j\in\mathcal{Z}_i$, let $P_{ij}(\bb)$ denote the probability that a dispatchable driver initially located in zone $i$ is available in zone $j$ after responding to the map.
Let $P_i^{\mathrm{Off}}(\bb)$ denote the probability of becoming unavailable for dispatch; since $i\in\mathcal Z_i$, the choice set includes remaining in the current zone.

We model these choices with an MNL model whose alternatives are the reachable service zones and an unavailable state.
Each driver therefore compares the expected attractiveness of remaining dispatchable in a reachable zone with that of leaving the dispatchable pool.
We use a location-normalized utility representation so that the coefficient combinations governing these choices can be identified from observed transitions.
For a driver initially located in zone $i$, the unnormalized deterministic utility of choosing destination zone $j\in\mathcal Z_i$ is
\(\beta^{(0)}+\beta^{(1)}b_j+\beta^{(2)}b_i+\beta^{(3)}d_{ij}\), where $b_j$ is the destination-zone bonus, $b_i$ is the origin-zone bonus, and $d_{ij}$ is the travel distance or relocation cost from $i$ to $j$.
The parameter $\beta^{(0)}$ is the common intercept of the dispatchable-zone alternatives, $\beta^{(1)}$ measures sensitivity to the destination-zone bonus, $\beta^{(2)}$ captures the effect of the origin-zone bonus on remaining dispatchable, and $\beta^{(3)}$ measures the disutility of relocation distance.
The unnormalized utility of entering the unavailable state is
\(\alpha_i^{\mathrm{Off}}-\gamma^{\mathrm{Off}}b_i\)
for each \(i\in\mathcal Z\), where $\alpha_i^{\mathrm{Off}}$ is the zone-specific baseline attractiveness of becoming unavailable and $\gamma^{\mathrm{Off}}\ge 0$ measures how the origin-zone bonus reduces the attractiveness of this state.

Because MNL choice probabilities are invariant to subtracting the same
origin-specific term from all alternatives, subtracting
\(\beta^{(0)}+\beta^{(2)}b_i\)
yields the normalized utilities
\(V_{ij}=\beta^{(1)}b_j+\beta^{(3)}d_{ij}\)
and
\(V_i^{\mathrm{Off}}
=\delta_i^{\mathrm{Off}}-\zeta^{\mathrm{Off}}b_i\), where
\(\delta_i^{\mathrm{Off}}
=\alpha_i^{\mathrm{Off}}-\beta^{(0)}\)
denotes the zone-specific relative intercept of the unavailable state and
\(\zeta^{\mathrm{Off}}
=\beta^{(2)}+\gamma^{\mathrm{Off}}\).
We impose \(\zeta^{\mathrm{Off}}\ge0\), so that a higher origin-zone bonus weakly reduces the relative attractiveness of becoming unavailable.
For simplicity, let
\(\bbeta=(\beta^{(1)},\beta^{(3)})\).
This reparameterization preserves all choice probabilities and
allows the identified coefficient combinations to be estimated
directly.
For each \(i\in\mathcal Z\) and \(j\in\mathcal Z_i\), the relocation and
unavailability probabilities are
\(P_{ij}(\bb)=\exp(V_{ij})/
[\sum_{m\in\mathcal Z_i}\exp(V_{im})+\exp(V_i^{\mathrm{Off}})]\)
and
\(P_i^{\mathrm{Off}}(\bb)=\exp(V_i^{\mathrm{Off}})/
[\sum_{m\in\mathcal Z_i}\exp(V_{im})+\exp(V_i^{\mathrm{Off}})]\),
respectively, and satisfy
\(\sum_{j\in\mathcal Z_i}P_{ij}(\bb)+P_i^{\mathrm{Off}}(\bb)=1\). For convenience, set \(P_{ij}(\bb)=0\) whenever
\(j\notin\mathcal Z_i\).

Bonuses may also activate latent drivers, including drivers who may otherwise serve competing platforms.
Let $P_i^{\mathrm{On}}(\bb)$ denote the corresponding bonus-dependent activation probability.
We model this probability using the binary logit model \(P_i^{\mathrm{On}}(\bb)=\sigma(\alpha_i^{\mathrm{On}}+\gamma^{\mathrm{On}}b_i)\) for each \(i\in\mathcal Z\), where $\sigma(z)=1/(1+\exp(-z))$, $\alpha_i^{\mathrm{On}}$ is the zone-specific baseline activation propensity, and $\gamma^{\mathrm{On}}\ge 0$ measures the sensitivity of latent-driver activation to the local bonus.
In deployment, the identifiable behavioral-response parameters \(\bbeta\), \(\boldsymbol{\delta}^{\mathrm{Off}}\), \(\zeta^{\mathrm{Off}}\), \(\boldsymbol{\alpha}^{\mathrm{On}}\), and \(\gamma^{\mathrm{On}}\) can be estimated from historical driver-state transitions using the calibration procedure in \ECsection{EC.1}.
They enter each real-time DSMDP instance as inputs.

Let
$\bP(\bb)=\{P_{ij}(\bb)\}_{i\in\mathcal{Z},\, j\in\mathcal{Z}_i}$
denote the collection of bonus-induced relocation probabilities.
The post-response driver supply in zone $i$ combines the initial dispatchable supply, relocation inflows and outflows, departures from the dispatchable pool, and activation of latent drivers.
Specifically, for each \(i\in\mathcal Z\), the post-response supply is
\begin{equation}
s_i(\bb)
=
s_i^{\mathrm{On}}
+\sum_{j\in\mathcal Z_i\setminus\{i\}}
s_j^{\mathrm{On}}P_{ji}(\bb)
-s_i^{\mathrm{On}}
\left[
\sum_{j\in\mathcal Z_i\setminus\{i\}}P_{ij}(\bb)
+P_i^{\mathrm{Off}}(\bb)
\right]
+s_i^{\mathrm{Off}}P_i^{\mathrm{On}}(\bb).
\label{eq1_supplybalance}
\end{equation}
The five terms in Eq.~\eqref{eq1_supplybalance} represent initial dispatchable supply, relocation inflows, relocation outflows, departures from the dispatchable pool, and activation of latent drivers, respectively.
Because the response probabilities are nonnegative and sum to one, the initial dispatchable supply in each zone is allocated among drivers who remain, relocate, or become unavailable.
The latent-supply term adds the expected number of drivers activated by the local bonus.
These two probability allocations imply that the post-response supply is nonnegative.
We write \(\bs(\bb):=(s_i(\bb))_{i\in\mathcal Z}\) for the post-response supply vector.

Based on the post-response supply, the platform assigns drivers to requests.
A match from driver zone \(i\) to request zone \(j\) has pre-bonus value \(p_{ij}\ge0\) and net contribution \(p_{ij}-b_i\).
The platform limits each zone-level bonus by \(b^{\max}>0\), realized bonus spending by \(B\), and neighboring-zone bonus differences by the spatial smoothness threshold \(\Delta^b\); the latter is consistent with ride-hailing incentive design and observed platform products \citep{ong2021driver}.
Figure~\ref{fig:driver_spatial_incentive_products} illustrates the spatial incentive maps that motivate this smoothness requirement.
The resulting model is
\begin{subequations}
\begin{align}
\text{(DSMDP)}\qquad
\max_{\bb,\bx}\quad
&\sum_{i,j\in\mathcal Z}(p_{ij}-b_i)x_{ij}
\label{eq1_obj}\\[-0.25ex]
\text{s.t.}\quad
&\sum_{i\in\mathcal Z}x_{ij}\le d_j,
&&j\in\mathcal Z,
\label{eq1_demand}\\
&\sum_{j\in\mathcal Z}x_{ij}\le s_i(\bb),
&&i\in\mathcal Z,
\label{eq1_supply}\\
&\sum_{i,j\in\mathcal Z}b_i x_{ij}\le B,
\label{eq1_budget}\\
&|b_i-b_j|\le\Delta^b,
&&i\in\mathcal Z,\quad j\in\mathcal Z_i\setminus\{i\},
\label{eq1_bonusdiff}\\
&x_{ij}\ge0,
&&i,j\in\mathcal Z,
\label{eq1_matchingDomain}\\
&0\le b_i\le b^{\max},
&&i\in\mathcal Z.
\label{eq1_domains}
\end{align}
\end{subequations}
Objective~\eqref{eq1_obj} maximizes the pre-bonus matching value net of bonus payments.
Constraint~\eqref{eq1_demand} limits assignments by rider demand, whereas Constraint~\eqref{eq1_supply} limits them by the supply available after driver response.
Constraint~\eqref{eq1_budget} enforces the realized bonus-spending budget, and Constraint~\eqref{eq1_bonusdiff} imposes spatial smoothness on the displayed map.
Constraints~\eqref{eq1_matchingDomain}--\eqref{eq1_domains} impose nonnegativity on matching flows and bound each zone bonus.
Directly solving DSMDP within the short decision window is challenging because bonuses enter the nonlinear driver-response functions and interact bilinearly with matching flows through the payment terms \(b_i x_{ij}\) in both the objective and the realized-spending constraint.

\section{Marginal Analysis of Zone-Level Bonuses}
\label{sec:matching_value}

Solving the DSMDP requires assessing whether changes in zone-level bonuses generate sufficient matching value to justify the associated payments.
We therefore characterize the marginal value of a zone bonus, identifying when introducing bonuses is beneficial and how spatial smoothness couples bonus decisions across zones.
This characterization also guides the bonus updates in the iterative solution approach developed next.

\subsection{Marginal Value of a Zone Bonus}

To characterize the marginal value of a zone bonus, we first examine how the downstream matching outcome depends on the driver surge map.
For a fixed map \(\bb\), the driver response model determines the available supply \(\bs(\bb)\) through Eq.~\eqref{eq1_supplybalance}.
The DSMDP then reduces to the following LP, which determines the optimal matching outcome conditional on \(\bb\):
\begin{subequations}
\begin{align}
\text{(DSMDP-M)}\qquad
Z(\bb)
=\max_{\bx}\quad
&\sum_{i,j\in\mathcal Z}(p_{ij}-b_i)x_{ij}
\label{eq:matchingObjective}\\
\text{s.t.}\quad
&\sum_{i\in\mathcal Z}x_{ij}\le d_j,
&& j\in\mathcal Z,
\label{eq:matchingDemand}\\
&\sum_{j\in\mathcal Z}x_{ij}\le s_i(\bb),
&& i\in\mathcal Z,
\label{eq:matchingSupply}\\
&\sum_{i,j\in\mathcal Z}b_i x_{ij}\le B
\label{eq:matchingBudget}\\
&x_{ij}\ge0,
&& i,j\in\mathcal Z.
\label{eq:matchingNonnegative}
\end{align}
\end{subequations}

Here, \(Z(\bb)\) denotes the optimal net operating profit under the fixed driver surge map \(\bb\).
For a selected optimal primal-dual solution, let \(\pi_j\ge0\), \(\eta_i\ge0\), and \(\lambda\ge0\) be the optimal dual variables for Constraints~\eqref{eq:matchingDemand}, \eqref{eq:matchingSupply}, and \eqref{eq:matchingBudget}, respectively, representing the shadow prices of rider demand in zone \(j\), post-response driver supply in zone \(i\), and the bonus budget.
The bonus vector \(\bb\) affects the net matching coefficients in Objective~\eqref{eq:matchingObjective}, the post-response supply in Constraint~\eqref{eq:matchingSupply}, and the bonus-payment coefficients in Constraint~\eqref{eq:matchingBudget}.
For the selected optimal primal-dual pair, differentiating the Lagrangian with respect to the zone bonus \(b_m\) gives
\begin{equation}
g_m
:=
\sum_{i\in\mathcal Z}
\eta_i
\frac{\partial s_i(\bb)}{\partial b_m}
-(1+\lambda)\sum_{j\in\mathcal Z}x_{mj},
\qquad \forall m\in\mathcal Z,
\label{eq:grad_phi}
\end{equation}
where $g_m$ is the marginal value of the zone-$m$ bonus associated with the selected optimal primal-dual solution.
We say that the stable optimal-basis condition holds at \(\bb\)
if, in a neighborhood of \(\bb\), DSMDP-M has a unique
primal-dual solution supported by a primal- and dual-nondegenerate
optimal basis, with twice continuously differentiable behavioral-response functions.
Under this condition, the envelope formula gives
\(\bg=\nabla Z(\bb)\).
At points where the optimal basis changes, \(Z(\bb)\) may be nondifferentiable, and different optimal primal-dual pairs may yield different vectors \(\bg\).
Since \(\eta_j\) is the shadow price of supply in zone \(j\), we therefore define \(\bar\eta_i:=\sum_{j\in\mathcal Z_i}P_{ij}(\bb)\eta_j\) as the expected downstream matching value of a dispatchable driver initially in zone \(i\), with zero value assigned to the unavailable state.
Expanding the supply derivatives gives
\(\sum_j\eta_j\frac{\partial s_j}{\partial b_m}
=\sum_i s_i^{\mathrm{On}}\sum_{j\in\mathcal Z_i}
\eta_j\frac{\partial P_{ij}}{\partial b_m}
+\eta_m s_m^{\mathrm{Off}}\frac{\partial P_m^{\mathrm{On}}}{\partial b_m}\).
Let \(\mathbbm{1}_{\{\cdot\}}\) denote the indicator function.
Holding the selected dual vector \(\boldsymbol\eta\) fixed, the MNL response derivatives in \ECsection{EC.2} imply
\begin{equation}
\begin{aligned}
\sum_{j\in\mathcal Z_i}
\eta_j\frac{\partial P_{ij}(\bb)}{\partial b_m}
=
\beta^{(1)}
\mathbbm{1}_{\{m\in\mathcal Z_i\}}
P_{im}(\bb)
\bigl(\eta_m-\bar\eta_i\bigr)
+
\zeta^{\mathrm{Off}}
\mathbbm{1}_{\{i=m\}}
P_i^{\mathrm{Off}}(\bb)\bar\eta_i.
\end{aligned}
\label{eq:weighted_dispatchable_response}
\end{equation}
The first term values attraction to zone \(m\) relative to the current choice-set average, whereas the second values retention of an origin-\(i\) dispatchable driver when the origin-zone bonus makes unavailability less attractive.
Aggregating Eq.~\eqref{eq:weighted_dispatchable_response} across origins and adding latent-driver activation yields the economically interpretable channels below.

For a fixed bonus vector \(\bb\) and a selected optimal primal-dual
solution of DSMDP-M, define the matched quantity
\(q_m:=\sum_{j\in\mathcal Z}x_{mj}\).
The destination-attraction, dispatchable-supply-retention, and
latent-driver-activation terms are defined as
\begin{equation}
\begin{aligned}
\mathsf{Dest}_m
&:=
\beta^{(1)}
\sum_{i:\,m\in\mathcal Z_i}
s_i^{\mathrm{On}}P_{im}(\bb)
\bigl(\eta_m-\bar\eta_i\bigr),\\
\mathsf{Ret}_m
&:=
\zeta^{\mathrm{Off}}
s_m^{\mathrm{On}}P_m^{\mathrm{Off}}(\bb)\bar\eta_m,\\
\mathsf{Act}_m
&:=
\gamma^{\mathrm{On}}s_m^{\mathrm{Off}}
P_m^{\mathrm{On}}(\bb)
\bigl(1-P_m^{\mathrm{On}}(\bb)\bigr)\eta_m.
\end{aligned}
\label{eq:matching_value_channels}
\end{equation}

\begin{proposition}[Decomposition of marginal bonus value]
\label{prop:matching_weighted_bonus_value}
For every selected optimal primal-dual solution of DSMDP-M,
the marginal bonus value decomposes as
\begin{equation}
g_m
=
\mathsf{Dest}_m
+
\mathsf{Ret}_m
+
\mathsf{Act}_m
-
(1+\lambda)q_m,
\qquad \forall m\in\mathcal Z.
\label{eq:matching_weighted_marginal}
\end{equation}
\end{proposition}

Proposition~\ref{prop:matching_weighted_bonus_value} separates the
matching value generated by driver responses from the additional
payments on existing matches and budget opportunity cost.
Destination attraction $\mathsf{Dest}_m$ compares the supply shadow value in zone
\(m\) with the driver's expected downstream matching value across
existing choices.
For \(\beta^{(1)}>0\), the contribution from origin \(i\) to zone
\(m\) is positive precisely when \(\eta_m>\bar\eta_i\), provided
\(s_i^{\mathrm{On}}P_{im}(\bb)>0\).
Dispatchable-supply retention $\mathsf{Ret}_m$ preserves the expected matching
value \(\bar\eta_m\) of drivers who would otherwise become
unavailable, whereas latent-driver activation $\mathsf{Act}_m$ adds supply
directly to zone \(m\), where it is valued at \(\eta_m\).
A higher bonus also increases payments on existing matches by
\(q_m\) per unit increase, with an additional budget opportunity
cost of \(\lambda q_m\).
The value of a bonus increase therefore depends on both driver
responsiveness and downstream matching opportunities, which
local shortage alone does not capture.

If the selected optimal primal-dual solution is supported by a
locally valid optimal basis, \(g_m\) is the sensitivity associated
with that basis.
Under the stable optimal-basis condition, this sensitivity equals
the derivative \(\partial Z(\bb)/\partial b_m\).

\subsection{Local Value of Bonus Intervention}

Proposition~\ref{prop:matching_weighted_bonus_value}
characterizes how a marginal change in a zone bonus affects
net operating profit.
We next use this decomposition to examine when the platform
should begin offering bonuses and how this decision depends
on the supply available without incentives and drivers'
responsiveness.
We use a superscript \(\mathrm{no}\) for quantities evaluated at
the no-bonus map \(\boldsymbol 0\), with
\(P_i^{\mathrm{On},\mathrm{no}}
:=\sigma(\alpha_i^{\mathrm{On}})\).
For each zone \(i\in\mathcal Z\), let \(\boldsymbol e_i\) have
entry 1 at \(i\) and 0 elsewhere, so that
\(b_i\boldsymbol e_i\) offers bonus \(b_i\) only in zone \(i\).
Taking the no-bonus map as the benchmark, the following
corollary identifies zones where introducing a small positive
bonus increases net operating profit.

\begin{corollary}[Local value of the first bonus]
\label{cor:no_bonus_intervention}
Suppose \(B>0\), \(b^{\max}>0\), and \(\Delta^b>0\),
and suppose that \(Z(\boldsymbol b)\) is differentiable
at \(\boldsymbol b=\boldsymbol 0\).
The one-sided marginal value of introducing a bonus is
\begin{equation}
\left.
\frac{\mathrm d}{\mathrm d b_i}
Z(b_i\boldsymbol e_i)
\right|_{b_i=0+}
=
\mathsf{Dest}_i^{\mathrm{no}}
+\mathsf{Ret}_i^{\mathrm{no}}
+\gamma^{\mathrm{On}}s_i^{\mathrm{Off}}
P_i^{\mathrm{On},\mathrm{no}}
\bigl(1-P_i^{\mathrm{On},\mathrm{no}}\bigr)\eta_i^{\mathrm{no}}
-q_i^{\mathrm{no}}.
\label{eq:no_bonus_intervention_threshold}
\end{equation}
\end{corollary}

If the right-hand side of
Eq.~\eqref{eq:no_bonus_intervention_threshold}
is positive for some zone \(i\), then a sufficiently small
positive bonus in that zone improves the platform objective,
and the no-bonus map is neither locally nor globally optimal.
The no-bonus map satisfies the first-order stationarity
condition if the right-hand side is nonpositive for every
zone and is a strict local maximizer if it is strictly
negative for every zone.
To make the intervention condition in
Corollary~\ref{cor:no_bonus_intervention} explicit in terms
of driver composition and activation responsiveness,
consider a single-zone setting with no relocation
or dispatchable-supply exit.
Suppose demand is sufficient to match all available supply
and the bonus budget remains slack for bonuses sufficiently
close to zero.
These assumptions isolate the trade-off between activating
additional drivers and paying bonuses on existing matches.
Each additional activated driver generates a match with
pre-bonus value \(p_{ii}>0\), and the platform's net operating
profit is
\begin{equation}
Z_i(b_i)
=
(p_{ii}-b_i)
\left[
s_i^{\mathrm{On}}
+s_i^{\mathrm{Off}}P_i^{\mathrm{On}}(b_i)
\right].
\label{eq:one_zone_objective}
\end{equation}
For \(s_i^{\mathrm{Off}}>0\) and \(\gamma^{\mathrm{On}}>0\),
differentiating Eq.~\eqref{eq:one_zone_objective} at zero shows
that the initial marginal bonus value is positive if and only if
\(p_{ii}>p_i^{\mathrm{crit}}\), where
\begin{equation}
p_i^{\mathrm{crit}}
:=
\frac{
s_i^{\mathrm{On}}/s_i^{\mathrm{Off}}
+P_i^{\mathrm{On},\mathrm{no}}
}{
\gamma^{\mathrm{On}}
P_i^{\mathrm{On},\mathrm{no}}
\bigl(1-P_i^{\mathrm{On},\mathrm{no}}\bigr)
}.
\label{eq:explicit_activation_threshold}
\end{equation}
The factor
\(P_i^{\mathrm{On},\mathrm{no}}(1-P_i^{\mathrm{On},\mathrm{no}})\)
is maximized at \(P_i^{\mathrm{On},\mathrm{no}}=1/2\).
Thus, holding \(\gamma^{\mathrm{On}}\) fixed, a small bonus
has the greatest effect on activation probability when
latent drivers are equally likely to participate or remain
inactive without incentives.

\subsection{Spatial Balance of Marginal Values}

Once bonuses are introduced, spatial smoothness limits how the
platform can adjust them across zones. If a zone's bonus already
exceeds a neighboring zone's bonus by \(\Delta^b\), increasing it further requires raising the neighboring bonus as well. We therefore characterize how marginal bonus values balance across zones whose bonus adjustments are linked by spatial smoothness.

Let \(
\mathcal E^b
=
\bigl\{\{i,j\}:i\in\mathcal Z,\;
j\in\mathcal Z_i\setminus\{i\}\bigr\}.
\)
Each unordered neighboring-zone pair appears once in \(\mathcal E^b\).
For each edge \(e=\{i,j\}\in\mathcal E^b\), choose an arbitrary orientation \((i,j)\).
Let \(\boldsymbol D\) be the oriented incidence matrix, with one
row per zone pair and one column per zone, and let
\(\boldsymbol 1\) denote an all-ones vector of conformable dimension.
For each edge \(e\) oriented as \((i,j)\), set
\(D_{ei}=1\), \(D_{ej}=-1\), and all other entries in that
row to zero.
Thus, \((\boldsymbol D\bb)_e=b_i-b_j\).
For \(0<\Delta^b<+\infty\), the bonus-feasible set is
\begin{equation}
\mathcal B
:=
\left\{
\bb\in\mathbb R^{|\mathcal Z|}:
\boldsymbol 0\le\bb\le b^{\max}\boldsymbol 1,\;
-\Delta^b\boldsymbol 1
\le\boldsymbol D\bb
\le\Delta^b\boldsymbol 1
\right\}.
\label{eq:bonus_feasible_set}
\end{equation}

Consider a locally optimal bonus map \(\bb^\star\in\mathcal B\),
i.e., there exists a radius \(r_{\mathrm{loc}}>0\) such that
\(Z(\bb^\star)\ge Z(\bb)\) for every \(\bb\in\mathcal B\)
satisfying \(\|\bb-\bb^\star\|_2<r_{\mathrm{loc}}\).
Assume \(Z\) is differentiable at \(\bb^\star\), where
\(Z(\bb^\star)\) denotes the optimal value of DSMDP-M
under bonus map \(\bb^\star\).
Let \(\bg^\star:=\nabla Z(\bb^\star)\), with entries
\(g_i^\star\) given by
Eq.~\eqref{eq:matching_weighted_marginal}.
We call zone \(i\) interior if \(0<b_i^\star<b^{\max}\).
For the constraints defining \(\mathcal B\), let
\(u_i^+\) and \(u_i^-\) denote the nonnegative Karush--Kuhn--Tucker (KKT) multipliers
for the upper and lower bonus bounds, respectively.
For each edge \(e\), oriented as \((i,j)\), let
\(\rho_e^+\) and \(\rho_e^-\) denote the nonnegative KKT
multipliers for \(b_i-b_j\le\Delta^b\) and
\(b_j-b_i\le\Delta^b\), respectively.
These multipliers satisfy complementary slackness,
while KKT stationarity yields the following
balance of marginal bonus values.

\begin{proposition}[Spatial balance of marginal bonus values]
\label{prop:spatial_marginal_balance}
At \(\bb^\star\), the first-order optimality condition can be
written as
\begin{equation}
g_i^\star
=
u_i^+-u_i^-
+
\sum_{e\in\mathcal E^b}
D_{ei}(\rho_e^+-\rho_e^-),
\qquad i\in\mathcal Z.
\label{eq:spatial_kkt_balance}
\end{equation}
If zone \(i\) is interior and all spatial smoothness multipliers
on edges incident to \(i\) are zero, then \(g_i^\star=0\).
Substituting this into Eq.~\eqref{eq:matching_weighted_marginal} yields
\begin{equation}
\mathsf{Dest}_i^\star
+\mathsf{Ret}_i^\star
+\mathsf{Act}_i^\star
=
(1+\lambda^\star)q_i^\star.
\label{eq:interior_zone_balance}
\end{equation}
In particular, this condition holds when zone \(i\) is interior
and all of its incident spatial smoothness constraints are slack.
\end{proposition}

Proposition~\ref{prop:spatial_marginal_balance} shows that
spatial smoothness couples marginal bonus values across neighboring
zones, so an interior zone need not satisfy \(g_i^\star=0\).
Define \(\mathcal E_\rho:=\{e\in\mathcal E^b:
\rho_e^++\rho_e^->0\}\) and consider the graph
\((\mathcal Z,\mathcal E_\rho)\).
The following corollary characterizes how these marginal values
balance over each connected component, i.e., a maximal group
of zones linked, directly or through other zones, by edges
in \(\mathcal E_\rho\).

\begin{corollary}[Balance across connected zones]
\label{cor:spatial_component_balance}
For any multiplier system satisfying
Proposition~\ref{prop:spatial_marginal_balance}, let \(C\) be a
connected component of \((\mathcal Z,\mathcal E_\rho)\). Then
\begin{equation}
\sum_{i\in C}g_i^\star
=
\sum_{i\in C}(u_i^+-u_i^-).
\label{eq:component_balance_with_bounds}
\end{equation}
If every zone in \(C\) is interior, then
\begin{equation}
\sum_{i\in C}
\left[
\mathsf{Dest}_i^\star
+\mathsf{Ret}_i^\star
+\mathsf{Act}_i^\star
-(1+\lambda^\star)q_i^\star
\right]
=0.
\label{eq:active_component_balance}
\end{equation}
\end{corollary}

Corollary~\ref{cor:spatial_component_balance} further shows that spatial smoothness links marginal bonus values across zones, so bonus adjustments must account for the active smoothness constraints that couple neighboring zones.
If every zone in a component is interior, positive and negative
marginal bonus values offset one another.
When \(\Delta^b=+\infty\), spatial smoothness constraints are
absent, and each interior zone balances independently.
This constrained balance motivates ISA, which combines marginal
sensitivities with active bonus-map constraints when constructing its update direction.

\section{Solution Approach and Analytical Results}
\label{theories}

The analysis in Section~\ref{sec:matching_value} provides
marginal bonus values that can guide adjustments to the
bonus map. We use this information to develop ISA, an
iterative solution approach for the DSMDP. At each iteration,
ISA evaluates driver responses and solves DSMDP-M, then
uses the resulting sensitivity information to update bonuses
subject to the bonus bounds and spatial smoothness constraints.
We first present the algorithm, then analyze its stationarity
properties and develop an upper bound on the optimal
net operating profit.

\subsection{Iterative Solution Approach}
\label{solution_approach}

Let \(\bb^k\in\mathcal B\) denote the current bonus map at
iteration \(k\), where \(\mathcal B\) is defined in
Eq.~\eqref{eq:bonus_feasible_set}.
As bonuses vary, the optimal matching basis may change,
making \(Z(\bb)\) nonsmooth.
We refer to the objective-value function defined by a fixed optimal matching basis as a matching branch.
An adjustment that improves one branch may not improve
another, so ISA uses sensitivities from nearby branches
to construct its update.

At the current map \(\bb^k\), referred to as the center,
ISA samples nearby bonus maps within a sampling radius
\(\epsilon_k>0\).
These ambient samples are drawn from the surrounding bonus
space without being restricted to \(\mathcal B\).
Let \(\mathcal A_k\) index the ambient samples, and let
\(\by^{k,a}\), \(a\in\mathcal A_k\), denote the sampled bonus maps.
Set \(\by^{k,0}:=\bb^k\) and
\(\mathcal I_k:=\{0\}\cup\mathcal A_k\).
At each \(\by^{k,a}\), \(a\in\mathcal I_k\), ISA evaluates
driver responses, solves DSMDP-M, and uses the optimal
primal-dual solution to compute the corresponding branch
sensitivity \(\bg^{k,a}\) from
Eq.~\eqref{eq:matching_weighted_marginal}.
For the center sensitivity, write \(\bg^k:=\bg^{k,0}\).
The matrix \(\bG_k=[\bg^{k,a}]_{a\in\mathcal I_k}\)
collects these sensitivities and forms the working bundle.
For later stationarity assessment, we use the ambient sensitivity
matrix \(\bG_k^{\mathrm{amb}}:=[\bg^{k,a}]_{a\in\mathcal A_k}\),
which excludes the center evaluation.

To construct a direction that respects the bonus constraints, let
\(\bN_k:=[\bn^{k,c}]_{c\in\mathcal C_k}\)
collect the outward normals of all constraints defining
\(\mathcal B\) that are active at \(\bb^k\),
where \(\mathcal C_k\) indexes these constraints and
\(\bn^{k,c}\) is the outward normal of constraint \(c\).
Let \(\bomega\) denote the convex-combination weights on sampled
sensitivities, with \(\bomega\ge0\) and
\(\boldsymbol 1^\top\bomega=1\), and let \(\bmu\ge0\) denote the
weights on the columns of \(\bN_k\).
ISA selects a convex combination \(\bG_k\bomega\) of sampled
sensitivities and a nonnegative combination \(\bN_k\bmu\) of active
constraint normals by solving the quadratic program (QP)
\begin{equation}
(\bomega^k,\bmu^k)
\in
\arg\min_{\substack{
\bomega\ge\boldsymbol{0},\;
\boldsymbol{1}^\top\bomega=1\\
\bmu\ge\boldsymbol{0}
}}
\frac12
\left\|
\bG_k\bomega-\bN_k\bmu
\right\|_2^2.
\label{eq:sampled_stationarity_qp}
\end{equation}
The QP finds the closest pair of these combinations.
Define
\(\br^k:=\bG_k\bomega^k-\bN_k\bmu^k\)
as the constraint-adjusted direction at iteration \(k\) and
\(\widehat R_k:=\|\br^k\|_2\)
as the working-bundle residual.
When \(\br^k\neq\boldsymbol{0}\), the QP optimality conditions give
\(\bN_k^\top\br^k\le\boldsymbol{0}\) and
\((\bg^{k,a})^\top\br^k\ge\|\br^k\|_2^2>0\)
for every \(a\in\mathcal I_k\).
Thus, the direction respects the active constraints and predicts a
positive improvement for every sampled branch.
The step size is then determined by evaluating the actual
change in net operating profit.

For a candidate step size \(\tau>0\), ISA projects
\(\bb^k+\tau\br^k\) onto \(\mathcal B\).
Let \(\Pi_{\mathcal B}\) denote this Euclidean projection
and define the feasible trial displacement as
$
\bh_k(\tau)
:=
\Pi_{\mathcal B}(\bb^k+\tau\br^k)-\bb^k.
$
The trial bonus map is therefore
\(\bb^k+\bh_k(\tau)\).
The projection enforces the bonus bounds and spatial
smoothness constraints.
ISA evaluates driver responses under the trial map and
solves DSMDP-M to obtain its net operating profit while
enforcing the realized bonus-spending budget.

To obtain sufficient improvement while avoiding a time-consuming exact line search, we select the step size using Armijo backtracking.
Let \(\chi\in(0,1)\) denote the Armijo sufficient-increase
parameter.
The projected bundle Armijo test accepts the trial when
\(\min_{a\in\mathcal I_k}
(\bg^{k,a})^\top\bh_k(\tau)>0\) and
\begin{equation}
Z\bigl(\bb^k+\bh_k(\tau)\bigr)
\ge
Z(\bb^k)
+
\chi
\min_{a\in\mathcal I_k}
(\bg^{k,a})^\top\bh_k(\tau).
\label{eq:sampled_armijo_main}
\end{equation}
Thus, all sampled sensitivities must predict an improvement,
and the realized improvement must be at least a fraction
\(\chi\) of the smallest predicted improvement.
ISA reduces the step size when a trial fails this test.
Otherwise, the trial becomes the new center.
If no trial is accepted within the line-search limits,
ISA retains the current center and adds new samples,
producing a null step.

To determine whether to terminate, ISA first checks the working-bundle residual \(\widehat R_k\) against a stationarity scale \(\nu_k>0\).
A small \(\widehat R_k\) indicates that the convex combination
\(\bG_k\bomega^k\) of sampled sensitivities nearly matches
the nonnegative combination \(\bN_k\bmu^k\) of active
constraint normals, i.e.,
\(\|\bG_k\bomega^k-\bN_k\bmu^k\|_2\approx0\).
However, this balance may be affected by the center sensitivity \(\bg^{k,0}\),
as \(Z\) can be nonsmooth at \(\bb^k\).
ISA therefore performs an ambient-only stationarity assessment
by replacing \(\bG_k\) with \(\bG_k^{\mathrm{amb}}\) in
Eq.~\eqref{eq:sampled_stationarity_qp}.
The resulting ambient residual \(\widehat R_k^{\mathrm{amb}}\),
defined in Eq.~\eqref{eq:ambient_stationarity_residual},
is used to determine whether to refine the sampling and stationarity scales or terminate.
Complete sampling, line-search, and stopping procedures are
provided in \ECsection{EC.3}.
Algorithm~\ref{alg:AO} summarizes ISA, which returns the best
feasible incumbent.

\begin{algorithm}[htbp]
\footnotesize
\renewcommand{\algorithmicrequire}{\textbf{Input:}}
\renewcommand{\algorithmicensure}{\textbf{Output:}}
\caption{Iterative solution approach for DSMDP.}
\label{alg:AO}
\begin{algorithmic}[1]
\Require Initial map \(\bb^0\in\mathcal B\), sampling scales, line-search parameters, and computation limits
\Ensure Final bonuses \(\bb^{\mathrm{best}}\) and matching decision \(\bx^{\mathrm{best}}\)
\State Initialize \(k\leftarrow0\), the sample pool, and \(Z^{\mathrm{best}}\leftarrow-\infty\)
\While{the computation limits remain}
    \State \textbf{(Center evaluation)} Compute driver response and supply at \(\bb^k\), and solve DSMDP-M in Eqs.~\eqref{eq:matchingObjective}--\eqref{eq:matchingNonnegative} for \(Z(\bb^k)\), \(\bx^k\), the dual multipliers, and \(\bg^k\)
    \State \textbf{if} \(Z(\bb^k)>Z^{\mathrm{best}}\) \textbf{then} set \((Z^{\mathrm{best}},\bb^{\mathrm{best}},\bx^{\mathrm{best}})\leftarrow(Z(\bb^k),\bb^k,\bx^k)\)
    \State \textbf{(Sensitivity sampling)} Evaluate nearby matching branches and form \(\bG_k\) and \(\bN_k\)
    \State \textbf{(Constrained direction)} Solve Eq.~\eqref{eq:sampled_stationarity_qp} for \(\br^k\) and \(\widehat R_k\)
    \State Perform the ambient-only assessment when required (see \ECsection{EC.3.3}).
    \State \textbf{(Line search)} Apply Eq.~\eqref{eq:sampled_armijo_main} with \(\bh_k(\tau)=\Pi_{\mathcal B}(\bb^k+\tau\br^k)-\bb^k\)
    \If{the trial is accepted}
        \State Set \(\bb^{k+1}\leftarrow\bb^k+\bh_k(\tau)\), update the best incumbent, and set \(k\leftarrow k+1\)
        \State Rebuild the sample pool around the new center
    \Else
        \State Retain the center and \(k\), and enrich the sample pool
    \EndIf
\EndWhile
\State \textbf{return} the best feasible incumbent \((\bb^{\mathrm{best}},\bx^{\mathrm{best}})\)
\end{algorithmic}
\end{algorithm}

\subsection{Theoretical Analysis}

We next examine the theoretical properties of the proposed solution approach.
We first establish the existence of an optimal solution and show how the sampling radius and ambient residual provide a finite-scale stationarity certificate.
We then study stationarity along ambient assessments
with vanishing sampling radii and residuals.
Proofs are provided in \ECsection{EC.4}.

\subsubsection{Existence and Stationarity.}
\label{subsec:theoretical_analysis}

For each fixed bonus vector, DSMDP-M is solved to optimality, and
\(Z(\bb)\) denotes its optimal value. The optimization problem and
Algorithm~\ref{alg:AO} use a positive budget \(B\) and solve
\(\max_{\bb\in\mathcal B}Z(\bb)\).
The set \(\mathcal B\) is nonempty, compact, and convex because it is
defined by bonus bounds and linear spatial smoothness constraints.

\begin{lemma}[Continuity and existence]
\label{lemma:continuity_existence}
Suppose that \(P_{ij}(\bb)\), \(P_i^{\mathrm{On}}(\bb)\), and
\(P_i^{\mathrm{Off}}(\bb)\) are continuous on an open set containing
\(\mathcal B\).
Then \(Z(\bb)\) is continuous on \(\mathcal B\), and the reduced
DSMDP \(\max_{\bb\in\mathcal B}Z(\bb)\) admits a globally optimal
solution.
\end{lemma}

Lemma~\ref{lemma:continuity_existence} shows that small
changes in bonuses lead to small changes in optimal net
operating profit, even when the optimal matching decisions
change.
We next examine how ISA can assess stationarity using
the sensitivities obtained from nearby bonus maps.

For the finite-scale analysis, suppose that \(Z\) admits a locally
Lipschitz extension around \(\mathcal B\).
For \(\epsilon>0\), define the enlarged Clarke generalized gradient
\begin{equation}
\partial_\epsilon^C Z(\bb)
:=
\operatorname{cl\,conv}
\left(
\bigcup_{\|\by-\bb\|_2\le\epsilon}
\partial^C Z(\by)
\right),
\label{eq:epsilon_clarke_gradient}
\end{equation}
where \(\operatorname{cl\,conv}\) denotes the closed convex hull,
the union is taken over points in the domain of the extension,
and \(\partial^C Z(\by)\) denotes its Clarke generalized gradient.
For two nonempty sets \(A_1\) and \(A_2\), let
\(
\operatorname{dist}(A_1,A_2)
:=
\inf\{\|\boldsymbol a_1-\boldsymbol a_2\|_2:
\boldsymbol a_1\in A_1,\boldsymbol a_2\in A_2\}.
\)
Let \(\mathcal N_{\mathcal B}(\bb)\) denote the outward normal cone
\(\{\bn\in\mathbb R^{|\mathcal Z|}:
\bn^\top(\widetilde{\bb}-\bb)\le0,\;
\forall\widetilde{\bb}\in\mathcal B\}\).
In particular, at the current map \(\bb^k\),
\(\mathcal N_{\mathcal B}(\bb^k)
=\{\bN_k\bmu:\bmu\ge\boldsymbol{0}\}\).
For \(\epsilon>0\) and \(\nu\ge0\), a feasible map \(\bb\) is
\((\epsilon,\nu)\)-stationary if
\(\operatorname{dist}(\partial_\epsilon^C Z(\bb),\mathcal N_{\mathcal B}(\bb))\le\nu\).
We next show that a sufficiently small ambient residual certifies
approximate stationarity of the bonus map at the given sampling
radius.

\begin{proposition}[Finite-scale stationarity certificate]
\label{prop:finite_scale_certificate}
At a center \(\bb^k\in\mathcal B\), let
\(\mathcal A_k\ne\varnothing\), and suppose the ambient samples
\(\by^{k,a}\) are differentiability points of the extension satisfying
\(\|\by^{k,a}-\bb^k\|_2\le\epsilon_k\) and
\(\bg^{k,a}=\nabla Z(\by^{k,a})\),
\(\forall a\in\mathcal A_k\).
Let \(\bG_k^{\mathrm{amb}}\) be the ambient sensitivity matrix, and let the columns of \(\bN_k\) be exactly the active outward normals of \(\mathcal B\) at \(\bb^k\), up to duplicate generators.
Define the ambient residual as the QP residual obtained from the ambient sensitivity matrix
\begin{equation}
\widehat R_k^{\mathrm{amb}}
:=
\min_{\substack{
\bomega\ge\boldsymbol 0,\;\boldsymbol 1^\top\bomega=1\\
\bmu\ge\boldsymbol 0
}}
\left\|
\bG_k^{\mathrm{amb}}\bomega-\bN_k\bmu
\right\|_2.
\label{eq:ambient_stationarity_residual}
\end{equation}
Then
\begin{equation}
\operatorname{dist}
\left(
\partial_{\epsilon_k}^C Z(\bb^k),
\mathcal N_{\mathcal B}(\bb^k)
\right)
\le
\widehat R_k^{\mathrm{amb}}.
\label{eq:finite_scale_certificate}
\end{equation}
Consequently, for any \(\nu_k\ge0\), \(\widehat R_k^{\mathrm{amb}}\le\nu_k\) certifies \((\epsilon_k,\nu_k)\)-stationarity of the current map.
\end{proposition}

Proposition~\ref{prop:finite_scale_certificate} shows that the ambient residual \(\widehat R_k^{\mathrm{amb}}\) provides
a certificate of approximate stationarity.
The sampling radius \(\epsilon_k\) determines how far the sampled bonus maps may lie from the current map \(\bb^k\),
while the stationarity scale \(\nu_k\) sets the acceptable upper bound on the ambient residual \(\widehat R_k^{\mathrm{amb}}\).
The result applies to the ambient-only assessment because
every sampled sensitivity is a gradient at a nearby
differentiability point and \(\bN_k\) contains all
active-constraint normals.
Since the working bundle also includes the center sensitivity,
which may not be a gradient at a nonsmooth matching transition,
a small working-bundle residual triggers this separate assessment.
We next examine the limiting behavior of the ambient-only
assessments as the sampling radius and residual vanish.

\begin{proposition}[Clarke stationarity of sampled limits]
\label{prop:sampled_clarke_limit}
Consider a sequence of ambient-only assessments satisfying
the conditions of
Proposition~\ref{prop:finite_scale_certificate}.
If
\(
\bb^\ell\to\bar{\bb},
\epsilon_\ell\downarrow0,
\widehat R_\ell^{\mathrm{amb}}\to0,
\)
then
\[
0\in-\partial^C Z(\bar{\bb})
+\mathcal N_{\mathcal B}(\bar{\bb}).
\]
Hence, \(\bar{\bb}\) is a constrained Clarke-stationary
point of DSMDP.
\end{proposition}

Proposition~\ref{prop:sampled_clarke_limit} shows that
a convergent sequence of bonus maps from ambient-only
assessments has a constrained Clarke-stationary limit
when the sampling radius and ambient residual vanish.
A smaller sampling radius restricts the assessment to
bonus maps closer to the current map, while a smaller
ambient residual provides a tighter stationarity
certificate at that sampling radius.
With Proposition~\ref{prop:finite_scale_certificate},
this result supports ISA's scale-refinement and stopping rules:
the ambient residual test
\(\widehat R_k^{\mathrm{amb}}\le\nu_k\)
certifies \((\epsilon_k,\nu_k)\)-stationarity,
and scale refinement reduces the sampling radius
\(\epsilon_k\) and stationarity scale \(\nu_k\).
Termination upon passing this test at the target sampling radius
\(\epsilon_{\mathrm{tar}}>0\) and stationarity scale
\(\nu_{\mathrm{tar}}>0\) provides an
\((\epsilon_{\mathrm{tar}},\nu_{\mathrm{tar}})\)-stationarity
certificate for the bonus map.

\subsection{Conic Outer Relaxation and Upper Bound}
\label{subsec:conic_upper_bound}

To assess the quality of the feasible solutions produced
by ISA, we construct a conic outer relaxation
\(\overline{\mathrm{DSMDP}}\) that provides an upper
bound on the optimal net operating profit.
The relaxation uses exponential-cone constraints and
McCormick envelopes to relax the driver-response
relations and bonus-payment products. 

For the MNL model, introduce variables \(t_{ij}\) and
\(t_i^{\mathrm{Off}}\) to represent
\(\exp(V_{ij})\) and \(\exp(V_i^{\mathrm{Off}})\),
respectively, for \(i\in\mathcal Z\) and
\(j\in\mathcal Z_i\).
Define the denominator variable
\(D_i=\sum_{m\in\mathcal Z_i}t_{im}+t_i^{\mathrm{Off}}\).
Let \(\underline V_{ij}\) and \(\overline V_{ij}\)
denote valid lower and upper bounds on \(V_{ij}\)
over \(\mathcal B\), and let
\(\mathcal S_{\exp}(v;\ell,u)\) denote the affine
secant of the exponential function over \([\ell,u]\),
with \(\mathcal S_{\exp}(v;\ell,\ell)=\exp(\ell)\).
For each \(i\in\mathcal Z\) and \(j\in\mathcal Z_i\),
we relax \(t_{ij}=\exp(V_{ij})\) using
\begin{equation}
\exp(V_{ij})
\le t_{ij}
\le
\mathcal S_{\exp}
\bigl(V_{ij};\underline V_{ij},\overline V_{ij}\bigr).
\label{eq:conic_exp_summary}
\end{equation}
The lower inequality is representable using the
exponential cone, while the upper inequality follows
from convexity of the exponential function.
The same construction applies to \(t_i^{\mathrm{Off}}\).

We treat \(P_{ij}\) and \(P_i^{\mathrm{Off}}\) as decision
variables and relax the bilinear equalities
\(t_{ij}=P_{ij}D_i\) and
\(t_i^{\mathrm{Off}}=P_i^{\mathrm{Off}}D_i\)
using McCormick envelopes over valid variable bounds.
The probability normalization
\(\sum_{j\in\mathcal Z_i}P_{ij}+P_i^{\mathrm{Off}}=1\)
and nonnegativity constraints are retained.
For latent-driver activation, use the same procedure
to bound the exponential term in the binary logit denominator
and use a McCormick envelope to relax the equality
between the probability--denominator product and one.

For bonus payments, recall the matched quantity
\(q_i=\sum_{j\in\mathcal Z}x_{ij}\), and introduce
\(w_i\) to represent \(b_iq_i\).
Let \(\bar q_i\) be a valid upper bound on \(q_i\),
obtained from total demand and the maximum possible
post-response supply in zone \(i\).
Over \((b_i,q_i)\in[0,b^{\max}]\times[0,\bar q_i]\),
the McCormick envelope is
\begin{equation}
\begin{aligned}
w_i &\ge 0,
&
w_i &\ge
b^{\max}q_i+\bar q_i b_i-b^{\max}\bar q_i,\\
w_i &\le b^{\max}q_i,
&
w_i &\le \bar q_i b_i,
\end{aligned}
\qquad i\in\mathcal Z.
\label{eq:conic_payment_summary}
\end{equation}
The relaxed objective is
\(\sum_{i,j\in\mathcal Z}p_{ij}x_{ij}
-\sum_{i\in\mathcal Z}w_i\),
and the budget constraint becomes
\(\sum_{i\in\mathcal Z}w_i\le B\).
The post-response supply constraints are linear in the probability variables, while the bonus bounds and spatial smoothness constraints remain
\(\bb\in\mathcal B\).
The variable bounds and complete conic formulation
are provided in
\textcolor{msomecblue}{Sections~EC.5.1--EC.5.5 of the
e-companion}.

Let \(Z^\star\) and \(\overline Z\) denote the optimal
objective values of DSMDP and
\(\overline{\mathrm{DSMDP}}\), respectively.
Every feasible DSMDP solution satisfies the exponential
bounds and McCormick envelopes when the auxiliary
variables equal the quantities they represent.
It therefore has a feasible representation in the
conic relaxation with the same objective value,
giving \(Z^\star\le\overline Z\), as proved in
\ECsection{EC.5.6}.
If the relaxation is not solved to optimality, the
solver-certified upper bound \(U^{\mathrm{solver}}\)
remains valid, so that
\(Z^\star\le\overline Z\le U^{\mathrm{solver}}\).
The difference between \(U^{\mathrm{solver}}\) and
the net operating profit of ISA's best feasible incumbent
therefore bounds its optimality gap.
To obtain a feasible Conic benchmark, we project the
relaxation incumbent's bonus map onto \(\mathcal B\),
recompute driver responses, and solve DSMDP-M, as
described in \ECsection{EC.5.7}.

\section{Computational Study}
\label{exps}

We use generated instances and a New York City case study
to address five questions.
First, we examine whether offering zone-level bonuses
improves the platform's net operating profit by comparing
the resulting solutions with the no-bonus benchmark
\(\mathbf B^0\).
Second, we investigate how the bonus budget affects
realized spending and net operating profit, thereby
assessing the value of allocating additional funds
to driver incentives.
Third, we examine how driver responsiveness,
behavioral-parameter estimation error, and spatial
smoothness affect the benefits of bonus intervention.
Fourth, we evaluate whether ISA achieves higher net
operating profit than the feasible benchmark methods
on instances of operational scale under the specified
timing protocol.
Finally, we use the New York City case study to examine
ISA's large-scale performance and the spatial allocation
of bonuses, using data-derived geography, demand,
and a proxy for dispatchable supply.
These analyses also examine how the marginal bonus values
and spatial balance characterized in
Section~\ref{sec:matching_value} help explain where bonuses
are allocated and when they improve net operating profit.

\subsection{Instance Generation}
\label{exp_gen}

To evaluate different system configurations under controlled conditions, we conduct experiments on generated ride-hailing instances.
Each instance represents a single real-time decision snapshot of a spatial ride-hailing market.
Following spatial ride-hailing models \citep{bimpikis2019spatial}, we construct each market on a \(\mathrm{rows}\times\mathrm{columns}\) grid, where each cell represents a service zone and \(d_{ij}\) is the Euclidean distance between zone centers. Rider demand is generated from a base level and nonoverlapping spatial hotspots following \citet{soza2024shareability}, with hotspot intensity and dispersion controlled by \(\theta^d\) and \(\sigma^d\). Initial dispatchable driver supply is distributed around activity centers and locally shifted toward demand-concentrated areas, producing spatial supply--demand mismatches \citep{ong2021driver}. We control the aggregate market condition through \(\kappa=\sum_{i\in\mathcal Z}s_i^{\mathrm{On}}/\sum_{j\in\mathcal Z}d_j\).
Latent supply is generated to be correlated with demand, while pre-bonus matching values combine pickup-distance decay and destination attractiveness. For each zone \(i\), \(\mathcal Z_i\) denotes the
\(3\times3\) neighborhood used to define short-term relocation options and spatial smoothness constraints.
Detailed construction procedures and parameter settings are provided in \ECsection{EC.6.1}.

\subsection{Benchmark Methods}
\label{exp_benchmarks}

We compare ISA with five feasible benchmarks motivated
by ride-hailing incentives, matching optimization,
and nonlinear optimization
\citep{lu2018surge,bimpikis2019spatial,besbes2021surge,xu2018large,ozkan2020dynamic}.
1) \(\mathbf B^0\): The platform posts no bonuses.
2) \(\mathbf H_1\): This heuristic constructs a one-pass bonus map from local shortage, matching value, and MNL responsiveness.
3) \(\mathbf H_2\): This heuristic constructs a one-pass bonus map from the dual sensitivities of the no-bonus matching LP.
4) NOMAD: We apply derivative-free local search to the bonus vector \citep{le2011nomad}.
5) Conic: We project the bonus vector obtained from the conic outer relaxation in
Section~\ref{subsec:conic_upper_bound} onto
\(\mathcal B\), recompute driver responses, and solve DSMDP-M to obtain a feasible solution.
We also report UB, the solver-certified upper bound for the conic outer relaxation.
In the result tables, Obj. denotes the objective value and Time reports wall-clock seconds.
The green percentages beside ISA's objective values report its improvement over the strongest non-ISA feasible benchmark.
All method entries are feasible; UB is excluded from method comparisons.
Detailed benchmark procedures are provided in \ECsection{EC.6}.

All experiments are run on a Windows server with AMD EPYC 9654 processors running at 2.4~GHz, 192 physical cores, and approximately 1.5~TiB of memory. Further experimental details are provided in \ECsection{EC.6.1}.
To reflect practical time constraints, we set a 120-second optimization time limit for \(\mathbf H_1\), \(\mathbf H_2\), ISA, and NOMAD.
The reported runtimes include all computational steps and may therefore slightly exceed this limit.
Table~\ref{case_study_convergence} reports cumulative runtime across the initial 120-second phase and the extended restart.

\subsection{Value of Intervention}
\label{exp_sd}

We fix \(|\mathcal Z|=150\) and vary the initial dispatchable-supply-to-demand ratio \(\kappa\) from \(0.20\) to \(1.20\), covering settings from severe shortage to oversupply while holding all other parameters fixed.
Table~\ref{tab_sd} supports two distinct comparisons.
Relative to the no-bonus policy \(\mathbf B^0\), the driver surge map obtained by ISA improves the objective by \(62.17\%\) at \(\kappa=0.20\); this intervention value decreases steadily to \(12.99\%\) at \(\kappa=1.20\).
Within these modeled settings, bonuses therefore generate the greatest incremental value when immediately dispatchable supply is scarce.
This pattern is consistent with Corollary~\ref{cor:no_bonus_intervention}.
When more drivers are already matched in a zone at the no-bonus map, a small bonus increase entails larger payments on existing matches.
Introducing a bonus is therefore worthwhile only when the matching value generated through destination attraction, retention, and activation exceeds these additional payments.
ISA also outperforms the strongest non-ISA feasible benchmark by \(5.10\%\)--\(19.00\%\) and terminates in approximately 20 seconds, demonstrating the practical value of matching-informed bonus allocation within real-time decision windows.

\begin{table*}[htbp]
\centering
\caption{Performance under different supply--demand regimes.}
\label{tab_sd}
\footnotesize
\setlength{\tabcolsep}{2pt}
\renewcommand{\arraystretch}{1.40}
\resizebox{\textwidth}{!}{
\begin{tabular}{@{}*{16}{c}@{}}
\toprule
\multirow{2}{*}{Regime}
& \multirow{2}{*}{$\kappa$}
& \multicolumn{2}{c}{Conic}
& \multicolumn{2}{c}{$\mathbf B^0$}
& \multicolumn{2}{c}{$\mathbf H_1$}
& \multicolumn{2}{c}{$\mathbf H_2$}
& \multicolumn{2}{c}{NOMAD}
& \multicolumn{2}{c}{UB}
& \multicolumn{2}{c}{ISA} \\
\cmidrule(lr){3-4}
\cmidrule(lr){5-6}
\cmidrule(lr){7-8}
\cmidrule(lr){9-10}
\cmidrule(lr){11-12}
\cmidrule(lr){13-14}
\cmidrule(lr){15-16}
& & Obj. & Time & Obj. & Time & Obj. & Time & Obj. & Time
& Obj. & Time & Obj. & Time & Obj. & Time \\
\midrule
\multirow{2}{*}{Severe shortage}
& 0.20 & 4432.36 & 3.04 & 5815.88 & 0.13 & 6453.35 & 0.31 & 6234.41 & 0.42 & 8811.53 & 119.81 & 14319.46 & 2.91 & \best{9431.46}~\gain{7.04} & 20.12 \\
& 0.30 & 4139.56 & 3.08 & 7736.52 & 0.15 & 7549.75 & 0.31 & 8111.77 & 0.45 & 10298.45 & 119.84 & 16335.39 & 2.94 & \best{11281.64}~\gain{9.55} & 17.24 \\

\midrule
\multirow{2}{*}{Moderate shortage}
& 0.40 & 3800.15 & 3.39 & 9353.77 & 0.12 & 8244.08 & 0.34 & 9698.14 & 0.46 & 10990.66 & 119.85 & 17982.31 & 3.25 & \best{12832.67}~\gain{16.76} & 19.78 \\
& 0.50 & 3469.64 & 3.23 & 10804.30 & 0.13 & 9037.60 & 0.31 & 11111.25 & 0.43 & 12111.39 & 119.83 & 19403.43 & 3.10 & \best{14183.46}~\gain{17.11} & 18.94 \\

\midrule
\multirow{2}{*}{Mild shortage}
& 0.60 & 3150.65 & 3.09 & 12101.25 & 0.15 & 10106.02 & 0.34 & 12370.68 & 0.48 & 13982.57 & 119.89 & 20554.12 & 2.95 & \best{15408.77}~\gain{10.20} & 19.04 \\
& 0.70 & 2861.47 & 4.30 & 13274.30 & 0.13 & 11075.36 & 0.39 & 13498.27 & 0.71 & 15094.14 & 119.92 & 21333.82 & 4.16 & \best{16459.56}~\gain{9.05} & 19.63 \\

\midrule
\multirow{3}{*}{Near balanced}
& 0.80 & 2622.43 & 2.74 & 14310.55 & 0.13 & 11943.39 & 0.34 & 14487.80 & 0.46 & 14596.16 & 119.60 & 21647.29 & 2.61 & \best{17368.96}~\gain{19.00} & 17.42 \\
& 0.90 & 2427.80 & 2.98 & 15220.59 & 0.15 & 12524.17 & 0.31 & 15351.59 & 0.45 & 16931.16 & 119.89 & 21875.30 & 2.85 & \best{18125.72}~\gain{7.06} & 20.38 \\
& 1.00 & 2286.11 & 3.01 & 16003.50 & 0.13 & 13145.11 & 0.32 & 16003.50 & 0.28 & 17581.28 & 119.90 & 22076.35 & 2.88 & \best{18630.80}~\gain{5.97} & 16.94 \\

\midrule
\multirow{2}{*}{Oversupplied}
& 1.10 & 2143.48 & 2.96 & 16567.67 & 0.13 & 13813.78 & 0.34 & 16567.67 & 0.30 & 17927.68 & 119.87 & 22250.89 & 2.81 & \best{18842.11}~\gain{5.10} & 18.34 \\
& 1.20 & 2024.37 & 2.97 & 16852.29 & 0.15 & 14082.46 & 0.36 & 16852.29 & 0.28 & 18101.64 & 119.93 & 22400.38 & 2.84 & \best{19040.79}~\gain{5.19} & 18.49 \\
\bottomrule
\end{tabular}
}
\end{table*}

\subsection{Value of the Bonus Budget}
\label{exp_value}

To quantify the operational value of driver surge maps and examine how this value changes with the available bonus budget, we vary the budget over \(B\in\{0,300,500,1000,1200\}\) for instances with
\(|\mathcal Z|=150\), holding all other parameters fixed.
For each setting, we report the objective value (Obj.), objective improvement (Obj. Imp.), defined as \((\mathrm{Obj}(B)-\mathrm{Obj}(0))/\mathrm{Obj}(0)\times100\%\), and realized bonus spending (Spend.).
We also report the match rate, defined as \(\sum_{i\in\mathcal Z}\sum_{j\in\mathcal Z}x_{ij}/\sum_{j\in\mathcal Z}d_j\times100\%\), and return on investment (ROI), defined as \((\mathrm{Obj}(B)-\mathrm{Obj}(0))/\mathrm{Spend}(B)\times100\).
Table~\ref{result_bonus_value} shows that larger bonus budgets improve the attained objective and match rate in the tested settings, while ROI declines as spending increases.
At \(B=1200\), the objective improvement is \(42.17\%\) under \(\kappa=0.50\), compared with \(16.63\%\) under \(\kappa=1.10\), indicating that incentives are most valuable when supply is scarce.
Moreover, realized bonus spending is below the nominal budget of \(B=1200\) in the balanced and oversupplied settings.
These results suggest that the platform can use the incremental net operating profit from budget expansion to select an appropriate bonus budget, identifying when further increases in spending deliver only modest additional gains.

\begin{table*}[htbp]
\centering
\caption{Value of driver surge maps under different bonus budgets.}
\label{result_bonus_value}

{\scriptsize
\compacttablerows
\setlength{\tabcolsep}{2pt}

\begin{tabular*}{\textwidth}{@{\extracolsep{\fill}}c c c c c c c@{}}
\toprule
\(\kappa\) & Bonus budget & Obj. & Obj. Imp. (\%) &
Spend. & Match rate (\%) & ROI (\%)\\
\midrule
\multirow{5}{*}{\(0.50\)}
& \(0\) & \(10804\) & \(0.00\) & \(0.00\) & \(45.59\) & -- \\
& \(300\) & \(13367\) & \(23.72\) & \(300.00\) & \(54.51\) & \(854.27\) \\
& \(500\) & \(14183\) & \(31.28\) & \(500.00\) & \(57.74\) & \(675.83\) \\
& \(1000\) & \(15190\) & \(40.60\) & \(1000.00\) & \(62.72\) & \(438.61\) \\
& \(1200\) & \(15360\) & \(42.17\) & \(1199.99\) & \(64.04\) & \(379.64\) \\
\midrule
\multirow{5}{*}{\(1.00\)}
& \(0\) & \(16004\) & \(0.00\) & \(0.00\) & \(87.61\) & -- \\
& \(300\) & \(18067\) & \(12.90\) & \(300.00\) & \(96.30\) & \(687.99\) \\
& \(500\) & \(18631\) & \(16.42\) & \(500.00\) & \(99.30\) & \(525.46\) \\
& \(1000\) & \(19106\) & \(19.39\) & \(1000.00\) & \(100.00\) & \(310.26\) \\
& \(1200\) & \(19144\) & \(19.63\) & \(1189.80\) & \(100.00\) & \(263.98\) \\
\midrule
\multirow{5}{*}{\(1.10\)}
& \(0\) & \(16568\) & \(0.00\) & \(0.00\) & \(96.01\) & -- \\
& \(300\) & \(18352\) & \(10.77\) & \(300.00\) & \(100.00\) & \(594.75\) \\
& \(500\) & \(18842\) & \(13.73\) & \(500.00\) & \(100.00\) & \(454.89\) \\
& \(1000\) & \(19296\) & \(16.47\) & \(1000.00\) & \(100.00\) & \(272.81\) \\
& \(1200\) & \(19323\) & \(16.63\) & \(1155.99\) & \(100.00\) & \(238.34\) \\
\bottomrule
\end{tabular*}
}
\compacttableend
\end{table*}

\subsection{Behavioral and Spatial Sensitivity Analyses}
\label{exp_sensi}

We next examine how latent-driver activation responsiveness, behavioral-parameter estimation error, and spatial smoothness affect the benefits of bonus intervention. Each analysis varies one factor while holding the remaining experimental settings fixed.

\subsubsection{Responsiveness of Latent Driver Activation.}
\label{result_response_sensitivity}

We vary \(\gamma^{\mathrm{On}}\) over \(\{1.0,1.2,1.5,1.8,2.0\}\) while holding all other parameters fixed; \(\gamma^{\mathrm{On}}=1.5\) is the baseline, and higher values represent stronger activation responses among latent drivers.
Table~\ref{result_response_intensity_scaling} shows that stronger activation responsiveness is associated with higher attained objective and ROI across the tested settings.
As \(\gamma^{\mathrm{On}}\) increases from \(1.0\) to \(2.0\), the objective improvement rises from \(23.07\%\) to \(37.54\%\) with \(\kappa=0.50\) and from \(9.95\%\) to \(16.55\%\) with \(\kappa=1.10\).
The objective can continue to increase even when the aggregate match rate has reached \(100\%\), because newly dispatchable supply may support higher-value matches.
Equation~\eqref{eq:explicit_activation_threshold} shows that, holding the baseline activation probability fixed, stronger
activation responsiveness lowers the intervention threshold.
A small initial bonus improves net operating profit when the pre-bonus matching value exceeds this threshold.
The value of activating latent drivers therefore depends on both their responsiveness and the matching opportunities available to the platform.

\begin{table*}[htbp]
\centering
\caption{Sensitivity analysis with respect to latent-driver activation responsiveness.}
\label{result_response_intensity_scaling}

{\scriptsize
\compacttablerows
\setlength{\tabcolsep}{2pt}

\begin{tabular*}{\textwidth}{@{\extracolsep{\fill}}c c c c c c c@{}}
\toprule
\(\kappa\) & \(\gamma^{\mathrm{On}}\) & Obj. & Obj. Imp. (\%) &
Spend. & Match rate (\%) & ROI (\%)\\
\midrule
\multirow{5}{*}{\(0.50\)}
& 1.0 & 13297 & 23.07 & 500.00 & 54.91 & 498.46 \\
& 1.2 & 13683 & 26.64 & 500.00 & 56.13 & 575.76 \\
& 1.5 & 14183 & 31.28 & 500.00 & 57.74 & 675.83 \\
& 1.8 & 14609 & 35.21 & 500.00 & 59.09 & 760.88 \\
& 2.0 & 14860 & 37.54 & 500.00 & 59.91 & 811.09 \\
\midrule
\multirow{5}{*}{\(1.00\)}
& 1.0 & 17936 & 12.08 & 500.00 & 96.50 & 386.49 \\
& 1.2 & 18251 & 14.05 & 500.00 & 97.76 & 449.59 \\
& 1.5 & 18631 & 16.42 & 500.00 & 99.30 & 525.46 \\
& 1.8 & 18931 & 18.29 & 500.00 & 100.00 & 585.50 \\
& 2.0 & 19103 & 19.37 & 500.00 & 100.00 & 619.99 \\
\midrule
\multirow{5}{*}{\(1.10\)}
& 1.0 & 18216 & 9.95 & 500.00 & 100.00 & 329.71 \\
& 1.2 & 18502 & 11.67 & 500.00 & 100.00 & 386.81 \\
& 1.5 & 18842 & 13.73 & 500.00 & 100.00 & 454.89 \\
& 1.8 & 19144 & 15.55 & 500.00 & 100.00 & 515.25 \\
& 2.0 & 19309 & 16.55 & 500.00 & 100.00 & 548.25 \\
\bottomrule
\end{tabular*}
}
\compacttableend
\end{table*}
\subsubsection{Error in Behavioral Parameter Estimates.} \label{para_robust}
In deployment, behavioral-response parameters calibrated from historical data may differ from current driver response because of limited observations and changing market conditions.
To examine the impact of these inaccuracies, we use the instances with \(\kappa\in\{0.50,1.00,1.10\}\) from Section~\ref{exp_sd} and independently perturb the six common behavioral coefficients \(\beta^{(1)},\beta^{(3)},\delta^{\mathrm{Off}},\zeta^{\mathrm{Off}},\alpha^{\mathrm{On}}\), and \(\gamma^{\mathrm{On}}\) by factors drawn from \(\mathrm{Uniform}(1-r,1+r)\), where \(r\in\{10\%,20\%,30\%,40\%,50\%\}\).
For these instances, \(\delta_i^{\mathrm{Off}}\equiv\delta^{\mathrm{Off}}\) and
\(\alpha_i^{\mathrm{On}}\equiv\alpha^{\mathrm{On}}\) across all zones.
For each \((\kappa,r)\), ISA constructs driver surge maps using 10 independently perturbed parameter vectors, which are then evaluated under the nominal response. Let \(Z^{\mathrm{nom}}_{\kappa}\) denote the nominal-parameter objective and \(Z^{\mathrm{deploy}}_{r,\kappa,n}\) the deployed objective under perturbation \(n\); we define \(\mathrm{Loss}_{r,\kappa,n}=(Z^{\mathrm{nom}}_{\kappa}-Z^{\mathrm{deploy}}_{r,\kappa,n})/|Z^{\mathrm{nom}}_{\kappa}|\times100\%\). Table~\ref{tab_rob_para} reports the mean deployed objective and the mean, standard deviation, and maximum loss across the 10 perturbations.
The mean loss remains between \(0.24\%\) and \(3.25\%\) across all settings. The largest mean losses occur under supply shortage and at the higher perturbation levels, indicating that accurate behavioral calibration is most valuable when supply is scarce.

\begin{table*}[htbp]
\centering
\caption{Sensitivity analysis with respect to behavioral-parameter estimation error.}
\label{tab_rob_para}

{\scriptsize
\compacttablerows
\setlength{\tabcolsep}{2pt}

\begin{tabular*}{\textwidth}{
@{\extracolsep{\fill}}
c c c c c c
@{}}
\toprule
\(\kappa\)
& Perturbation level
& Mean Obj.
& Mean loss (\%)
& Std. loss (\%)
& Max loss (\%) \\
\midrule
\multirow{6}{*}{\(0.50\)}
& \(0\%\)      & 14183 & 0.00 & --   & --   \\
& \(\pm10\%\) & 14118 & 0.46 & 0.59 & 2.08 \\
& \(\pm20\%\) & 14058 & 0.89 & 0.71 & 2.12 \\
& \(\pm30\%\) & 13833 & 2.47 & 1.70 & 5.20 \\
& \(\pm40\%\) & 13840 & 2.42 & 2.04 & 7.63 \\
& \(\pm50\%\) & 13722 & 3.25 & 2.49 & 7.64 \\
\midrule
\multirow{6}{*}{\(1.00\)}
& \(0\%\)      & 18631 & 0.00 & --   & --   \\
& \(\pm10\%\) & 18586 & 0.24 & 0.16 & 0.49 \\
& \(\pm20\%\) & 18540 & 0.49 & 0.28 & 0.99 \\
& \(\pm30\%\) & 18502 & 0.69 & 0.50 & 1.55 \\
& \(\pm40\%\) & 18313 & 1.70 & 1.09 & 3.80 \\
& \(\pm50\%\) & 18223 & 2.19 & 1.14 & 4.45 \\
\midrule
\multirow{6}{*}{\(1.10\)}
& \(0\%\)      & 18842 & 0.00 & --   & --   \\
& \(\pm10\%\) & 18788 & 0.29 & 0.31 & 0.96 \\
& \(\pm20\%\) & 18717 & 0.66 & 0.48 & 1.33 \\
& \(\pm30\%\) & 18740 & 0.54 & 0.50 & 1.81 \\
& \(\pm40\%\) & 18594 & 1.32 & 1.02 & 3.26 \\
& \(\pm50\%\) & 18676 & 0.88 & 0.52 & 2.13 \\
\bottomrule
\end{tabular*}
}
\compacttableend
\end{table*}
\subsubsection{Effect of Spatial Smoothness.}
\label{exp_sm}

To examine how spatial smoothness affects solution performance, we vary the spatial smoothness threshold over \(\Delta^b\in\{0.5,1.0,2.0,+\infty\}\), where larger values permit greater spatial differentiation. We measure spatial bonus variation as the mean of \(|b_i-b_j|\) over neighboring zone pairs. Table~\ref{tab:smoothness_sensitivity} shows that relaxing the threshold initially improves performance, but the gains diminish quickly.
When \(\kappa=0.50\), increasing \(\Delta^b\) from \(0.5\) to \(2.0\) raises the objective improvement from \(28.61\%\) to \(31.99\%\), whereas removing the constraint increases bonus variation without further improvement. Similar patterns arise when \(\kappa=1.00\) and \(1.10\), suggesting that a moderate threshold preserves most of the targeting value while avoiding sharp bonus differences across neighboring zones.
Corollary~\ref{cor:spatial_component_balance} explains how smoothness changes the first-order targeting rule.
Edges with positive smoothness multipliers can sustain opposing local marginal-value pressures, but those pressures must balance
over the associated component unless a bonus bound absorbs them.
The numerical results measure the attained value of relaxing the threshold, whereas the corollary identifies the
spatial mechanism through which this relaxation can matter.

\begin{table*}[htbp]
\centering
\caption{Sensitivity analysis with respect to the spatial smoothness threshold.}
\label{tab:smoothness_sensitivity}

{\scriptsize
\compacttablerows
\setlength{\tabcolsep}{2pt}

\begin{tabular*}{\textwidth}
{@{\extracolsep{\fill}}c c c c c c c c@{}}
\toprule
\(\kappa\) & \(\Delta^b\) & Obj. & Obj. Imp. (\%) &
Match rate (\%) & Mean \(|b_i-b_j|\) & ROI (\%) & Time \\
\midrule
\multirow{4}{*}{\(0.50\)}
& \(0.5\)      & 13895 & 28.61 & 56.94 & 0.221 & 618.12 & 19.11 \\
& \(1.0\)      & 14183 & 31.28 & 57.74 & 0.249 & 675.83 & 18.94 \\
& \(2.0\)      & 14261 & 31.99 & 58.02 & 0.262 & 691.31 & 15.67 \\
& \(+\infty\)  & 13936 & 28.98 & 56.86 & 0.533 & 626.31 & 20.89 \\
\midrule
\multirow{4}{*}{\(1.00\)}
& \(0.5\)      & 18336 & 14.57 & 98.12 & 0.210 & 466.43 & 16.06 \\
& \(1.0\)      & 18631 & 16.42 & 99.30 & 0.251 & 525.46 & 16.94 \\
& \(2.0\)      & 18714 & 16.94 & 99.76 & 0.232 & 542.15 & 21.36 \\
& \(+\infty\)  & 18713 & 16.93 & 99.76 & 0.232 & 541.99 & 19.44 \\
\midrule
\multirow{4}{*}{\(1.10\)}
& \(0.5\)      & 18568 & 12.07 & 100.00 & 0.212 & 400.06 & 15.94 \\
& \(1.0\)      & 18842 & 13.73 & 100.00 & 0.266 & 454.89 & 18.34 \\
& \(2.0\)      & 18929 & 14.25 & 100.00 & 0.240 & 472.21 & 18.88 \\
& \(+\infty\)  & 18929 & 14.25 & 100.00 & 0.241 & 472.30 & 20.55 \\
\bottomrule
\end{tabular*}
}
\compacttableend
\end{table*}
\subsection{Anytime Performance}
\label{exp_anytime}

To assess performance within the 120-second operating window,
we compare the incumbent trajectories of ISA and NOMAD across
the 26 unique generated instances used in the supply--demand
regime and scalability experiments.
For method \(m\), let \(Z_m(t)\) denote the objective value of
its latest feasible incumbent available by time \(t\), with
the incumbent carried forward if the method terminates early.
We measure ISA's objective lead over NOMAD at time \(t\) as
\(L(t):=[Z_{\mathrm{ISA}}(t)-Z_{\mathrm{NOMAD}}(t)]/
|Z_{\mathrm{NOMAD}}(t)|\times100\%\).
We also compare how quickly the two methods reach a common
objective target.
Let \(Z_{\max}^{120}:=\max\{Z_{\mathrm{ISA}}(120),
Z_{\mathrm{NOMAD}}(120)\}\), and define
\(Z_\alpha:=Z^{\mathrm{no}}
+\alpha(Z_{\max}^{120}-Z^{\mathrm{no}})\)
for \(\alpha\in\{0.90,0.95,0.99\}\), where
\(Z^{\mathrm{no}}\) is the no-bonus objective.
Thus, the target represents a fraction \(\alpha\) of the
improvement from the no-bonus objective to the better
120-second incumbent.
The target time and number of matching-LP solves are recorded
when a method first obtains an incumbent satisfying
\(Z_m(t)\ge Z_\alpha\).
These targets measure progress relative to the observed
incumbents and are distinct from the solver-certified upper bound.

\begin{figure}[htbp]
    \centering
    \includegraphics[width=0.78\textwidth]{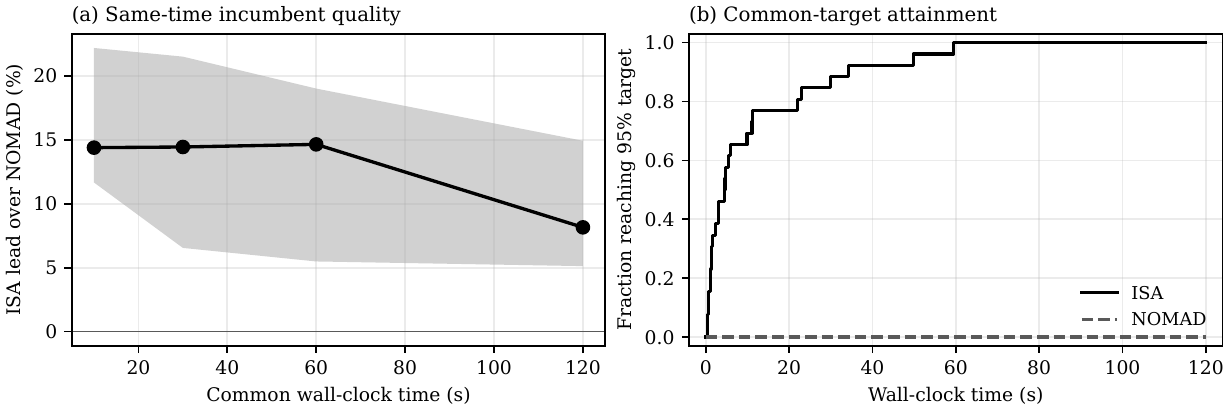}
    \caption{Anytime performance on 26 unique generated instances.}
    \label{fig:anytime_operational}
    \compactfigureend
\end{figure}

Figure~\ref{fig:anytime_operational} summarizes the incumbent
trajectories and attainment of the common \(95\%\) target.
ISA reaches this target in a median of \(4.52\) seconds,
whereas NOMAD remains below it throughout the operating window
on all 26 instances.
By 30 seconds, ISA's incumbent also exceeds the best completed
non-ISA feasible solution on every instance.
Detailed results, including matching-LP solve counts,
are reported in \ECsection{EC.6.2}.

\subsection{Scalability}
\label{exp_scale}
To evaluate the scalability of ISA, we conduct experiments by varying the number of zones over
\(|\mathcal Z|\in\{50,100,150,200,250,300\}\) under moderate-shortage, balanced, and oversupplied regimes, while holding all other parameters fixed.
Table~\ref{result_scale} shows that ISA achieves the highest feasible objective in all instances, with improvements of \(9.60\%\)--\(15.71\%\) even at 300 zones. ISA terminates by its stopping rule within 120 seconds in all but one instance, where it retains the best feasible incumbent. These results support real-time incentive optimization at larger spatial scales.

\begin{table*}[htbp]
\centering
\caption{Scalability across problem sizes.}
\label{result_scale}
\scriptsize
\setlength{\tabcolsep}{2pt}
\renewcommand{\arraystretch}{1.10}
\begin{tabular*}{\textwidth}{@{\extracolsep{\fill}}*{16}{c}@{}}
\toprule
\multirow{2}{*}{$|\mathcal{Z}|$}
& \multirow{2}{*}{$\kappa$}
& \multicolumn{2}{c}{Conic}
& \multicolumn{2}{c}{$\mathbf B^0$}
& \multicolumn{2}{c}{$\mathbf H_1$}
& \multicolumn{2}{c}{$\mathbf H_2$}
& \multicolumn{2}{c}{NOMAD}
& \multicolumn{2}{c}{UB}
& \multicolumn{2}{c}{ISA} \\
\cmidrule(lr){3-4}
\cmidrule(lr){5-6}
\cmidrule(lr){7-8}
\cmidrule(lr){9-10}
\cmidrule(lr){11-12}
\cmidrule(lr){13-14}
\cmidrule(lr){15-16}
& & Obj. & Time & Obj. & Time & Obj. & Time & Obj. & Time
& Obj. & Time & Obj. & Time & Obj. & Time \\
\midrule

\multirow{3}{*}{50}
& 0.50 & 820 & 0.59 & 7052 & 0.01 & 7374 & 0.04 & 7258 & 0.05 & 7950 & 49.75 & 12540 & 0.57 & \best{8393}~\gain{5.58} & 1.76 \\
& 1.00 & 1151 & 0.31 & 13490 & 0.01 & 13495 & 0.04 & 13646 & 0.05 & 14098 & 74.66 & 15409 & 0.29 & \best{14186}~\gain{0.63} & 9.88 \\
& 1.10 & 1151 & 0.27 & 14644 & 0.01 & 14523 & 0.04 & 14734 & 0.05 & 14881 & 74.31 & 15409 & 0.26 & \best{14934}~\gain{0.36} & 1.66 \\

\midrule
\multirow{3}{*}{100}
& 0.50 & 2957 & 1.46 & 6628 & 0.10 & 4909 & 0.14 & 6854 & 0.18 & 8732 & 119.22 & 12892 & 1.40 & \best{9367}~\gain{7.26} & 9.36 \\
& 1.00 & 1951 & 1.62 & 10837 & 0.06 & 6316 & 0.19 & 10929 & 0.22 & 12416 & 119.85 & 15190 & 1.53 & \best{12850}~\gain{3.49} & 6.04 \\
& 1.10 & 1856 & 1.48 & 11386 & 0.06 & 6373 & 0.15 & 11386 & 0.13 & 12655 & 119.84 & 15310 & 1.43 & \best{13030}~\gain{2.96} & 6.26 \\

\midrule
\multirow{3}{*}{150}
& 0.50 & 3470 & 3.23 & 10804 & 0.13 & 9038 & 0.31 & 11111 & 0.43 & 12111 & 119.83 & 19403 & 3.10 & \best{14183}~\gain{17.11} & 18.94 \\
& 1.00 & 2286 & 3.01 & 16004 & 0.17 & 13145 & 0.32 & 16004 & 0.28 & 17581 & 119.90 & 22076 & 2.88 & \best{18631}~\gain{5.97} & 16.94 \\
& 1.10 & 2143 & 2.96 & 16568 & 0.14 & 13814 & 0.34 & 16568 & 0.30 & 17928 & 119.87 & 22251 & 2.81 & \best{18842}~\gain{5.10} & 18.34 \\

\midrule
\multirow{3}{*}{200}
& 0.50 & 2410 & 6.02 & 14676 & 0.24 & 12745 & 0.59 & 15086 & 0.78 & 15407 & 119.86 & 26239 & 5.79 & \best{18392}~\gain{19.37} & 35.28 \\
& 1.00 & 1667 & 5.06 & 19373 & 0.28 & 17841 & 0.64 & 19373 & 0.55 & 21832 & 119.98 & 28697 & 4.83 & \best{22295}~\gain{2.12} & 36.07 \\
& 1.10 & 1578 & 5.10 & 19792 & 0.29 & 18510 & 0.60 & 19792 & 0.55 & 22185 & 120.04 & 28994 & 4.87 & \best{22613}~\gain{1.93} & 41.40 \\

\midrule
\multirow{3}{*}{250}
& 0.50 & 2820 & 9.55 & 17455 & 0.44 & 16799 & 0.93 & 17845 & 1.32 & 18249 & 119.92 & 29753 & 9.13 & \best{21166}~\gain{15.98} & 62.50 \\
& 1.00 & 2009 & 9.26 & 21319 & 0.47 & 20193 & 1.00 & 21319 & 0.95 & 21319 & 120.12 & 32432 & 8.85 & \best{24530}~\gain{15.06} & 87.11 \\
& 1.10 & 1959 & 8.46 & 21741 & 0.40 & 20656 & 0.97 & 21741 & 0.91 & 21741 & 120.09 & 32755 & 8.06 & \best{24895}~\gain{14.51} & 103.57 \\

\midrule
\multirow{3}{*}{300}
& 0.50 & 2623 & 13.54 & 21825 & 0.66 & 20427 & 1.52 & 22454 & 2.12 & 22777 & 120.00 & 40192 & 12.90 & \best{26354}~\gain{15.71} & 106.66 \\
& 1.00 & 1714 & 12.04 & 31032 & 0.62 & 29925 & 1.56 & 31032 & 1.38 & 31032 & 120.07 & 44870 & 11.43 & \best{34285}~\gain{10.48} & 98.02 \\
& 1.10 & 1636 & 13.79 & 31862 & 0.62 & 30910 & 1.57 & 31862 & 1.40 & 31862 & 120.08 & 45306 & 13.17 & \best{34922}~\gain{9.60} & 120.08 \\
\bottomrule
\end{tabular*}
\end{table*}

\subsection{New York City Study}
\label{exp_case}
The generated instances allow us to evaluate the model under controlled operating conditions. 
We next complement these experiments with a data-based case study to assess ISA in a realistic large-scale setting. Specifically, we use 2015 New York City Yellow Taxi trip records released by the Taxi and Limousine Commission \citep{nyc_tlc_trip_data}. 
The full preprocessing region comprises 2,350 Geohash-6 cells; each daily instance retains cells appearing as pickup or dropoff locations, leaving approximately 1,400 active zones.
Demand equals the number of pickups in a two-minute target window, and initial dispatchable supply is approximated by the number of dropoffs in the preceding minute.
Latent supply is set to \(0.3\) times total demand and allocated across zones in proportion to demand.
The pre-bonus matching value is \(p_{ij}=15+10\exp(-d_{ij}/3)\), where \(d_{ij}\) is the Haversine distance in kilometers.

For a daily instance with \(n=|\mathcal Z|\), let \(m=\lceil\sqrt n\rceil\).
We partition the latitude--longitude bounding box of the observed zone centers into an \(m\times m\) equal-width grid and connect zones whose cells coincide or are horizontally, vertically, or diagonally adjacent.
For these instances, the \(3\times3\) cell stencil defines the sparse graph used for driver relocation, spatial smoothness, and matching.
The study combines taxi-based geography, demand, and a dispatchable-supply proxy with specified latent supply, matching values, and driver response parameters.

We select three dates for each of the intervals 07:00--07:02, 10:00--10:02, and 17:00--17:02 to represent morning peak, off-peak, and evening peak conditions. The nine instances contain 1,396--1,486 active zones and preserve temporal variation and heterogeneous spatial imbalances.
Table~\ref{details_case} summarizes these instances and uses M, O, and E to denote the three market conditions, where
\(D^{\mathrm{tot}}=\sum_{i\in\mathcal Z}d_i\),
\(S^{\mathrm{tot}}=\sum_{i\in\mathcal Z}s_i^{\mathrm{On}}\), and
\(D^{\mathrm{gap}}=\sum_{i\in\mathcal Z}\max\{d_i-s_i^{\mathrm{On}},0\}\)
denote total rider demand, the preceding-minute dropoff proxy for dispatchable supply, and the total local shortage gap, respectively.
\begin{table}[htbp]
\centering
\caption{Details of the New York City instances.}
\label{details_case}
\scriptsize
\setlength{\tabcolsep}{2pt}
\compacttablerows
\begin{tabular*}{\textwidth}{@{\extracolsep{\fill}}c c c c c c c c c@{}}
\toprule
Type
& ID
& Date
& Period
& \(|\mathcal Z|\)
& \(D^{\mathrm{tot}}\)
& \(S^{\mathrm{tot}}\)
& \(D^{\mathrm{gap}}\)
& \(\kappa\) \\
\midrule
\multirow{3}{*}{Morning peak}
& M1 & 2015-01-20 & \multirow{3}{*}{07:00--07:02} & 1417 & 521 & 278 & 338 & 0.53 \\
& M2 & 2015-01-22 &                              & 1486 & 574 & 274 & 384 & 0.48 \\
& M3 & 2015-01-29 &                              & 1459 & 570 & 267 & 377 & 0.47 \\
\midrule
\multirow{3}{*}{Off-peak}
& O1 & 2015-01-08 & \multirow{3}{*}{10:00--10:02} & 1485 & 680 & 380 & 363 & 0.56 \\
& O2 & 2015-01-18 &                              & 1464 & 539 & 292 & 286 & 0.54 \\
& O3 & 2015-01-24 &                              & 1418 & 550 & 314 & 281 & 0.57 \\
\midrule
\multirow{3}{*}{Evening peak}
& E1 & 2015-01-03 & \multirow{3}{*}{17:00--17:02} & 1396 & 648 & 343 & 350 & 0.53 \\
& E2 & 2015-01-07 &                              & 1428 & 614 & 323 & 341 & 0.53 \\
& E3 & 2015-01-17 &                              & 1421 & 680 & 386 & 359 & 0.57 \\
\bottomrule
\end{tabular*}
\compacttableend
\end{table}
We compare ISA with the benchmarks in Section~\ref{exp_benchmarks} using a total bonus budget of \(B=1500\) for every instance; both one-pass heuristics complete within the 120-second optimization time limit.
Table~\ref{case_study_benchmark} shows that, under the nominal 120-second optimization deadline, ISA attains the best objective in all nine instances and improves upon the strongest non-ISA feasible benchmark by \(2.38\%\)--\(8.43\%\).
NOMAD retains its \(\mathbf H_1\) starting solution throughout the run, whereas the Conic benchmark provides a recovered feasible solution and the conic outer relaxation provides a solver-certified upper bound for every instance.

\begin{table}[htbp]
\centering
\caption{Performance comparison in the New York City study.}
\label{case_study_benchmark}
\footnotesize
\setlength{\tabcolsep}{2pt}
\renewcommand{\arraystretch}{0.8}
\resizebox{\textwidth}{!}{
\begin{tabular}{@{}*{16}{c}@{}}
\toprule
\multirow{2}{*}{Type}
& \multirow{2}{*}{ID}
& \multicolumn{2}{c}{Conic}
& \multicolumn{2}{c}{$\mathbf B^0$}
& \multicolumn{2}{c}{$\mathbf H_1$}
& \multicolumn{2}{c}{$\mathbf H_2$}
& \multicolumn{2}{c}{NOMAD}
& \multicolumn{2}{c}{UB}
& \multicolumn{2}{c}{ISA} \\
\cmidrule(lr){3-4}
\cmidrule(lr){5-6}
\cmidrule(lr){7-8}
\cmidrule(lr){9-10}
\cmidrule(lr){11-12}
\cmidrule(lr){13-14}
\cmidrule(lr){15-16}
& & Obj. & Time & Obj. & Time & Obj. & Time & Obj. & Time
& Obj. & Time & Obj. & Time & Obj. & Time \\
\midrule

\multirow{3}{*}{Morning peak}
& M1 & 5796 & 157.15 & 6576 & 0.38 & 8224 & 2.82 & 6730 & 1.44
& 8224 & 119.85 & 10495 & 154.49
& \best{8703}~\gain{5.83} & 120.27 \\
& M2 & 5251 & 128.49 & 6548 & 0.39 & 8276 & 2.75 & 6727 & 1.50
& 8276 & 120.16 & 10828 & 128.10
& \best{8974}~\gain{8.43} & 120.47 \\
& M3 & 5107 & 107.96 & 6466 & 0.39 & 8371 & 2.72 & 6643 & 1.46
& 8371 & 120.17 & 10688 & 107.58
& \best{8881}~\gain{6.10} & 120.72 \\

\midrule
\multirow{3}{*}{Off-peak}
& O1 & 5284 & 136.89 & 9157 & 0.42 & 11348 & 2.88 & 9354 & 1.58
& 11348 & 120.10 & 14322 & 136.48
& \best{11863}~\gain{4.54} & 120.37 \\
& O2 & 5009 & 126.37 & 6994 & 0.39 & 8777 & 2.71 & 7154 & 1.47
& 8777 & 120.16 & 11050 & 124.75
& \best{9180}~\gain{4.60} & 120.39 \\
& O3 & 5781 & 232.43 & 7616 & 0.41 & 9555 & 2.88 & 7771 & 1.47
& 9555 & 119.96 & 11798 & 231.77
& \best{9785}~\gain{2.40} & 120.18 \\

\midrule
\multirow{3}{*}{Evening peak}
& E1 & 5284 & 131.10 & 8191 & 0.37 & 10479 & 2.53 & 8383 & 1.30
& 10479 & 119.79 & 13113 & 130.73
& \best{10829}~\gain{3.34} & 120.14 \\
& E2 & 5095 & 104.61 & 7760 & 0.39 & 10015 & 2.69 & 7941 & 1.42
& 10015 & 120.16 & 12359 & 104.22
& \best{10253}~\gain{2.38} & 120.18 \\
& E3 & 5231 & 130.39 & 8969 & 0.38 & 11241 & 2.73 & 9165 & 1.44
& 11241 & 119.96 & 14193 & 129.99
& \best{11715}~\gain{4.22} & 120.50 \\

\bottomrule
\end{tabular}
}
\compacttableend
\end{table}
To assess objective stability beyond the 120-second operating window, we compare the cutoff incumbents with those from extended runs.
Table~\ref{case_study_convergence} shows that the extended runs yield additional objective gains of up to \(0.17\%\), indicating limited improvement beyond the operating window for the tested instances.
The columns Obj. (120 s) and Obj. (Ext.) report the cutoff and extended-run objectives, respectively, and Ext. gain reports their relative difference.
Cum. time denotes the total wall-clock time across both phases.

\begin{figure}[H]
	\centering
\captionsetup[subfigure]{font=small}
\begin{subfigure}{0.46\linewidth}
		\centering
		\includegraphics[width=\linewidth]{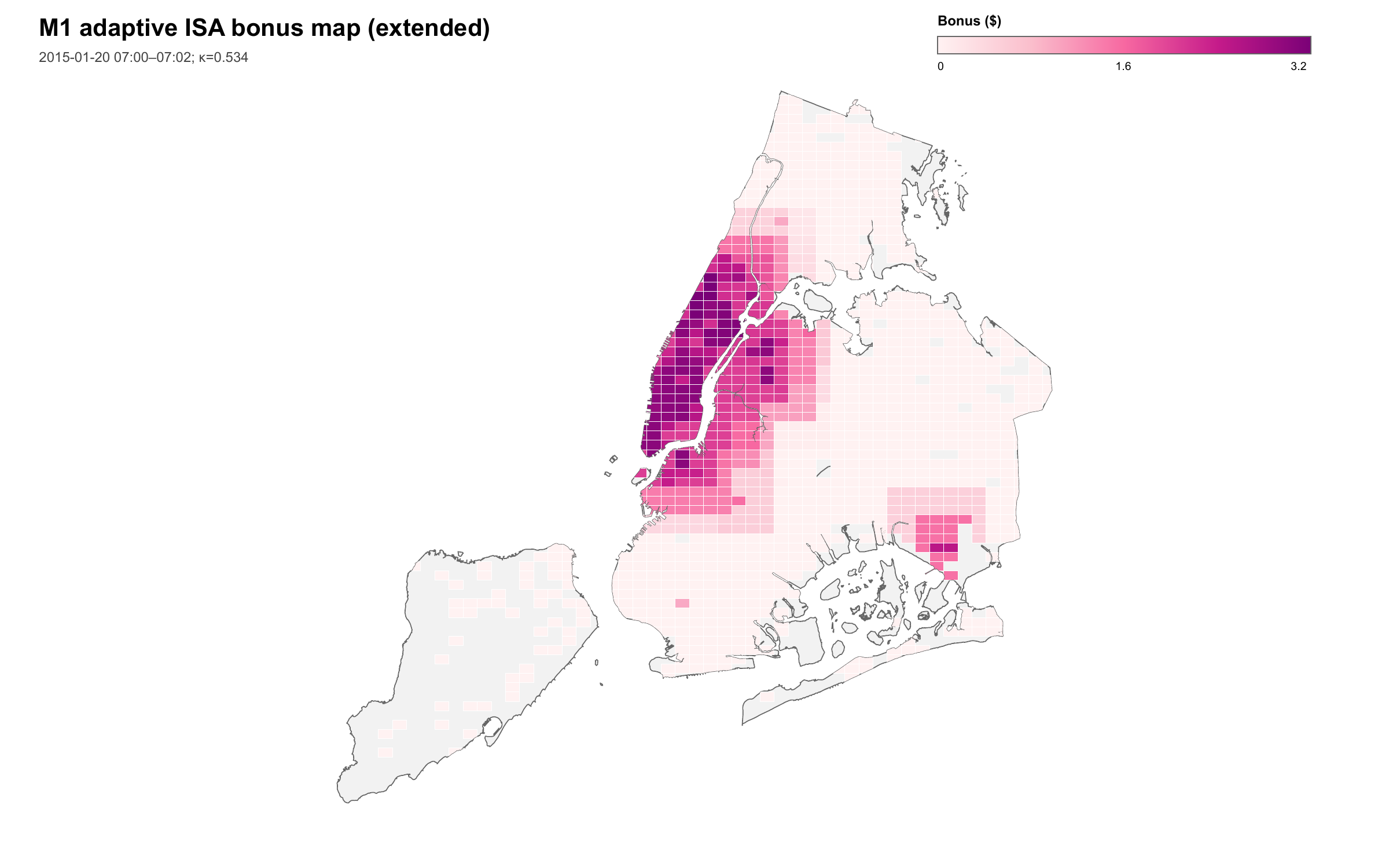}
	\caption{M1}
	\label{M1}
\end{subfigure}
\hspace{0.01\linewidth}
    \begin{subfigure}{0.46\linewidth}
		\centering
		\includegraphics[width=\linewidth]{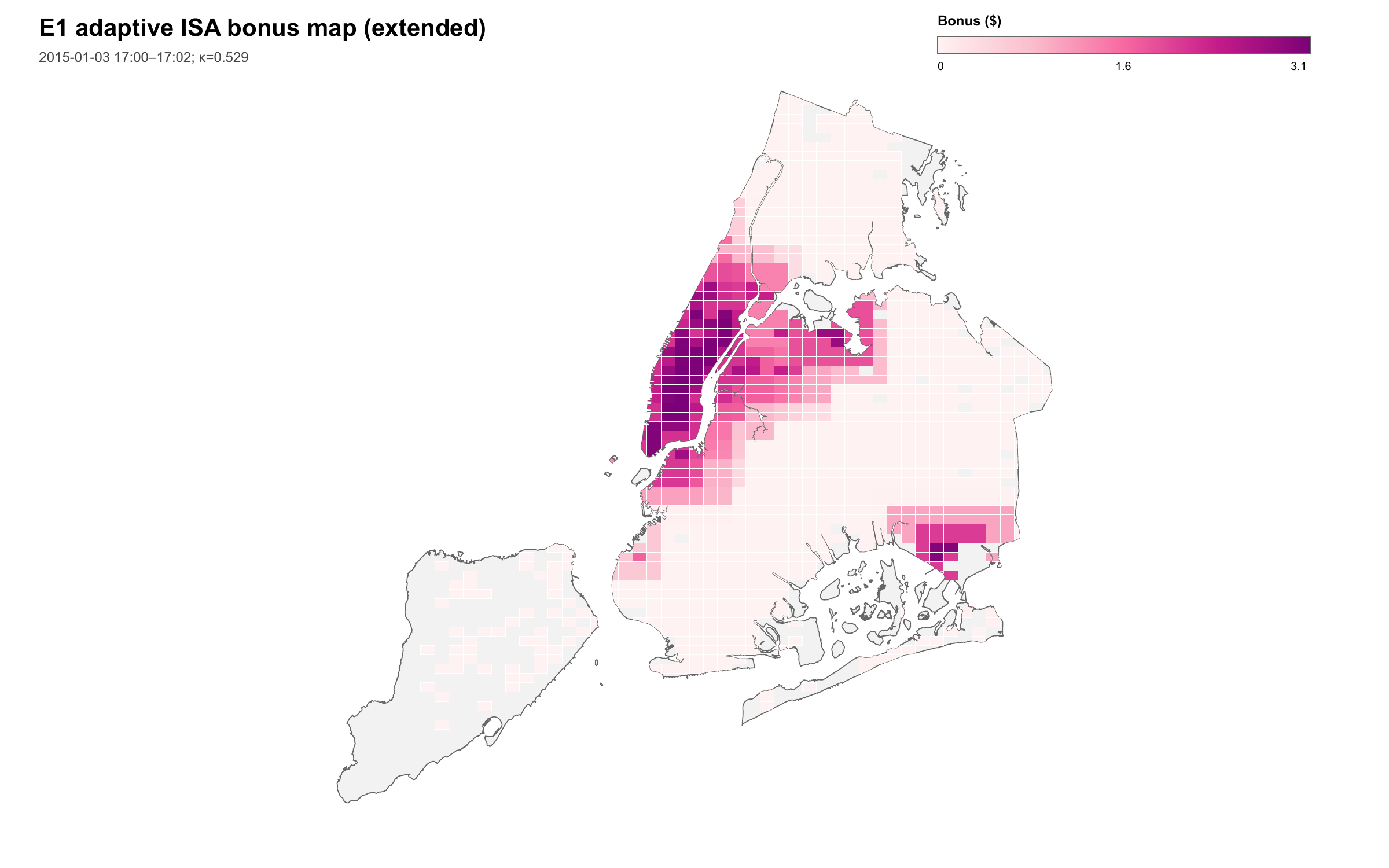}
		\caption{E1}
		\label{E1}
	\end{subfigure}
\caption{Optimized driver surge maps for two New York City instances.}
	\label{fig_nyc}
\compactfigureend
\end{figure}

\begin{table*}[htbp]
\centering
\caption{Comparison of 120-second and extended ISA incumbents in the New York City study.}
\label{case_study_convergence}
{\scriptsize
\compacttablerows
\setlength{\tabcolsep}{2pt}

\begin{tabular*}{\textwidth}{
@{\extracolsep{\fill}}*{6}{c}@{}
}
\toprule
Type
& ID
& Obj. (120 s)
& Obj. (Ext.)
& Ext. gain (\%)
& Cum. time \\
\midrule

\multirow{3}{*}{Morning peak}
& M1 & 8703 & 8708 & 0.05 & 1413.16 \\
& M2 & 8974 & 8989 & 0.16 & 1484.69 \\
& M3 & 8881 & 8885 & 0.04 & 1497.22 \\

\midrule
\multirow{3}{*}{Off-peak}
& O1 & 11863 & 11883 & 0.17 & 1530.00 \\
& O2 & 9180 & 9185 & 0.05 & 1485.62 \\
& O3 & 9785 & 9792 & 0.07 & 1425.45 \\

\midrule
\multirow{3}{*}{Evening peak}
& E1 & 10829 & 10831 & 0.02 & 1394.55 \\
& E2 & 10253 & 10254 & 0.01 & 1432.34 \\
& E3 & 11715 & 11718 & 0.02 & 1448.75 \\
\bottomrule
\end{tabular*}
}
\compacttableend
\end{table*}
Beyond the objective comparisons, the optimized maps show how ISA allocates bonuses across space. Figure~\ref{fig_nyc} shows the ISA maps for the New York City M1 and E1 instances. Within each panel, darker shading denotes a higher driver bonus, whereas blank or very light zones correspond to zero or low bonuses. In Figure~\ref{M1}, high bonuses are mainly concentrated along the Manhattan street grid, especially in central and lower Manhattan. In Figure~\ref{E1}, high-bonus zones appear across Manhattan and extend into nearby surrounding areas, reflecting a more spatially dispersed evening supply--demand imbalance.

\subsection{Stochastic Driver Response Evaluation}
\label{exp_stochastic}

In practice, driver choices can cause realized supply and net operating profit to differ from their expected-supply estimates. We therefore conduct 300 independent Monte Carlo (MC) realizations for each fixed bonus map in Table~\ref{case_study_benchmark} on the nine New York City instances. In each realization, drivers respond according to the choice probabilities in Section~\ref{subsec:mathematical_formulation}, and DSMDP-M is re-solved with realized supply. 
Motivated by \citet{ashkrof2024relocation}, we also use matching-aware Monte Carlo (MA-MC) simulation to capture how drivers adjust their choices based on their perceived likelihood of receiving the advertised bonus.
Specifically, each \(b_i\) in the choice utilities is replaced by \(b_i\min\{1,q_i/s_i\}\) for \(s_i>0\), and zero otherwise, where \(q_i\) and \(s_i\) are the expected matched quantity and supply in zone \(i\). 
Table~\ref{tab:stochastic_response} compares the original, mean MC, and mean MA-MC objectives across methods and instances.
The gray percentages report changes relative to
Original for MC and relative to MC for MA-MC.
NOMAD is omitted because it returns the same map as \(\mathbf H_1\).

\begin{table}[htbp]
\centering
\caption{Performance under stochastic driver responses in the New York City study.}
\label{tab:stochastic_response}
{\scriptsize
\setlength{\tabcolsep}{2pt}
\renewcommand{\arraystretch}{0.86}
\newcommand{\modechg}[1]{%
  {\color[gray]{0.45}\scriptsize(#1\%)}%
}
\newcommand{\instsep}{%
  \noalign{%
    \vskip 1pt
    \begingroup
      \color[gray]{0.85}
      \hrule height 0.24pt
    \endgroup
    \vskip 1pt
  }%
}
\newsavebox{\stochbox}
\def\stochrows{%
\toprule
ID & Mode & Conic & $\mathbf B^0$ & $\mathbf H_1$
& $\mathbf H_2$ & ISA \\
\midrule
\multirow{3}{*}{M1}
& Original & 5796 & 6576 & 8224 & 6730 & \best{8703} \\
& MC
& 5785\,\modechg{-0.19}
& 6534\,\modechg{-0.64}
& 8174\,\modechg{-0.60}
& 6687\,\modechg{-0.63}
& \best{8647}\,\modechg{-0.64}\,\gain{5.79} \\
& MA-MC
& 5808\,\modechg{+0.39}
& 6534\,\modechg{0.00}
& 8174\,\modechg{0.00}
& 6688\,\modechg{+0.01}
& \best{8648}\,\modechg{+0.01}\,\gain{5.79} \\
\instsep
\multirow{3}{*}{M2}
& Original & 5251 & 6548 & 8276 & 6727 & \best{8974} \\
& MC
& 5226\,\modechg{-0.48}
& 6504\,\modechg{-0.68}
& 8229\,\modechg{-0.57}
& 6672\,\modechg{-0.82}
& \best{8922}\,\modechg{-0.58}\,\gain{8.43} \\
& MA-MC
& 5267\,\modechg{+0.79}
& 6504\,\modechg{0.00}
& 8229\,\modechg{0.00}
& 6672\,\modechg{+0.01}
& \best{8924}\,\modechg{+0.02}\,\gain{8.45} \\
\instsep
\multirow{3}{*}{M3}
& Original & 5107 & 6466 & 8371 & 6643 & \best{8881} \\
& MC
& 5094\,\modechg{-0.26}
& 6422\,\modechg{-0.67}
& 8306\,\modechg{-0.78}
& 6595\,\modechg{-0.72}
& \best{8811}\,\modechg{-0.79}\,\gain{6.08} \\
& MA-MC
& 5132\,\modechg{+0.74}
& 6422\,\modechg{0.00}
& 8306\,\modechg{0.00}
& 6596\,\modechg{+0.01}
& \best{8811}\,\modechg{0.00}\,\gain{6.08} \\
\midrule
\multirow{3}{*}{O1}
& Original & 5284 & 9157 & 11348 & 9354 & \best{11863} \\
& MC
& 5268\,\modechg{-0.30}
& 9119\,\modechg{-0.42}
& 11290\,\modechg{-0.51}
& 9316\,\modechg{-0.41}
& \best{11791}\,\modechg{-0.61}\,\gain{4.43} \\
& MA-MC
& 5292\,\modechg{+0.46}
& 9119\,\modechg{0.00}
& 11291\,\modechg{+0.01}
& 9316\,\modechg{0.00}
& \best{11791}\,\modechg{+0.01}\,\gain{4.43} \\
\instsep
\multirow{3}{*}{O2}
& Original & 5009 & 6994 & 8777 & 7154 & \best{9180} \\
& MC
& 4978\,\modechg{-0.63}
& 6945\,\modechg{-0.70}
& 8724\,\modechg{-0.60}
& 7104\,\modechg{-0.70}
& \best{9113}\,\modechg{-0.73}\,\gain{4.46} \\
& MA-MC
& 5024\,\modechg{+0.93}
& 6945\,\modechg{0.00}
& 8724\,\modechg{0.00}
& 7104\,\modechg{0.00}
& \best{9113}\,\modechg{0.00}\,\gain{4.46} \\
\instsep
\multirow{3}{*}{O3}
& Original & 5781 & 7616 & 9555 & 7771 & \best{9785} \\
& MC
& 5770\,\modechg{-0.20}
& 7589\,\modechg{-0.35}
& 9520\,\modechg{-0.37}
& 7743\,\modechg{-0.35}
& \best{9732}\,\modechg{-0.54}\,\gain{2.23} \\
& MA-MC
& 5786\,\modechg{+0.29}
& 7589\,\modechg{0.00}
& 9520\,\modechg{0.00}
& 7744\,\modechg{0.00}
& \best{9733}\,\modechg{+0.01}\,\gain{2.24} \\
\midrule
\multirow{3}{*}{E1}
& Original & 5284 & 8191 & 10479 & 8383 & \best{10829} \\
& MC
& 5271\,\modechg{-0.24}
& 8136\,\modechg{-0.67}
& 10431\,\modechg{-0.46}
& 8329\,\modechg{-0.64}
& \best{10783}\,\modechg{-0.42}\,\gain{3.38} \\
& MA-MC
& 5294\,\modechg{+0.44}
& 8136\,\modechg{0.00}
& 10431\,\modechg{0.00}
& 8330\,\modechg{+0.01}
& \best{10785}\,\modechg{+0.01}\,\gain{3.39} \\
\instsep
\multirow{3}{*}{E2}
& Original & 5095 & 7760 & 10015 & 7941 & \best{10253} \\
& MC
& 5081\,\modechg{-0.27}
& 7739\,\modechg{-0.28}
& 9974\,\modechg{-0.41}
& 7917\,\modechg{-0.30}
& \best{10186}\,\modechg{-0.66}\,\gain{2.13} \\
& MA-MC
& 5113\,\modechg{+0.62}
& 7739\,\modechg{0.00}
& 9974\,\modechg{0.00}
& 7917\,\modechg{0.00}
& \best{10186}\,\modechg{0.00}\,\gain{2.13} \\
\instsep
\multirow{3}{*}{E3}
& Original & 5231 & 8969 & 11241 & 9165 & \best{11715} \\
& MC
& 5215\,\modechg{-0.30}
& 8921\,\modechg{-0.53}
& 11196\,\modechg{-0.40}
& 9120\,\modechg{-0.49}
& \best{11666}\,\modechg{-0.42}\,\gain{4.20} \\
& MA-MC
& 5274\,\modechg{+1.13}
& 8921\,\modechg{0.00}
& 11197\,\modechg{+0.01}
& 9121\,\modechg{+0.01}
& \best{11668}\,\modechg{+0.02}\,\gain{4.21} \\
\bottomrule
}
\sbox{\stochbox}{%
  \begin{tabular}{@{}c c*{5}{c}@{}}
    \stochrows
  \end{tabular}%
}
\resizebox{\linewidth}{!}{%
  \begin{tabular*}{\dimexpr\wd\stochbox*5/4\relax}
    {@{\extracolsep{\fill}}c c*{5}{c}@{}}
    \stochrows
  \end{tabular*}%
}
}
\compacttableend
\end{table}

Table~\ref{tab:stochastic_response} shows that stochastic driver responses reduce realized net operating profit by less than 1\%, while ISA remains the best-performing method among the original benchmarks in all nine instances under both MC and MA-MC. 
The results also show that accounting for drivers' perceived likelihood of receiving bonuses has little effect on profit. The adjustment primarily affects zones with unmatched supply, since additional drivers have limited marginal matching value.
These findings demonstrate the robustness of surge maps optimized using expected supply, even in the presence of stochastic and matching-aware driver responses.
Finally, changes in bonuses can shift the optimal matching basis, making current-map marginal values less reliable for guiding updates. A comparison with the center-gradient specialization of ISA (ISA-CG) suggests that sampling nearby matching regimes is particularly beneficial under tighter bonus budgets (\ECsection{EC.6.4}).


\section{Conclusion}\label{sec_conclusion}
This paper studies how ride-hailing platforms can design driver surge maps in real time.
We formulate the DSMDP to optimize zone-level bonuses while anticipating voluntary driver relocation and participation and the resulting downstream matching value.
Exploiting the linear-program structure of the matching
problem, ISA samples dual sensitivities across nearby matching regimes to guide feasible bonus-map updates.
Theoretical analysis provides finite-scale stationarity
certificates and identifies sufficient conditions for
constrained Clarke stationarity of sampled limits.
A conic outer relaxation provides a valid upper bound for evaluating feasible solutions.

The computational results show that ISA improves net operating profit over the strongest non-ISA feasible benchmark by up to 19.37\% on the generated instances and by 2.38\%--8.43\% in the New York City study.
The marginal analysis further suggests that platforms should direct bonuses toward responsive supply whose downstream matching value justifies the associated payments, rather than relying on local shortage alone.
Intervention generates the greatest incremental value under supply shortage.
Higher bonuses increase payments on existing matches and can also incur a budget opportunity cost when the bonus budget is binding, leading to diminishing returns to additional spending.
Latent-driver activation can remain valuable even when all demand is served by supporting higher-value matches.
Moderate spatial smoothness preserves most of the gains from spatial targeting, suggesting that effective bonus maps need not display sharp differences between neighboring zones.
Finally, under stochastic driver responses, the expected-supply objective overstates mean realized net operating profit by less than 1\% under MC, while ISA's improvements over the original benchmarks persist under both MC and MA-MC, suggesting that expected-supply optimization provides a robust basis for real-time surge map design.

We conclude with several directions for extending the current framework.
The current analysis is based on a single-period model, pre-estimated behavioral-response parameters, and a platform net operating profit objective.
Future research could extend the framework to dynamic settings in which current bonuses affect future supply and matching opportunities, incorporate heterogeneous
and strategic driver behavior such as multi-homing and trip acceptance, and develop robust or data-driven methods for uncertain demand, driver responses, and real-time
supply observations.
Extensions could also examine how service quality, driver earnings, and spatial equity alter the design of driver surge maps.

\bibliographystyle{informs2014}
\bibliography{mybibfile}

\ECRUNAUTHOR{Xu et al.}
\ECRUNTITLE{Driver Surge Maps with Downstream Matching}
\ECSwitch
\ECHead{E-Companion to Driver Surge Maps for Supply Rebalancing with Downstream Matching}
\makeatletter
\@ifundefined{theHsection}{}{\renewcommand{\theHsection}{ec.\arabic{section}}}
\@ifundefined{theHsubsection}{}{\renewcommand{\theHsubsection}{ec.\arabic{section}.\arabic{subsection}}}
\@ifundefined{theHsubsubsection}{}{\renewcommand{\theHsubsubsection}{ec.\arabic{section}.\arabic{subsection}.\arabic{subsubsection}}}
\@ifundefined{theHequation}{}{\renewcommand{\theHequation}{ec.\arabic{equation}}}
\@ifundefined{theHtheorem}{}{\renewcommand{\theHtheorem}{ec.\arabic{theorem}}}
\@ifundefined{theHlemma}{}{\renewcommand{\theHlemma}{ec.\arabic{lemma}}}
\@ifundefined{theHproposition}{}{\renewcommand{\theHproposition}{ec.\arabic{proposition}}}
\@ifundefined{theHcorollary}{}{\renewcommand{\theHcorollary}{ec.\arabic{corollary}}}
\@ifundefined{theHfigure}{}{\renewcommand{\theHfigure}{ec.\arabic{figure}}}
\@ifundefined{theHtable}{}{\renewcommand{\theHtable}{ec.\arabic{table}}}
\@ifundefined{theHalgorithm}{}{\renewcommand{\theHalgorithm}{ec.\arabic{algorithm}}}
\makeatother
\setcounter{theorem}{0}
\setcounter{proposition}{0}
\setcounter{corollary}{0}
\begingroup
\renewcommand{\compactequationspacing}{\compactalignspacing}
\setlength{\jot}{0pt}

\section{Calibration of Driver Response Parameters}\label{cal_para}

In deployment, the behavioral-response parameters introduced in the mathematical formulation can be estimated from historical data before optimization and treated as fixed inputs when solving the DSMDP.
Specifically, the platform can construct a large sample of driver transitions between consecutive decision periods using historical platform records.
For each dispatchable driver, the data record the current zone, the bonuses displayed in nearby reachable zones, and the location or dispatch-availability state in the next period.
Transitions associated with dispatched trips can be excluded to isolate voluntary responses.

For observation \(n\), let \(i_n\) denote the current zone, \(\bb_n\) the displayed bonus vector, and \(y_n\) the observed next-period choice, which is either a reachable zone \(j\in\mathcal Z_{i_n}\) or the unavailable state.
The parameters \(\bbeta\), \(\boldsymbol{\delta}^{\mathrm{Off}}\), and \(\zeta^{\mathrm{Off}}\) can then be estimated by maximizing the multinomial logit (MNL) log-likelihood \(\sum_n\left[\sum_{j\in\mathcal Z_{i_n}}\mathbbm{1}_{\{y_n=j\}}\log P_{i_nj}(\bb_n)+\mathbbm{1}_{\{y_n=\mathrm{Off}\}}\log P_{i_n}^{\mathrm{Off}}(\bb_n)\right]\) subject to \(\zeta^{\mathrm{Off}}\ge0\).
The likelihood is expressed in terms of the normalized utilities
\(V_{ij}=\beta^{(1)}b_j+\beta^{(3)}d_{ij}\) and
\(V_i^{\mathrm{Off}}=\delta_i^{\mathrm{Off}}-\zeta^{\mathrm{Off}}b_i\),
obtained by subtracting the common dispatchable-state utility term
\(\beta^{(0)}+\beta^{(2)}b_i\).
The normalization yields
\(\delta_i^{\mathrm{Off}}=\alpha_i^{\mathrm{Off}}-\beta^{(0)}\) and
\(\zeta^{\mathrm{Off}}=\beta^{(2)}+\gamma^{\mathrm{Off}}\), allowing the identified coefficient combinations to be estimated
directly while preserving the original choice probabilities.

Similarly, for latent drivers, the platform observes their current zone, local bonus, and whether they enter the dispatchable pool in the next period.
Let \(a_n=1\) if observation \(n\) records activation into the dispatchable pool and 0 otherwise.
The parameters
\(\boldsymbol{\alpha}^{\mathrm{On}}\) and \(\gamma^{\mathrm{On}}\)
can be estimated using the binary logit likelihood
\(
\sum_n
\left[
a_n\log P_{i_n}^{\mathrm{On}}(\bb_n)
+
(1-a_n)\log
\left(1-P_{i_n}^{\mathrm{On}}(\bb_n)\right)
\right].
\)
The estimated parameters may be updated periodically using a rolling historical sample and are held fixed during each real-time decision period.

\section{Derivatives of Driver Response Probabilities}
\label{app:response_derivatives}

The utility derivatives used by the Lagrangian-based gradient update are
\begin{equation}
\frac{\partial V_{ij}}{\partial b_m}
=
\beta^{(1)}\mathbbm{1}_{\{j=m\}},
\qquad
\frac{\partial V_i^{\mathrm{Off}}}{\partial b_m}
=
-\zeta^{\mathrm{Off}}\mathbbm{1}_{\{i=m\}},
\qquad
i\in\mathcal Z,\ j\in\mathcal Z_i,\ m\in\mathcal Z.
\label{eq:app_utility_derivatives}
\end{equation}
Substituting the utility derivatives in Eq.~\eqref{eq:app_utility_derivatives} into the MNL derivative gives, for \(i,m\in\mathcal Z\) and \(j\in\mathcal Z_i\),
\begin{equation}
\frac{\partial P_{ij}(\bb)}{\partial b_m}
=
P_{ij}(\bb)\!\left[
\beta^{(1)}
\left(
\mathbbm{1}_{\{j=m\}}
-\mathbbm{1}_{\{m\in\mathcal Z_i\}}P_{im}(\bb)
\right)
+\zeta^{\mathrm{Off}}P_i^{\mathrm{Off}}(\bb)
\mathbbm{1}_{\{i=m\}}
\right].
\label{eq:dPij_mnl_spec}
\end{equation}
The parameter \(\zeta^{\mathrm{Off}}\) represents the identified relative effect of the origin-zone bonus on the dispatchable-versus-unavailable margin.
The derivative of the unavailability probability is
\begin{equation}
\frac{\partial P_i^{\mathrm{Off}}(\bb)}{\partial b_m}
=
-P_i^{\mathrm{Off}}(\bb)
\left[
\beta^{(1)}
\mathbbm{1}_{\{m\in\mathcal Z_i\}}P_{im}(\bb)
+\zeta^{\mathrm{Off}}
\left(1-P_i^{\mathrm{Off}}(\bb)\right)
\mathbbm{1}_{\{i=m\}}
\right],
\qquad i,m\in\mathcal Z.
\label{eq:dPoff_closed}
\end{equation}
The derivative of the latent-driver activation probability is
\begin{equation}
\frac{\partial P_i^{\mathrm{On}}(\bb)}{\partial b_m}
=
\gamma^{\mathrm{On}}
P_i^{\mathrm{On}}(\bb)
\left(1-P_i^{\mathrm{On}}(\bb)\right)
\mathbbm{1}_{\{i=m\}},
\qquad \forall i,m\in\mathcal Z.
\label{eq:dPon_closed}
\end{equation}

Combining Eqs.~\eqref{eq:dPij_mnl_spec}, \eqref{eq:dPoff_closed}, and \eqref{eq:dPon_closed}, we obtain the supply derivative
\begin{equation}
\begin{aligned}
\frac{\partial s_i(\bb)}{\partial b_m}
={}&
\sum_{j\in\mathcal Z_i\setminus\{i\}}
s_j^{\mathrm{On}}
\frac{\partial P_{ji}(\bb)}{\partial b_m}
\\
&-
s_i^{\mathrm{On}}
\left[
\sum_{j\in\mathcal Z_i\setminus\{i\}}
\frac{\partial P_{ij}(\bb)}{\partial b_m}
+
\frac{\partial P_i^{\mathrm{Off}}(\bb)}{\partial b_m}
\right]
\\
&+
s_i^{\mathrm{Off}}
\frac{\partial P_i^{\mathrm{On}}(\bb)}{\partial b_m},
\qquad \forall i,m\in\mathcal Z.
\end{aligned}
\label{eq:dsdb_generalized}
\end{equation}
\subsection{Derivation of the Marginal Bonus Value}

The supply sensitivity in Eq.~\eqref{eq:dsdb_generalized} enters the Lagrangian derivative below.
Holding the selected optimal primal-dual solution fixed, differentiating the Lagrangian of Objective~\eqref{eq:matchingObjective} and Constraints~\eqref{eq:matchingDemand}--\eqref{eq:matchingBudget} with multipliers \((\bpi,\boldsymbol{\eta},\lambda)\) gives Eq.~\eqref{eq:grad_phi}; the direct payout and budget derivatives combine as \(-(1+\lambda)q_m\).
\begin{proof}{Proof of Proposition~\ref{prop:matching_weighted_bonus_value}.}
Given
\(\sum_{j\in\mathcal Z_i}P_{ij}(\bb)
+
P_i^{\mathrm{Off}}(\bb)
=
1\),
the supply-balance equation can be written as
\begin{equation}
s_i(\bb)
=
\sum_{j:\,i\in\mathcal Z_j}
s_j^{\mathrm{On}}P_{ji}(\bb)
+
s_i^{\mathrm{Off}}P_i^{\mathrm{On}}(\bb).
\label{eq:ec_supply_compact}
\end{equation}
Using \(\bar\eta_i\) as defined before Proposition~\ref{prop:matching_weighted_bonus_value}, Eq.~\eqref{eq:dPij_mnl_spec} gives Eq.~\eqref{eq:weighted_dispatchable_response}.
Differentiating Eq.~\eqref{eq:ec_supply_compact}, weighting the result by the supply-constraint dual multipliers, and applying Eqs.~\eqref{eq:weighted_dispatchable_response} and \eqref{eq:dPon_closed} therefore yields
\begin{align}
\sum_{j\in\mathcal Z}
\eta_j\frac{\partial s_j(\bb)}{\partial b_m}
={}&
\beta^{(1)}
\sum_{i:\,m\in\mathcal Z_i}
s_i^{\mathrm{On}}P_{im}(\bb)
\bigl(\eta_m-\bar\eta_i\bigr)
\nonumber\\*
&+
\zeta^{\mathrm{Off}}
s_m^{\mathrm{On}}P_m^{\mathrm{Off}}(\bb)\bar\eta_m
\nonumber\\*
&+
\gamma^{\mathrm{On}}s_m^{\mathrm{Off}}
P_m^{\mathrm{On}}(\bb)
\bigl(1-P_m^{\mathrm{On}}(\bb)\bigr)\eta_m.
\label{eq:ec_matching_weighted_supply}
\end{align}
Substituting Eq.~\eqref{eq:ec_matching_weighted_supply} into Eq.~\eqref{eq:grad_phi} proves
Eq.~\eqref{eq:matching_weighted_marginal}.
\Halmos
\end{proof}

\subsection{Proofs of the Intervention and Spatial Balance Results}

\begin{proof}{Proof of Corollary~\ref{cor:no_bonus_intervention}.}
At the no-bonus map, realized bonus spending is zero and the budget constraint is strictly slack because \(B>0\).
Complementary slackness therefore gives \(\lambda^{\mathrm{no}}=0\).
Differentiability of \(Z\) at \(\boldsymbol 0\) and Proposition~\ref{prop:matching_weighted_bonus_value} then give Eq.~\eqref{eq:no_bonus_intervention_threshold}.
Moreover, \(b_i\boldsymbol e_i\in\mathcal B\) for every sufficiently small \(b_i>0\) because \(b^{\max}>0\) and \(\Delta^b>0\).
A positive one-sided derivative consequently yields
\(Z(b_i\boldsymbol e_i)>Z(\boldsymbol 0)\) for every
sufficiently small positive \(b_i\).

At \(\boldsymbol 0\), the tangent cone of \(\mathcal B\) is the nonnegative orthant.
The first-order stationarity condition is therefore equivalent to \(g_i^{\mathrm{no}}\le0\) for every \(i\in\mathcal Z\).
If \(g_i^{\mathrm{no}}<0\) for every \(i\in\mathcal Z\), differentiability at \(\boldsymbol 0\) gives
\[
Z(\bb)-Z(\boldsymbol 0)
=
(\bg^{\mathrm{no}})^\top\bb+o(\|\bb\|)<0
\]
for every sufficiently small nonzero feasible \(\bb\), since \(\bb\ge\boldsymbol 0\).
Thus, the no-bonus map is a strict local maximizer.
\Halmos
\end{proof}

\begin{proof}{Derivation of Eq.~\eqref{eq:explicit_activation_threshold}.}
In the single-zone setting,
Eq.~\eqref{eq:one_zone_objective} is the reduced objective on a right
neighborhood of zero. 
The right derivative of \(Z_i(b_i)\) at zero,
denoted by \(Z_{i,+}'(0)\), is
\[
Z_{i,+}'(0)
=
p_{ii}s_i^{\mathrm{Off}}\gamma^{\mathrm{On}}
P_i^{\mathrm{On},\mathrm{no}}
\bigl(1-P_i^{\mathrm{On},\mathrm{no}}\bigr)
-
\left(
s_i^{\mathrm{On}}
+
s_i^{\mathrm{Off}}P_i^{\mathrm{On},\mathrm{no}}
\right).
\]
Rearranging the condition \(Z_{i,+}'(0)>0\) yields
Eq.~\eqref{eq:explicit_activation_threshold}.
\Halmos
\end{proof}

\begin{proof}{Proof of Proposition~\ref{prop:spatial_marginal_balance} and Corollary~\ref{cor:spatial_component_balance}.}
Since \(\bg^\star=\nabla Z(\bb^\star)\),
local optimality over the convex set \(\mathcal B\) implies
\begin{equation}
(\bg^\star)^\top(\bb-\bb^\star)
\le0,
\qquad
\forall\bb\in\mathcal B,
\label{eq:ec_spatial_variational_inequality}
\end{equation}
The variational inequality in Eq.~\eqref{eq:ec_spatial_variational_inequality}
gives
\(\bg^\star\in\mathcal N_{\mathcal B}(\bb^\star)\).
The normal cone of the polyhedral set in
Eq.~\eqref{eq:bonus_feasible_set} is the nonnegative hull of its
active outward normals, which yields
Eq.~\eqref{eq:spatial_kkt_balance}.

For an interior zone with no incident spatial multiplier, every term on the right-hand side of the corresponding coordinate in Eq.~\eqref{eq:spatial_kkt_balance} is zero.
Equation~\eqref{eq:interior_zone_balance} then follows from Eq.~\eqref{eq:matching_weighted_marginal}.
For a component \(C\), every edge outside \(\mathcal E_\rho\) has zero multiplier, and every edge in \(\mathcal E_\rho\) has both endpoints in the same component.
Let \(\boldsymbol\rho^\pm:=(\rho_e^\pm)_{e\in\mathcal E^b}\).
Summing the coordinates of \(\boldsymbol D^\top(\boldsymbol\rho^+-\boldsymbol\rho^-)\) over \(C\) gives zero because each incidence row contributes once with sign \(+1\) and once with sign \(-1\).
Equations~\eqref{eq:component_balance_with_bounds} and \eqref{eq:active_component_balance} follow.
\Halmos
\end{proof}

\section{Adaptive Constrained Gradient Sampling}
\label{ec:armijo_stationarity}

Following the gradient-sampling principle
\citep{burke2002approximating,curtis2012sequential}, ISA samples sensitivities from nearby matching branches to address transitions between locally smooth regimes.
We specify the sampling, null-step, scale-refinement, and stopping procedures omitted from Algorithm~\ref{alg:AO}.
The construction links the sampled sensitivity bundle to feasible ascent, stationarity assessment, and adaptive scale refinement.

\subsection{Sampled Constrained Direction}
\label{app:sec:3.1}

At the center \(\bb^k\in\mathcal B\), the sampling radius
\(\epsilon_k\), stationarity scale \(\nu_k\), and sample pool may
be updated over several attempts while the center remains fixed,
with all existing evaluations reused. 
Define
\(\mathcal I_k:=\{0\}\cup\mathcal A_k\), where \(a=0\) denotes the
evaluation at the center and \(\mathcal A_k\) indexes the ambient samples.
Set \(\by^{k,0}:=\bb^k\) and draw each \(\by^{k,a}\),
\(a\in\mathcal A_k\), uniformly from
\(\mathbb B_{\epsilon_k}(\bb^k):=
\{\by:\|\by-\bb^k\|_2\le\epsilon_k\}\).

Suppose \(Z\) admits the locally Lipschitz extension in
Propositions~\ref{prop:finite_scale_certificate} and \ref{prop:sampled_clarke_limit}.
The ambient samples are then differentiability points with probability one whenever their sampling balls lie within this neighborhood.
We assume that, at each ambient sample \(\by^{k,a}\), the sensitivity computed from Eq.~\eqref{eq:grad_phi} equals \(\nabla Z(\by^{k,a})\), so that \(\bg^{k,a}=\nabla Z(\by^{k,a})\) for each \(a\in\mathcal A_k\).
The center sensitivity \(\bg^{k,0}\) is the sensitivity of the selected center branch and is retained in the working bundle even when differentiability at the center has not been established.
Let \(\bG_k:=[\bg^{k,a}]_{a\in\mathcal I_k}\) denote this working bundle.

The constraints defining \(\mathcal B\) enter through \(\bN_k\), while the realized bonus-spending budget enters each sampled sensitivity through the multiplier \(\lambda\) in Eq.~\eqref{eq:grad_phi}. 
ISA then solves
\begin{gather}
(\bomega^k,\bmu^k)
\in
\arg\min_{\substack{
\bomega\ge\boldsymbol{0},\;
\boldsymbol{1}^\top\bomega=1\\
\bmu\ge\boldsymbol{0}
}}
\frac12
\left\|
\bG_k\bomega-\bN_k\bmu
\right\|_2^2
\label{ec:eq:sampled_direction_qp}
\\
\bar{\bg}^k
:=\bG_k\bomega^k,\qquad
\bar{\bn}^k:=\bN_k\bmu^k,\qquad
\br^k:=\bar{\bg}^k-\bar{\bn}^k,\qquad
\widehat R_k:=\|\br^k\|_2.
\label{ec:eq:sampled_direction}
\end{gather}
The QP optimality conditions imply
\(\bN_k^\top\br^k\le\boldsymbol{0}\) and
\(\min_{a\in\mathcal I_k}
(\bg^{k,a})^\top\br^k\ge\|\br^k\|_2^2\).
Hence, \(\br^k\) is a feasible common ascent direction for the sampled
sensitivities whenever \(\widehat R_k>0\).
The projected line search below tests this direction while preserving the bonus-feasible set \(\mathcal B\).

\subsection{Projected Bundle Line Search and Null Steps}
\label{appen:sec:3.2}

Fix an initial step \(\bar\tau>0\), contraction factor
\(\delta\in(0,1)\), Armijo sufficient-increase parameter \(\chi\in(0,1)\), maximum
number of trials \(L_k\), and minimum trial step
\(\tau_{\min}>0\). For \(\ell=0,\ldots,L_k-1\), set
\(\tau_{k,\ell}:=\bar\tau\delta^\ell\),
\(\bh^{k,\ell}:=
\Pi_{\mathcal B}(\bb^k+\tau_{k,\ell}\br^k)-\bb^k\), and
\(\Gamma_k(\bh^{k,\ell})
:=\min_{a\in\mathcal I_k}
(\bg^{k,a})^\top\bh^{k,\ell}\).
To avoid vanishingly small steps, trials with step sizes below \(\tau_{\min}\) are not evaluated.
A trial is accepted if
\(\Gamma_k(\bh^{k,\ell})>0\) and
\(Z(\bb^k+\bh^{k,\ell})
\ge Z(\bb^k)+\chi\Gamma_k(\bh^{k,\ell})\).

An accepted trial sets
\(\bb^{k+1}:=\bb^k+\bh^{k,\ell}\) and resets the sample pool around
the new center. If no trial is accepted, the center is retained and a
batch of new ambient evaluations is added. Once the pool reaches its
capacity, ambient sensitivities with the smallest QP weights are
removed, while the center sensitivity is retained. When a time,
accepted-step, or sampling-attempt limit is reached, ISA
returns the best center found.
The working-bundle residual \(\widehat R_k\) and accumulated ambient samples feed into the scale-refinement and stationarity assessment described next.

\subsection{Adaptive Sampling Radius and Stationarity Assessment}
\label{appen:sec:3.3}

ISA uses an attempt
window \(W\) and an objective-stability threshold \(\eta_Z>0\). Let
\(Z^{[u]}\) denote the objective value at the center for target-scale
attempt \(u\). Once \(W+1\) observations are available, define
\(\Delta_Z^{[u]}:=
\max\{0,Z^{[u]}-Z^{[u-W]}\}/
\max\{1,|Z^{[u]}|\}\). Objective stability holds when
\(\Delta_Z^{[u]}\le\eta_Z\).

When \(\widehat R_k\le\nu_k\) and either the current scales are coarse or the
target-scale objective-stability test is satisfied, ISA enriches the
ambient pool to contain at least \(|\mathcal Z|+1\) points. It then
re-solves Eq.~\eqref{ec:eq:sampled_direction_qp} using only
\(\bG_k^{\mathrm{amb}}:=
[\bg^{k,a}]_{a\in\mathcal A_k}\), thereby omitting the deterministic
center sensitivity \(\bg^{k,0}\).
Let \(\widehat R_k^{\mathrm{amb}}\) and
\(\br^{k,\mathrm{amb}}\) denote the resulting residual and direction,
respectively. At the target scales, this ambient assessment is
deferred until objective stability holds; earlier attempts continue
with the working bundle \(\bG_k\).
Proposition~\ref{prop:finite_scale_certificate} shows that an ambient residual \(\widehat R_k^{\mathrm{amb}}\le\nu_k\), paired with radius \(\epsilon_k\), gives an \((\epsilon_k,\nu_k)\)-stationarity certificate.

Scale refinement or termination occurs only when
\(\widehat R_k^{\mathrm{amb}}\le\nu_k\). If this condition holds at
\((\epsilon_{\mathrm{tar}},\nu_{\mathrm{tar}})\), ISA returns the
current center. At a coarser scale, it sets
\(\epsilon_k\leftarrow
\max\{\epsilon_{\mathrm{tar}},\kappa_\epsilon\epsilon_k\}\) and
\(\nu_k\leftarrow
\max\{\nu_{\mathrm{tar}},\kappa_\nu\nu_k\}\), where
\(\kappa_\epsilon,\kappa_\nu\in(0,1)\), and rebuilds the pool at the
same center. If \(\widehat R_k^{\mathrm{amb}}>\nu_k\), the subsequent line search uses \(\br^{k,\mathrm{amb}}\) and the
ambient bundle \(\bG_k^{\mathrm{amb}}\). 
For this line search, replace \(\br^k\) by
\(\br^{k,\mathrm{amb}}\) and \(\mathcal I_k\) by
\(\mathcal A_k\) in the line-search formulas.
Proposition~\ref{prop:sampled_clarke_limit} establishes
constrained Clarke stationarity along convergent ambient assessments satisfying
\(\epsilon_\ell\downarrow0\) and \(\widehat R_\ell^{\mathrm{amb}}\to0\).

\subsection{Smooth Projected Gradient Specialization}
\label{appen:sec:3.4}

For the smooth-case analysis, we consider
a specialization that uses only the center gradient.
Under the stable optimal-basis condition in Section~4.1,
the Lagrangian sensitivity equals \(\nabla Z(\bb)\).
Fix an initial step \(\bar\tau>0\), a contraction factor
\(\delta\in(0,1)\), and an Armijo parameter \(\chi\in(0,1)\).
For \(\tau_{k,\ell}:=\bar\tau\delta^\ell\),
\(\ell=0,1,\ldots\), define
\(
\bh_{k,\ell}
:=
\Pi_{\mathcal B}
\bigl(\bb^k+\tau_{k,\ell}\nabla Z(\bb^k)\bigr)-\bb^k.
\)
The procedure accepts the first trial satisfying
\(
Z(\bb^k+\bh_{k,\ell})
\ge
Z(\bb^k)
+\chi\nabla Z(\bb^k)^\top \bh_{k,\ell},
\)
and sets \(\bb^{k+1}:=\bb^k+\bh_{k,\ell}\).
The convergence and accepted-step complexity of this
specialization are established in \ECsection{EC.4.4}.

\section{Proofs of Theoretical Results}
\label{app:proofs}

The proofs use the normal-cone geometry of \(\mathcal B\) and parametric-LP sensitivity to connect the sampled directions to stationarity.
We establish continuity and the finite-scale certificate before analyzing smooth accepted steps and shrinking-scale Clarke limits.

\subsection{Proof of Lemma~\ref{lemma:continuity_existence}}
\label{app:sec:proof_lemma1}

\begin{proof}{Proof.}
Let \(\mathcal X(\bb)\) denote the feasible set of DSMDP-M, and define
\(q_i(\bx):=\sum_{j\in\mathcal Z}x_{ij}\) and
\(\bq(\bx):=(q_i(\bx))_{i\in\mathcal Z}\).
The zero matching is feasible for every \(\bb\), and the demand
constraints provide the uniform bound
\(
\sum_{i\in\mathcal Z}\sum_{j\in\mathcal Z}x_{ij}
\le
\sum_{j\in\mathcal Z}d_j.
\)
Thus, \(\mathcal X(\bb)\) is nonempty and compact and lies within
a common compact set.
Because \(\bs(\bb)\) is continuous and every defining inequality is
continuous in \((\bb,\bx)\), \(\mathcal X(\cdot)\) has a closed graph
and is upper hemicontinuous.

It remains to establish lower hemicontinuity.
Let \(\bb^\nu\in\mathcal B\) converge to \(\bb\), and take
any \(\bx\in\mathcal X(\bb)\).
For each \(i\), let \(a_i^\nu\) be the fraction of its matching
flows retained to satisfy the supply limit at \(\bb^\nu\). 
Set \(a_i^\nu:=\min\{1,s_i(\bb^\nu)/q_i(\bx)\}\) if \(q_i(\bx)>0\)
and \(a_i^\nu:=1\) otherwise, and define
\(\widetilde x_{ij}^\nu:=a_i^\nu x_{ij}\).
If \(q_i(\bx)>0\), feasibility implies
\(s_i(\bb)\ge q_i(\bx)>0\), and the continuity of \(s_i\) gives
\(a_i^\nu\to1\).
The same conclusion follows immediately when \(q_i(\bx)=0\).
Hence, \(\widetilde{\bx}^\nu\to\bx\).
Row scaling preserves the demand constraints and ensures
\(q_i(\widetilde{\bx}^\nu)\le s_i(\bb^\nu)\).

To restore the parameter-dependent budget constraint, adopt the
convention \(B/0=+\infty\) and define
\(
\theta_\nu
:=
\min\left\{
1,\,
\frac{B}
{(\bb^\nu)^\top\bq(\widetilde{\bx}^\nu)}
\right\},
\bx^\nu:=\theta_\nu\widetilde{\bx}^\nu.
\)
The denominator converges to
\(\bb^\top\bq(\bx)\le B\).
Since \(B>0\), it follows that \(\theta_\nu\to1\).
Therefore, \(\bx^\nu\in\mathcal X(\bb^\nu)\) and
\(\bx^\nu\to\bx\), which proves lower hemicontinuity.
The matching objective is continuous in \((\bb,\bx)\).
Berge's maximum theorem \citep{berge1963topological} therefore implies
that \(Z(\bb)\) is continuous on \(\mathcal B\).
Since \(\mathcal B\) is compact, existence follows from the
Weierstrass theorem \citep{rudin1976principles}.
\Halmos
\end{proof}

\subsection{Proof of Proposition~\ref{prop:finite_scale_certificate}}
\label{app:sec:proof_finite_scale_certificate}

\begin{proof}{Proof.}
Define the convex hull
\begin{equation}
\mathcal H_k
:=
\operatorname{conv}
\left\{
\nabla Z(\by^{k,a}):a\in\mathcal A_k
\right\}.
\label{eq:ec_ambient_gradient_hull}
\end{equation}
Every ambient gradient in Eq.~\eqref{eq:ec_ambient_gradient_hull} belongs to \(\partial^C Z(\by^{k,a})\), and every ambient point lies within distance \(\epsilon_k\) of \(\bb^k\).
It follows from Eq.~\eqref{eq:epsilon_clarke_gradient} that \(\mathcal H_k\subseteq\partial_{\epsilon_k}^C Z(\bb^k)\).
Because the columns of \(\bN_k\) are exactly the active outward normals, up to duplicate generators, \(\{\bN_k\bmu:\bmu\ge\boldsymbol 0\}=\mathcal N_{\mathcal B}(\bb^k)\).
Consequently, Eq.~\eqref{eq:ambient_stationarity_residual} equals \(\operatorname{dist}(\mathcal H_k,\mathcal N_{\mathcal B}(\bb^k))\).
Expanding the first set from \(\mathcal H_k\) to \(\partial_{\epsilon_k}^C Z(\bb^k)\) can only decrease the set distance, which proves Eq.~\eqref{eq:finite_scale_certificate}.
\Halmos
\end{proof}

\subsection{Sampled Feasible Ascent Direction}
\label{app:sec:sampled_ascent_direction}
The finite-scale certificate uses only ambient gradients.
The next lemma returns to the working bundle and establishes the geometry of its search direction.

\begin{lemma}[Sampled feasible ascent direction]
\label{app:lemma:sampled_ascent_direction}
Let \(\bb^k\in\mathcal B\), and let
\((\bomega^k,\bmu^k)\) solve
Eq.~\eqref{ec:eq:sampled_direction_qp}. Then the residual direction defined in
Eq.~\eqref{ec:eq:sampled_direction} satisfies
\begin{gather}
\begin{aligned}
\bN_k^\top\br^k
&\le\boldsymbol{0},
&
(\bmu^k)^\top\bN_k^\top\br^k
&=0,
\end{aligned}
\label{app:eq:sampled_tangent_conditions}
\\
\min_{a\in\mathcal I_k}
(\bg^{k,a})^\top\br^k
\ge
\|\br^k\|_2^2.
\label{app:eq:sampled_common_ascent}
\end{gather}
Consequently, \(\br^k\) belongs to the tangent cone
\(T_{\mathcal B}(\bb^k)\). If
\(\br^k\ne\boldsymbol{0}\), every sensitivity in the working bundle has
strictly positive first-order gain along \(\br^k\).
\end{lemma}

\begin{proof}{Proof.}
The optimality conditions with respect to \(\bmu\) give
\(\bN_k^\top\br^k\le\boldsymbol{0}\), and complementary slackness gives the
second relation in Eq.~\eqref{app:eq:sampled_tangent_conditions}. 
The simplex optimality conditions for \(\bomega\) give
\(
(\bg^{k,a})^\top\br^k
\ge
(\bar{\bg}^k)^\top\br^k
\)
for all \(a\in\mathcal I_k\).
Moreover,
\(
(\bar{\bg}^k)^\top\br^k
=
(\br^k+\bar{\bn}^k)^\top\br^k
=
\|\br^k\|_2^2,
\)
because
\((\bar{\bn}^k)^\top\br^k
=(\bmu^k)^\top \bN_k^\top\br^k=0\).
This proves Eq.~\eqref{app:eq:sampled_common_ascent}.
For the polyhedral set \(\mathcal B\), the tangent cone at \(\bb^k\)
is
\(
T_{\mathcal B}(\bb^k)
=
\{\boldsymbol{v}:\bN_k^\top\boldsymbol{v}\le\boldsymbol{0}\}.
\)
\Halmos
\end{proof}

\subsection{Convergence of the Smooth Specialization}
\label{app:sec:proof_proposition1}

We analyze the center-gradient specialization defined in
\ECsection{EC.3.4}.
Assume the conditions of
Lemma~\ref{lemma:continuity_existence} hold, and let
\[
\mathcal L_0
:=
\{\bb\in\mathcal B:Z(\bb)\ge Z(\bb^0)\}.
\]
Throughout this subsection, assume \(B>0\) and that the
stable optimal-basis condition in Section~4.1 holds on
an open neighborhood of \(\mathcal L_0\).
Let \(\bar\tau>0\), \(\delta\in(0,1)\), and \(\chi\in(0,1)\)
denote the initial trial step, backtracking contraction
factor, and Armijo parameter, respectively.

A point \(\bar{\bb}\in\mathcal B\) is stationary if
\[
\nabla Z(\bar{\bb})^\top(\bb-\bar{\bb})
\le0,
\qquad \forall\bb\in\mathcal B.
\]
For a fixed reference step \(\tau_{\mathrm{ref}}>0\),
define the projected-gradient residual
\begin{equation}
R_{\tau_{\mathrm{ref}}}(\bb)
:=
\frac{1}{\tau_{\mathrm{ref}}}
\left\|
\Pi_{\mathcal B}
\bigl(\bb+\tau_{\mathrm{ref}}\nabla Z(\bb)\bigr)
-\bb
\right\|_2.
\label{app:eq:smooth_residual_definition}
\end{equation}

\begin{proposition}[Convergence of the smooth specialization]
\label{prop:convergence_stationary_final}
The projected Armijo procedure in \ECsection{EC.3.4}
terminates finitely at each iteration.
For the resulting infinite sequence,
\begin{equation}
\lim_{k\to\infty}
R_{\tau_{\mathrm{ref}}}(\bb^k)
=0
\label{app:eq:residual_convergence}
\end{equation}
for every fixed \(\tau_{\mathrm{ref}}>0\),
and every accumulation point is stationary.
\end{proposition}

\begin{proof}{Proof.}
The stable optimal-basis condition and standard
parametric-LP sensitivity analysis
\citep{gal1997advances} imply that \(Z\) is continuously
differentiable with a Lipschitz-continuous gradient on
a closed neighborhood of \(\mathcal L_0\).
The envelope formula identifies the center sensitivity
with \(\nabla Z(\bb)\).

For a trial step \(\tau>0\), define
\(
\bg^k:=\nabla Z(\bb^k),
\bh^k(\tau)
:=
\Pi_{\mathcal B}(\bb^k+\tau\bg^k)-\bb^k.
\)
Projection optimality then gives
\begin{equation}
(\bg^k)^\top\bh^k(\tau)
\ge
\frac{1}{\tau}\|\bh^k(\tau)\|_2^2.
\label{app:eq:projection_ascent}
\end{equation}
Since \(\mathcal L_0\) is compact and contained in the
open stable-basis neighborhood, choose \(r>0\) such that
its closed \(r\)-neighborhood remains within that
neighborhood.
Let \(L\) be a Lipschitz constant of \(\nabla Z\) on this
closed neighborhood, and let
\(
M
:=
\max_{\bb\in\mathcal L_0}
\|\nabla Z(\bb)\|_2.
\)
For \(\bb^k\in\mathcal L_0\), the nonexpansiveness of
the projection gives
\(
\|\bh^k(\tau)\|_2
\le
\tau\|\bg^k\|_2
\le
\tau M.
\)
The ascent lemma and
Eq.~\eqref{app:eq:projection_ascent} imply
\[
\begin{aligned}
Z(\bb^k+\bh^k(\tau))
&\ge
Z(\bb^k)
+(\bg^k)^\top\bh^k(\tau)
-\frac{L}{2}\|\bh^k(\tau)\|_2^2\\
&\ge
Z(\bb^k)
+\chi(\bg^k)^\top\bh^k(\tau)
\end{aligned}
\]
whenever
\(
\tau
\le
\widehat\tau
:=
\min\left\{
\frac{2(1-\chi)}{L},
\frac{r}{M}
\right\}.
\)
Here, \(r/M:=+\infty\) if \(M=0\), and the first ratio
is interpreted as \(+\infty\) if \(L=0\).
Let \(\tau_k\) denote the first accepted trial step.
Thus, backtracking terminates finitely, and every accepted
step satisfies
\begin{equation}
\tau_k
\ge
\underline\tau
:=
\min\left\{
\bar\tau,
\delta\widehat\tau
\right\}
>0.
\label{app:eq:step_lower_bound}
\end{equation}

For the accepted step, set
\(\bh^k:=\bh^k(\tau_k)=\bb^{k+1}-\bb^k\).
The Armijo condition,
Eq.~\eqref{app:eq:projection_ascent}, and
\(\tau_k\le\bar\tau\) yield
\begin{equation}
Z(\bb^{k+1})-Z(\bb^k)
\ge
\chi(\bg^k)^\top\bh^k
\ge
\frac{\chi}{\bar\tau}\|\bh^k\|_2^2.
\label{app:eq:sufficient_ascent}
\end{equation}
Since \(\bb^0\in\mathcal L_0\), the accepted iterates
remain in \(\mathcal L_0\) by induction.
By Lemma~\ref{lemma:continuity_existence}, \(Z\) is bounded
above on the compact set \(\mathcal B\).
Summing Eq.~\eqref{app:eq:sufficient_ascent} gives
\(
\sum_{k=0}^{\infty}\|\bh^k\|_2^2<\infty,
\)
and therefore \(\bh^k\to\boldsymbol{0}\).

Let \(\bb^{k_\nu}\to\bar{\bb}\) be any convergent subsequence.
By Eq.~\eqref{app:eq:step_lower_bound}, a further
subsequence satisfies
\(
\tau_{k_\nu}
\to
\bar\tau_0\in[\underline\tau,\bar\tau].
\)
The continuity of the gradient and projection, together
with
\(
\bb^{k_\nu+1}
=
\Pi_{\mathcal B}
\bigl(
\bb^{k_\nu}
+
\tau_{k_\nu}\nabla Z(\bb^{k_\nu})
\bigr)
\)
and \(\bh^{k_\nu}\to\boldsymbol{0}\), gives
\(
\bar{\bb}
=
\Pi_{\mathcal B}
\bigl(
\bar{\bb}
+
\bar\tau_0\nabla Z(\bar{\bb})
\bigr).
\)
Since \(\bar\tau_0>0\), the projection variational
inequality yields
\(
\nabla Z(\bar{\bb})^\top(\bb-\bar{\bb})
\le0,
\forall\bb\in\mathcal B.
\)
Thus, every accumulation point is stationary.

For any fixed \(\tau_{\mathrm{ref}}>0\), the residual
\(R_{\tau_{\mathrm{ref}}}\) is continuous on
\(\mathcal L_0\) and vanishes exactly at stationary points.
Suppose, to the contrary, that
Eq.~\eqref{app:eq:residual_convergence} fails.
Compactness then yields a subsequence with residuals
bounded away from zero and a stationary accumulation
point, contradicting the continuity of
\(R_{\tau_{\mathrm{ref}}}\).
Therefore, Eq.~\eqref{app:eq:residual_convergence} holds
for every fixed \(\tau_{\mathrm{ref}}>0\).
\Halmos
\end{proof}

\begin{corollary}[Complexity of accepted smooth steps]
\label{cor:smooth_accepted_step_complexity}
Under the assumptions of
Proposition~\ref{prop:convergence_stationary_final}, define
the projected-gradient mapping
\begin{equation}
\mathcal G_{\tau_k}(\bb^k)
:=
\frac{
\Pi_{\mathcal B}
\bigl(\bb^k+\tau_k\nabla Z(\bb^k)\bigr)-\bb^k
}{\tau_k}.
\label{eq:ec_smooth_gradient_mapping}
\end{equation}
For every integer \(K\ge1\),
\begin{equation}
\min_{0\le k<K}
\left\|
\mathcal G_{\tau_k}(\bb^k)
\right\|_2^2
\le
\frac{
Z^\star-Z(\bb^0)
}{
\chi\underline\tau K
},
\label{eq:ec_smooth_complexity}
\end{equation}
where \(\underline\tau>0\) is defined in
Eq.~\eqref{app:eq:step_lower_bound}.
Thus, the smooth specialization reaches an iterate with
gradient-mapping norm at most \(\varepsilon\) within
\(O(\varepsilon^{-2})\) accepted steps.
\end{corollary}

\begin{proof}{Proof.}
For an accepted step, the Armijo condition and
Eqs.~\eqref{app:eq:projection_ascent}
and~\eqref{app:eq:step_lower_bound} give
\begin{equation}
Z(\bb^{k+1})-Z(\bb^k)
\ge
\chi\tau_k
\left\|
\mathcal G_{\tau_k}(\bb^k)
\right\|_2^2
\ge
\chi\underline\tau
\left\|
\mathcal G_{\tau_k}(\bb^k)
\right\|_2^2.
\label{eq:ec_mapping_ascent}
\end{equation}
Summing Eq.~\eqref{eq:ec_mapping_ascent}
for \(k=0,\ldots,K-1\) and using
\(Z(\bb^K)\le Z^\star\) proves
Eq.~\eqref{eq:ec_smooth_complexity}.
\Halmos
\end{proof}

\subsection{Proof of Proposition~\ref{prop:sampled_clarke_limit}}
\label{app:sec:proof_proposition2}

\begin{proof}{Proof.}
For each \(\ell\), choose minimizers
\((\bomega^\ell,\bmu^\ell)\) defining
\(\widehat R_\ell^{\mathrm{amb}}\), and set
\(\bar{\bg}^{\,\ell}:=\bG_\ell^{\mathrm{amb}}\bomega^\ell\) and
\(\bar{\bn}^{\,\ell}:=\bN_\ell\bmu^\ell\).
Since \(\bb^\ell\to\bar{\bb}\) and \(\epsilon_\ell\to0\),
local Lipschitz continuity bounds the sampled gradients
for all sufficiently large \(\ell\), so
\(\{\bar{\bg}^{\,\ell}\}\) is bounded.
Moreover,
\(
\|\bar{\bg}^{\,\ell}-\bar{\bn}^{\,\ell}\|_2
=
\widehat R_\ell^{\mathrm{amb}}
\to0,
\)
so \(\{\bar{\bn}^{\,\ell}\}\) is also bounded.
Passing to a subsequence, we obtain
\(\bar{\bg}^{\,\ell}\to\bar{\bg}\),
\(\bar{\bn}^{\,\ell}\to\bar{\bn}\), and
\(\bar{\bg}=\bar{\bn}\).

By Carath\'eodory's theorem, each
\(\bar{\bg}^{\,\ell}\) can be represented using at most
\(\lvert\mathcal Z\rvert+1\) columns of
\(\bG_\ell^{\mathrm{amb}}\).
The corresponding sample points satisfy
\(
\|\by^{\ell,a}-\bar{\bb}\|_2
\le
\epsilon_\ell+\|\bb^\ell-\bar{\bb}\|_2
\to0.
\)
Boundedness of the gradients and compactness of the coefficient
simplex yield a subsequence along which all these terms converge.
The gradient-limit characterization and convexity of the Clarke
generalized gradient therefore give
\(\bar{\bg}\in\partial^C Z(\bar{\bb})\).

For every \(\widetilde{\bb}\in\mathcal B\), the representation
in terms of active outward normals gives
\(
(\bar{\bn}^{\,\ell})^\top
(\widetilde{\bb}-\bb^\ell)\le0.
\)
Taking limits yields
\(\bar{\bn}^\top(\widetilde{\bb}-\bar{\bb})\le0\)
for every \(\widetilde{\bb}\in\mathcal B\), and hence
\(\bar{\bn}\in\mathcal N_{\mathcal B}(\bar{\bb})\).
Since \(\bar{\bg}=\bar{\bn}\), we obtain
\(
\partial^C Z(\bar{\bb})\cap\mathcal N_{\mathcal B}(\bar{\bb})
\ne\varnothing,
\)
equivalently,
\(
\boldsymbol{0}\in-\partial^C Z(\bar{\bb})
+\mathcal N_{\mathcal B}(\bar{\bb}).
\)
Thus, \(\bar{\bb}\) is a constrained Clarke-stationary
point of DSMDP.
\Halmos
\end{proof}

\begin{remark}[Finite and limiting roles of the ambient residual]
For a fixed sampling radius, Proposition~\ref{prop:finite_scale_certificate} uses the ambient residual \(\widehat R_k^{\mathrm{amb}}\) to upper-bound the enlarged-Clarke stationarity distance.
The certificate applies at a fixed sampling scale and is therefore reported together with its radius.
Proposition~\ref{prop:sampled_clarke_limit} yields exact Clarke stationarity only along convergent ambient assessments whose radii and residuals both vanish.
\end{remark}

\section{Conic DSMDP Relaxation}
\label{app:conic_upper_bound}

Using valid bounds on the MNL utilities and denominators, we construct a conic outer relaxation based on exponential-cone epigraphs and affine secant upper bounds for the response terms and McCormick envelopes for the remaining products.
The complete formulation contains a lift of every feasible DSMDP solution and distinguishes the solver-certified upper bound from the feasible Conic benchmark used in the experiments.

\subsection{Bounds for the MNL Utilities and Denominators}
\label{app:utility_bounds}

Each normalized MNL utility depends affinely on a single bonus
coordinate, so valid bounds over \(\mathcal B\) follow by evaluating
the utility at \(0\) and \(b^{\max}\). For \(i\in\mathcal Z\) and
\(j\in\mathcal Z_i\), this gives
\(\underline V_{ij}
=\beta^{(3)}d_{ij}
+\min\{0,\beta^{(1)}b^{\max}\}\) and
\(\overline V_{ij}
=\beta^{(3)}d_{ij}
+\max\{0,\beta^{(1)}b^{\max}\}\), while
\(\underline V_i^{\mathrm{Off}}
=\delta_i^{\mathrm{Off}}
+\min\{0,-\zeta^{\mathrm{Off}}b^{\max}\}\) and
\(\overline V_i^{\mathrm{Off}}
=\delta_i^{\mathrm{Off}}
+\max\{0,-\zeta^{\mathrm{Off}}b^{\max}\}\)
for \(i\in\mathcal Z\).
These bounds imply
\(\underline t_{ij}=\exp(\underline V_{ij})\) and
\(\overline t_{ij}=\exp(\overline V_{ij})\) for
\(i\in\mathcal Z\) and \(j\in\mathcal Z_i\), and
\(\underline t_i^{\mathrm{Off}}
=\exp(\underline V_i^{\mathrm{Off}})\) and
\(\overline t_i^{\mathrm{Off}}
=\exp(\overline V_i^{\mathrm{Off}})\) for
\(i\in\mathcal Z\).
The denominator variable satisfies
\(D_i\in[\underline D_i,\overline D_i]\), where
\(\underline D_i=\sum_{m\in\mathcal Z_i}\underline t_{im}
+\underline t_i^{\mathrm{Off}}\) and
\(\overline D_i=\sum_{m\in\mathcal Z_i}\overline t_{im}
+\overline t_i^{\mathrm{Off}}\) for all \(i\in\mathcal Z\).

For any finite interval \([\ell,u]\), the affine secant
introduced in Section~\ref{subsec:conic_upper_bound} is
\begin{equation}
\mathcal S_{\exp}(v;\ell,u)
:=
\begin{cases}
\displaystyle
\frac{u-v}{u-\ell}\exp(\ell)
+
\frac{v-\ell}{u-\ell}\exp(u),
& \ell<u,\\[4pt]
\exp(\ell),
& \ell=u.
\end{cases}
\label{eq:app_exp_secant_definition}
\end{equation}
By convexity,
\(\exp(v)\le\mathcal S_{\exp}(v;\ell,u)\) for all
\(v\in[\ell,u]\).
When \(\ell=u\), the exponential lower bound and
secant upper bound coincide.

\subsection{Conic Relaxation of the MNL Probabilities}
\label{app:mnl_conic}

Using the variables introduced in
Section~\ref{subsec:conic_upper_bound}, impose
\(D_i=\sum_{m\in\mathcal Z_i}t_{im}+t_i^{\mathrm{Off}}
\in[\underline D_i,\overline D_i]\)
for each \(i\in\mathcal Z\).
For each \(i\in\mathcal Z\) and \(j\in\mathcal Z_i\),
the exponential and box constraints are
\begin{align}
t_{ij}
&\ge \exp(V_{ij}),
\label{eq:app_exp_epi}
\\
t_{ij}
&\le
\mathcal S_{\exp}\bigl(
V_{ij};
\underline V_{ij},
\overline V_{ij}
\bigr),
\label{eq:app_exp_secant}
\\
\underline t_{ij}
&\le t_{ij}\le\overline t_{ij}.
\label{eq:app_t_box}
\end{align}

The probability variables satisfy
\begin{equation}
\begin{aligned}
\sum_{j\in\mathcal Z_i}P_{ij}+P_i^{\mathrm{Off}}
&=1,
&&i\in\mathcal Z,\\
P_{ij}
&\ge0,
&&i\in\mathcal Z,\quad j\in\mathcal Z_i,\\
P_i^{\mathrm{Off}}
&\ge0,
&&i\in\mathcal Z.
\end{aligned}
\label{eq:app_P_norm}
\end{equation}

The exact probability--denominator relations are
\begin{equation}
\begin{aligned}
t_{ij}
&=
P_{ij}D_i,
\qquad
i\in\mathcal Z,\;j\in\mathcal Z_i,\\
t_i^{\mathrm{Off}}
&=
P_i^{\mathrm{Off}}D_i,
\qquad
i\in\mathcal Z.
\end{aligned}
\label{eq:app_probability_products}
\end{equation}
Equation~\eqref{eq:app_P_norm} implies
\(P_{ij}\in[0,1]\) and \(P_i^{\mathrm{Off}}\in[0,1]\).
Together with \(D_i\in[\underline D_i,\overline D_i]\),
these bounds yield the following McCormick envelopes.

For each \(i\in\mathcal Z\) and \(j\in\mathcal Z_i\),
the envelope for \(t_{ij}=P_{ij}D_i\) is
\begin{align}
t_{ij}
&\ge \underline D_iP_{ij},
\label{eq:app_mcc_t_1}
\\
t_{ij}
&\ge D_i+\overline D_i(P_{ij}-1),
\\
t_{ij}
&\le \overline D_iP_{ij},
\\
t_{ij}
&\le D_i+\underline D_i(P_{ij}-1).
\end{align}

For each \(i\in\mathcal Z\), the envelope for
\(t_i^{\mathrm{Off}}=P_i^{\mathrm{Off}}D_i\) is
\begin{align}
t_i^{\mathrm{Off}}
&\ge \underline D_iP_i^{\mathrm{Off}},
\\
t_i^{\mathrm{Off}}
&\ge D_i+\overline D_i(P_i^{\mathrm{Off}}-1),
\\
t_i^{\mathrm{Off}}
&\le \overline D_iP_i^{\mathrm{Off}},
\\
t_i^{\mathrm{Off}}
&\le D_i+\underline D_i(P_i^{\mathrm{Off}}-1).
\label{eq:app_mcc_toff_4}
\end{align}

Equations~\eqref{eq:app_mcc_t_1}--\eqref{eq:app_mcc_toff_4}
define the convex hulls of the graphs
\(t_{ij}=P_{ij}D_i\) and
\(t_i^{\mathrm{Off}}=P_i^{\mathrm{Off}}D_i\)
over their respective boxes.

\subsection{Conic Relaxation of Latent Driver Activation and the Unavailable State Exponential Term}
\label{app:onoff_conic}

We express the latent-driver activation probability using a binary logit model.
Define
\(
U_i^{\mathrm{On}}
:=
\alpha_i^{\mathrm{On}}
+
\gamma^{\mathrm{On}}b_i
\)
and the negative utility
\(
\widetilde U_i^{\mathrm{On}}
:=
-U_i^{\mathrm{On}}.
\)
Introduce the exponential variable \(e_i^{\mathrm{On}}\), the
denominator variable \(D_i^{\mathrm{On}}\), and the probability variable
\(\varphi_i^{\mathrm{On}}\)
representing \(P_i^{\mathrm{On}}(\bb)\). Their exact relations are
\(
e_i^{\mathrm{On}}
=
\exp(\widetilde U_i^{\mathrm{On}}),
\),
\(
D_i^{\mathrm{On}}
=
1+e_i^{\mathrm{On}},
\), and
\(
1
=
\varphi_i^{\mathrm{On}}D_i^{\mathrm{On}}
\).
The endpoint utility bounds are
\(
\underline U_i^{\mathrm{On}}
=
\min\left\{
\alpha_i^{\mathrm{On}},
\alpha_i^{\mathrm{On}}
+
\gamma^{\mathrm{On}}b^{\max}
\right\},
\) and
\(
\overline U_i^{\mathrm{On}}
=
\max\left\{
\alpha_i^{\mathrm{On}},
\alpha_i^{\mathrm{On}}
+
\gamma^{\mathrm{On}}b^{\max}
\right\},
\forall i\in\mathcal Z
\).
For \(\gamma^{\mathrm{On}}\ge0\), these bounds reduce to
\(
\underline U_i^{\mathrm{On}}
=
\alpha_i^{\mathrm{On}},
\qquad
\overline U_i^{\mathrm{On}}
=
\alpha_i^{\mathrm{On}}
+
\gamma^{\mathrm{On}}b^{\max}.
\)
The negative-utility bounds are
\(
\underline{\widetilde U}_i^{\mathrm{On}}
=
-\overline U_i^{\mathrm{On}},
\) and
\(
\overline{\widetilde U}_i^{\mathrm{On}}
=
-\underline U_i^{\mathrm{On}}.
\)
The corresponding exponential, denominator, and probability bounds are
\(
\underline e_i^{\mathrm{On}}
:=
\exp(\underline{\widetilde U}_i^{\mathrm{On}})
=
\exp(-\overline U_i^{\mathrm{On}})
\)
and
\(
\overline e_i^{\mathrm{On}}
:=
\exp(\overline{\widetilde U}_i^{\mathrm{On}})
=
\exp(-\underline U_i^{\mathrm{On}})
\);
\(
\underline D_i^{\mathrm{On}}
:=
1+\underline e_i^{\mathrm{On}}
\)
and
\(
\overline D_i^{\mathrm{On}}
:=
1+\overline e_i^{\mathrm{On}}
\);
and
\(
\underline \varphi_i^{\mathrm{On}}
:=
\sigma(\underline U_i^{\mathrm{On}})
\)
and
\(
\overline \varphi_i^{\mathrm{On}}
:=
\sigma(\overline U_i^{\mathrm{On}})
\),
for all \(i\in\mathcal Z\).

For each \(i\in\mathcal Z\), we use an exponential-cone epigraph, an affine secant upper bound, and a McCormick envelope.
The epigraph and secant constraints enclose the exponential relation over the valid utility interval.
The McCormick inequalities relax the reciprocal product between the activation probability and its denominator.
\begin{align}
&e_i^{\mathrm{On}}
\ge
\exp(\widetilde U_i^{\mathrm{On}}),
\label{eq:app_on_exp_epi}
\\
&e_i^{\mathrm{On}}
\le
\mathcal S_{\exp}
\bigl(
\widetilde U_i^{\mathrm{On}};
\underline{\widetilde U}_i^{\mathrm{On}},
\overline{\widetilde U}_i^{\mathrm{On}}
\bigr),
\label{eq:app_on_exp_secant}
\\
&D_i^{\mathrm{On}}
=
1+e_i^{\mathrm{On}},
\\
&\underline D_i^{\mathrm{On}}
\le
D_i^{\mathrm{On}}
\le
\overline D_i^{\mathrm{On}},
\\
&\underline \varphi_i^{\mathrm{On}}
\le
\varphi_i^{\mathrm{On}}
\le
\overline \varphi_i^{\mathrm{On}},
\\
&1
\ge
\underline \varphi_i^{\mathrm{On}}D_i^{\mathrm{On}}
+
\underline D_i^{\mathrm{On}}\varphi_i^{\mathrm{On}}
-
\underline \varphi_i^{\mathrm{On}}
\underline D_i^{\mathrm{On}},
\label{eq:app_on_mcc_1}
\\
&1
\ge
\overline \varphi_i^{\mathrm{On}}D_i^{\mathrm{On}}
+
\overline D_i^{\mathrm{On}}\varphi_i^{\mathrm{On}}
-
\overline \varphi_i^{\mathrm{On}}
\overline D_i^{\mathrm{On}},
\\
&1
\le
\overline \varphi_i^{\mathrm{On}}D_i^{\mathrm{On}}
+
\underline D_i^{\mathrm{On}}\varphi_i^{\mathrm{On}}
-
\overline \varphi_i^{\mathrm{On}}
\underline D_i^{\mathrm{On}},
\\
&1
\le
\underline \varphi_i^{\mathrm{On}}D_i^{\mathrm{On}}
+
\overline D_i^{\mathrm{On}}\varphi_i^{\mathrm{On}}
-
\underline \varphi_i^{\mathrm{On}}
\overline D_i^{\mathrm{On}}.
\label{eq:app_on_mcc_4}
\end{align}
Every exact binary logit point satisfies
Eqs.~\eqref{eq:app_on_exp_epi}--\eqref{eq:app_on_mcc_4}. These constraints
therefore define an outer relaxation of the latent-driver activation response.
For each \(i\in\mathcal Z\), the exponential term
\(t_i^{\mathrm{Off}}\) satisfies the following bounds.
\begin{align}
t_i^{\mathrm{Off}}
&\ge
\exp(V_i^{\mathrm{Off}}),
\label{eq:app_off_exp_epi}
\\
t_i^{\mathrm{Off}}
&\le
\mathcal S_{\exp}\bigl(
V_i^{\mathrm{Off}};
\underline V_i^{\mathrm{Off}},
\overline V_i^{\mathrm{Off}}
\bigr),
\label{eq:app_off_exp_secant}
\\
\underline t_i^{\mathrm{Off}}
&\le
t_i^{\mathrm{Off}}
\le
\overline t_i^{\mathrm{Off}}.
\label{eq:app_off_t_box}
\end{align}
With the McCormick relaxation of
\(t_i^{\mathrm{Off}}=P_i^{\mathrm{Off}}D_i\), these constraints retain
every exact departure from the dispatchable pool and latent-driver activation response within the lifted feasible
region.

\subsection{McCormick Relaxation of Bonus Payments}
\label{app:payment_mccormick}

For the matched quantity \(q_i\) defined in
Section~\ref{subsec:conic_upper_bound}, combining total
demand with the maximum possible post-response supply
gives the upper bound
\begin{align}
\bar q_i
&=
\min\left\{
\sum_{j\in\mathcal Z}d_j,\;
s_i^{\mathrm{On}}
+
\sum_{j\in\mathcal Z_i\setminus\{i\}}s_j^{\mathrm{On}}
+
s_i^{\mathrm{Off}}
\right\},
\qquad
\forall i\in\mathcal Z,
\label{eq:app_qbar}
\\
0
&\le
q_i
\le
\bar q_i,
\qquad
\forall i\in\mathcal Z.
\label{eq:app_q_box}
\end{align}
Under the symmetric reachability relation and
\(i\in\mathcal Z_i\), the dispatchable-supply term in
Eq.~\eqref{eq:app_qbar} can equivalently be written as
\(
\sum_{j\in\mathcal Z:\,i\in\mathcal Z_j}s_j^{\mathrm{On}}
\).
For the New York City instances, every sum involving
\(x_{ij}\) is restricted to the sparse matching arc set
described in Section~\ref{exp_case}; the original model
and its outer relaxation use the same arc set.

Over
\(
(b_i,q_i)\in[0,b^{\max}]\times[0,\bar q_i]
\),
the McCormick envelope for \(w_i=b_iq_i\) is
\begin{align}
w_i
&\ge
0,
\label{eq:app_w_mcc_1}
\\
w_i
&\ge
b^{\max}q_i
+
\bar q_i b_i
-
b^{\max}\bar q_i,
\\
w_i
&\le
b^{\max}q_i,
\\
w_i
&\le
\bar q_i b_i,
\qquad \forall i\in\mathcal Z.
\label{eq:app_w_mcc_4}
\end{align}
Equations~\eqref{eq:app_w_mcc_1}--\eqref{eq:app_w_mcc_4} give the McCormick envelope of the payment product.

\subsection{Full Conic Relaxation}
\label{app:full_conic_relaxation}

Collect the decision variables in
\(
\boldsymbol{\vartheta}
:=
(\bx,\bq,\bw,\bb,\bP,\bP^{\mathrm{Off}},
\bt,\bt^{\mathrm{Off}},(D_i)_{i\in\mathcal Z},\boldsymbol{\varphi}^{\mathrm{On}},
\boldsymbol{e}^{\mathrm{On}},\bD^{\mathrm{On}})
\).
The full conic relaxation is
\begingroup
\small
\setlength{\jot}{0pt}
\begin{subequations}
\label{eq:app_conic_DSMDP}
\begin{align}
(\overline{\mathrm{DSMDP}})\qquad
\overline Z
=
\max_{\boldsymbol{\vartheta}}
\quad
&
\sum_{i,j\in\mathcal Z}p_{ij}x_{ij}
-
\sum_{i\in\mathcal Z}w_i
\label{eq:app_conic_obj}
\\
\text{s.t.}\qquad
&
\sum_{i\in\mathcal Z}x_{ij}\le d_j,
\qquad
\forall j\in\mathcal Z,
\label{eq:app_conic_demand}
\\
&
\begin{aligned}
\sum_{j\in\mathcal Z}x_{ij}
&\le
s_i^{\mathrm{On}}
+
\sum_{j\in\mathcal Z_i\setminus\{i\}}
s_j^{\mathrm{On}}P_{ji} \\
&\quad-
s_i^{\mathrm{On}}
\left(
\sum_{j\in\mathcal Z_i\setminus\{i\}}P_{ij}
+
P_i^{\mathrm{Off}}
\right)
+
s_i^{\mathrm{Off}}\varphi_i^{\mathrm{On}},
\qquad
\forall i\in\mathcal Z.
\end{aligned}
\label{eq:app_conic_supply}
\\
&
\sum_{i\in\mathcal Z}w_i\le B,
\label{eq:app_conic_budget_q}
\\
&
q_i=\sum_{j\in\mathcal Z}x_{ij},
&&i\in\mathcal Z,
\label{eq:app_conic_q_definition}
\\
&
\eqref{eq:app_q_box},
\qquad
\eqref{eq:app_w_mcc_1}--\eqref{eq:app_w_mcc_4},
\label{eq:app_conic_q_w_bounds}
\\
&
\sum_{j\in\mathcal Z_i}P_{ij}+P_i^{\mathrm{Off}}=1,
&&i\in\mathcal Z,
\label{eq:app_conic_P_norm}
\\
&
P_{ij}\ge0,
&&i\in\mathcal Z,\quad j\in\mathcal Z_i,
\label{eq:app_conic_P_domain}
\\
&
P_i^{\mathrm{Off}}\ge0,
&&i\in\mathcal Z,
\label{eq:app_conic_Poff_domain}
\\
&
D_i=\sum_{m\in\mathcal Z_i}t_{im}+t_i^{\mathrm{Off}},
&&i\in\mathcal Z,
\label{eq:app_conic_D}
\\
&
\underline D_i\le D_i\le\overline D_i,
&&i\in\mathcal Z,
\label{eq:app_conic_D_bounds}
\\
&
\eqref{eq:app_exp_epi}--\eqref{eq:app_t_box},
\qquad
\eqref{eq:app_off_exp_epi}--\eqref{eq:app_off_t_box},
\label{eq:app_conic_exp_outer}
\\
&
\eqref{eq:app_mcc_t_1}--\eqref{eq:app_mcc_toff_4},
\label{eq:app_conic_mnl_mcc}
\\
&
\eqref{eq:app_on_exp_epi}--\eqref{eq:app_on_mcc_4},
\label{eq:app_conic_on_outer}
\\
&\bb\in\mathcal B,
\label{eq:app_conic_bonus_domain}
\\
&x_{ij}\ge0,
&&i,j\in\mathcal Z.
\label{eq:app_conic_domains}
\end{align}
\end{subequations}
\endgroup

Objective~\eqref{eq:app_conic_obj} maximizes the relaxed net operating profit.
Constraints~\eqref{eq:app_conic_demand} and \eqref{eq:app_conic_supply} impose demand capacity and the relaxed post-response supply balance.
Constraint~\eqref{eq:app_conic_budget_q} enforces the bonus budget, Constraint~\eqref{eq:app_conic_q_definition} defines matched quantity, and Constraint~\eqref{eq:app_conic_q_w_bounds} invokes the bounds on matched quantities and payment envelopes.
Constraint~\eqref{eq:app_conic_P_norm} normalizes the relocation and unavailability probabilities, while Constraints~\eqref{eq:app_conic_P_domain} and \eqref{eq:app_conic_Poff_domain} impose their domains.
Constraints~\eqref{eq:app_conic_D} and \eqref{eq:app_conic_D_bounds} define and bound the MNL denominators.
Constraints~\eqref{eq:app_conic_exp_outer}--\eqref{eq:app_conic_mnl_mcc} impose the exponential and product relaxations for relocation and unavailability, whereas Constraint~\eqref{eq:app_conic_on_outer} imposes the corresponding relaxation for latent-driver activation.
Constraints~\eqref{eq:app_conic_bonus_domain} and \eqref{eq:app_conic_domains} retain the feasible driver surge maps and nonnegative matching flows.

\subsection{Validity of the Upper Bound}
\label{app:upper_bound_validity}

The upper-bound proof constructs a feasible lift of an arbitrary DSMDP solution into \(\overline{\mathrm{DSMDP}}\) without changing its objective value.
The feasible-set inclusion establishes the bound independently of the feasible-recovery procedure used for the Conic benchmark.

\begin{proposition}[Validity of the conic upper bound]
\label{prop:conic_upper_bound_validity1}

Let \(Z^\star\) denote the optimal objective value of the original DSMDP,
and let \(\overline Z\) denote the optimal objective value of the conic
outer relaxation \(\overline{\mathrm{DSMDP}}\) in
Eq.~\eqref{eq:app_conic_DSMDP}. Then
\(
Z^\star
\le
\overline Z.
\)
\end{proposition}

\begin{proof}{Proof.}
Consider any feasible solution \((\bb,\bx)\) of the original DSMDP.
For each \(i\in\mathcal Z\) and \(j\in\mathcal Z_i\), set
\(t_{ij}=\exp(V_{ij})\),
\(t_i^{\mathrm{Off}}=\exp(V_i^{\mathrm{Off}})\), and
\(D_i=\sum_{m\in\mathcal Z_i}t_{im}+t_i^{\mathrm{Off}}\).
Take \(P_{ij}=P_{ij}(\bb)=t_{ij}/D_i\) and
\(P_i^{\mathrm{Off}}=P_i^{\mathrm{Off}}(\bb)
=t_i^{\mathrm{Off}}/D_i\).
These assignments satisfy the probability normalization, denominator
definition, exponential-cone epigraph constraints, and all prescribed
bounds. By convexity of the exponential function, the exact
exponential values also satisfy the secant upper bounds in
Eqs.~\eqref{eq:app_exp_secant} and \eqref{eq:app_off_exp_secant}.
The assignments further satisfy
\(t_{ij}=P_{ij}D_i\) and
\(t_i^{\mathrm{Off}}=P_i^{\mathrm{Off}}D_i\).
Every point on either bilinear graph satisfies its McCormick envelope,
so all MNL outer-relaxation constraints hold.
For the latent-driver activation response, set
\(U_i^{\mathrm{On}}
=\alpha_i^{\mathrm{On}}+\gamma^{\mathrm{On}}b_i\),
\(e_i^{\mathrm{On}}=\exp(-U_i^{\mathrm{On}})\),
\(D_i^{\mathrm{On}}=1+e_i^{\mathrm{On}}\), and
\(\varphi_i^{\mathrm{On}}
=P_i^{\mathrm{On}}(\bb)=1/D_i^{\mathrm{On}}\).
The exact exponential value satisfies
Eq.~\eqref{eq:app_on_exp_epi}, the secant upper bound in
Eq.~\eqref{eq:app_on_exp_secant}, and its box bounds, while
\(1=\varphi_i^{\mathrm{On}}D_i^{\mathrm{On}}\) satisfies every McCormick
inequality in
Eqs.~\eqref{eq:app_on_mcc_1}--\eqref{eq:app_on_mcc_4}.
Hence, the latent-driver activation outer constraints contain the exact binary logit
response.

Finally, set
\(q_i=\sum_{j\in\mathcal Z}x_{ij}\) and \(w_i=b_iq_i\).
The point \((b_i,q_i,w_i)\) lies on the original product graph over
\([0,b^{\max}]\times[0,\bar q_i]\) and therefore satisfies the payment
McCormick envelope. Under these assignments, the relaxed supply
constraint coincides with the original supply constraint, and the
relaxed objective and budget coincide with their original counterparts.
Thus, every feasible DSMDP solution has a feasible lift into
\(\overline{\mathrm{DSMDP}}\) with the same objective value, which
gives \(Z^\star\le\overline Z\).
\Halmos
\end{proof}

\subsection{Solver Bound and Feasible Recovery}
\label{app:reported_bound_recovery}

We report the solver-certified upper bound
\(U^{\mathrm{solver}}\) only when a valid bound is available
after an optimal or time-limit termination.
It satisfies
\(
Z^\star
\le
\overline Z
\le
U^{\mathrm{solver}}.
\)
If the relaxation is solved to global optimality,
then \(U^{\mathrm{solver}}=\overline Z\).

To construct the feasible Conic benchmark, let
\(\bb^{\mathrm{raw}}\) and \(\bx^{\mathrm{raw}}\) denote
the bonus and matching variables from the relaxation
incumbent, and define
\(
q_i^{\mathrm{raw}}
=
\sum_{j\in\mathcal Z}x_{ij}^{\mathrm{raw}}.
\)
We project \(\bb^{\mathrm{raw}}\) onto \(\mathcal B\)
using outflow-based weights \(q_i^{\mathrm{raw}}\) when
\(q_i^{\mathrm{raw}}>10^{-9}\) and unit weights otherwise.
We then recompute driver responses at the projected
bonus vector and solve DSMDP-M, which enforces the
realized bonus-spending constraint
\(\sum_{i,j\in\mathcal Z}b_i x_{ij}\le B\).
The resulting solution is feasible for the original
DSMDP, and its net operating profit
\(Z^{\mathrm{Conic}}\) is reported separately from
\(U^{\mathrm{solver}}\).
Thus, whenever \(U^{\mathrm{solver}}\) is available,
\(
Z^{\mathrm{Conic}}
\le
Z^\star
\le
\overline Z
\le
U^{\mathrm{solver}}.
\)

\section{Experimental Settings and Benchmarks}\label{baselines}

All feasible methods are compared under a common protocol, while the solver-certified upper bound is reported separately.
We specify the shared data generation, initialization, randomization, solver, timing, and ISA settings, then define the benchmark implementations and the within-protocol fixed-time and time-to-target comparisons.
The resulting protocol separates feasible method performance from the quality information supplied by the relaxation bound.

\subsection{Parameter Settings}
\label{app:parameter_settings}
We vary the system configurations across experiments while keeping the
following parameters fixed. For rider demand generation, the base
demand level is \(d^{\mathrm{base}}=5.0\), the demand-spread parameter
is \(\sigma^d=2.0\), and the hotspot intensity is
\(\theta^d=25.0\).
For initial dispatchable-supply generation, the baseline dispatchable-supply level
is \(s^{\mathrm{On,base}}=3.0\), the activity-center intensity is
\(\theta^s=8.0\), and the supply-spread parameter is
\(\sigma^s=1.5\). The supply-shift radius is \(r^s=4.0\), and the
supply-shift strength is \(\phi^s=0.005\).

For latent-supply generation, we set the latent-supply-to-demand ratio to $\rho^{\mathrm{Off}}=0.30$ and the lognormal noise parameter to $\sigma^{\mathrm{Off}}=0.3$. We generate the supply from a demand-weighted raw latent-supply vector and then rescale it so that total latent supply equals $\rho^{\mathrm{Off}}$ times total demand.
For pre-bonus matching-value generation, the pickup-distance term uses an
exponential distance-decay function with base value
\(\theta^p=15.0\) and decay rate \(\xi^p=0.5\). 
Destination
attractiveness uses spread \(2.5\), intensity \(12.0\), and weight
\(0.7\). The behavioral-response parameters are set to
\(
\beta^{(1)}=0.1,
\beta^{(3)}=-0.25,
\delta_i^{\mathrm{Off}}\equiv0.03,
\zeta^{\mathrm{Off}}=0.15,
\alpha_i^{\mathrm{On}}\equiv-2.0,
\gamma^{\mathrm{On}}=1.5.
\)
Empirical studies provide evidence for both behavioral-response channels: monetary incentives influence drivers' relocation decisions \citep{lu2018surge,ashkrof2024relocation} and labor supply \citep{chen2022incentive}. These findings motivate our behavioral-response specification, while the effects of parameter choices are examined through the activation-responsiveness and behavioral-parameter estimation error analyses in Section~\ref{exp_sensi}.
Unless the corresponding sensitivity analysis varies one of these
values, the generated experiments use \(B=500\), \(b^{\max}=15\), and
\(\Delta^b=1\).
All generated instances use data seed 42.
The behavioral-parameter estimation error experiment uses the three initial dispatchable-supply-to-demand ratios
\(\kappa\in\{0.50,1.00,1.10\}\) and 10 independent parameter draws for
each combination of ratio and perturbation level; its draw streams are indexed by the ratio, perturbation level, and draw number.

For ISA, the initial map allocates \(0.3B\) across zones
in proportion to demand and is projected onto \(\mathcal B\). 
The initial sampling radius is
\(\epsilon_0=0.01b^{\max}=0.15\), and the initial stationarity
scale is \(\nu_0=0.1\max\{1,\|\bg^0\|_2\}\).
Both scales contract by a factor of \(\kappa_\epsilon=\kappa_\nu=0.5\).
The sample pool contains at most \(2|\mathcal Z|+1\)
sensitivities, begins with the center and up to eight ambient samples,
and adds ambient samples in batches of at most eight. A scale
contraction is preceded by a residual test using at least
\(|\mathcal Z|+1\) ambient sensitivities. The projected bundle Armijo
search uses initial step \(1\), minimum trial step \(10^{-14}\), step contraction \(0.5\), Armijo
parameter \(10^{-4}\), and at most 70 backtracking trials.
The main comparisons use sampling seed 2026.
Each ordinary ISA phase has a 120-second optimization deadline, 300 accepted steps, and 1,500 sampling attempts; the extended New York City restart has an additional 1,800-second optimization deadline with the same accepted-step and attempt limits.
NOMAD has a 120-second optimization deadline and 2,000 black-box evaluations, and the conic outer relaxation has a 1,800-second optimization deadline.
Unless stated otherwise, the target sampling and stationarity scales are one third of their initial values, and the target objective-stability test uses a window of 20 attempts and a relative threshold of \(10^{-4}\).
Every ISA execution returns the best center reached before its stopping condition or computation limit.

In the New York City study, ISA, ISA-CG, and NOMAD start from the \(\mathbf H_1\) map.
All experiments use Python 3.13.12, Gurobi 13.0.2, and PyNomadBBO 4.5.1 on the Windows server described in Section~\ref{exp_benchmarks}.
In each realization, every dispatchable driver in zone \(i\) independently selects a reachable zone \(j\in\mathcal Z_i\) or the unavailable state with probabilities \(P_{ij}(\bb)\) and \(P_i^{\mathrm{Off}}(\bb)\).
The latent supply \(s_i^{\mathrm{Off}}\) is rounded to \(\lfloor s_i^{\mathrm{Off}}\rfloor\) plus a Bernoulli variable with mean \(s_i^{\mathrm{Off}}-\lfloor s_i^{\mathrm{Off}}\rfloor\), and each latent driver becomes available independently with probability \(P_i^{\mathrm{On}}(\bb)\).
The matching-aware fixed point starts from the matching probabilities of the nominal expected-supply solution, uses damping weight \(0.5\), and stops when \(\max_{i\in\mathcal Z}|a_i^{\mathrm{new}}-a_i|\le10^{-6}\), where \(a_i=\min\{1,q_i/s_i\}\).
ISA-CG applies the projected Armijo procedure in \ECsection{EC.3.4} with the line-search parameters and computation limits of ISA, and stops when \(\|\Pi_{\mathcal B}(\bb+\nabla Z(\bb))-\bb\|_2\le10^{-6}\) or no trial step is accepted.
The binding-budget comparison uses \(B=1000\).

\subsection{Anytime Performance and Efficiency of the Matching LP}
\label{ec:anytime_benchmark}

Table~\ref{ec:tab_anytime} summarizes ISA's anytime performance
on 26 generated instances.
ISA reaches the common \(95\%\) target in a median of \(4.52\)
seconds and outperforms the best completed non-ISA feasible
solution on all instances by 30 seconds.
As higher targets require progressively more computation,
these results suggest that platforms can obtain high-quality
feasible bonus maps early in the decision window and use
the remaining time for further refinement.

\begin{table}[htbp]
\centering
\caption{Anytime performance on 26 unique generated instances.}
\label{ec:tab_anytime}
\scriptsize
\setlength{\tabcolsep}{2pt}
\compacttablerows
\begin{tabular*}{\textwidth}{@{\extracolsep{\fill}}rrrrrr@{}}
\toprule
\multicolumn{6}{c}{Panel A: Fixed-time ISA objective lead over NOMAD} \\
\cmidrule(lr){1-6}
Time (s) & Wins & Losses & Median (\%) & IQR (\%) & Range (\%) \\
\midrule
10  & 26 & 0 & 14.39 & [11.66, 22.18] & [0.02, 33.84] \\
30  & 26 & 0 & 14.44 & [6.55, 21.50]  & [0.58, 33.84] \\
60  & 26 & 0 & 14.64 & [5.49, 19.01]  & [0.41, 37.62] \\
120 & 26 & 0 &  8.15 & [5.12, 14.92]  & [0.36, 19.37] \\
\midrule
\multicolumn{6}{c}{Panel B: Attainment of common within-protocol targets} \\
\cmidrule(lr){1-6}
\(\alpha\)
& \makecell{ISA\\instances}
& \makecell{NOMAD\\instances}
& \makecell{ISA time (s)\\median [IQR]}
& \makecell{Max. time\\(s)}
& \makecell{ISA LP solves\\median [IQR]} \\
\midrule
0.90 & 26 & 0 & 2.37 [1.15, 10.68]  & 48.63
& 132.00 [54.00, 301.25] \\
0.95 & 26 & 0 & 4.52 [1.22, 11.22]  & 59.44
& 209.50 [114.75, 318.25] \\
0.99 & 26 & 0 & 9.31 [3.21, 22.89] & 73.00
& 468.50 [315.00, 640.75] \\
\midrule
\multicolumn{6}{c}{Panel C: ISA lead over the best completed non-ISA feasible solution} \\
\cmidrule(lr){1-6}
ISA time (s) & Wins & Losses & Median (\%) & IQR (\%) & Range (\%) \\
\midrule
30 & 26 & 0 & 7.05 & [5.12, 13.72] & [0.36, 19.37] \\
60 & 26 & 0 & 8.15 & [5.12, 14.91] & [0.36, 19.37] \\
\bottomrule
\end{tabular*}
\compacttableend
\end{table}
\FloatBarrier

\subsection{Benchmark Implementations}
\label{app:benchmark_implementations}

\paragraph{Baseline 1: MNL-Response-Aware Greedy Heuristic.}
\label{app_h1}

\(\mathbf H_1\) constructs a one-pass driver surge map from the no-bonus map using the score
\begin{equation}
\psi_m
=
\max\{d_m-s_m(\boldsymbol{0}),0\}
\max_{j\in\mathcal Z}p_{mj}
\sum_{i:\,m\in\mathcal Z_i}
s_i^{\mathrm{On}}P_{im}(\boldsymbol{0})
\bigl(1-P_{im}(\boldsymbol{0})\bigr).
\label{ec:eq:h1_score}
\end{equation}
After solving the no-bonus matching LP for
\(\bx^{\mathrm{no}}\), \(\mathbf H_1\) uses the score in Eq.~\eqref{ec:eq:h1_score} to set
\(\hat b_m=\max\{\psi_m,0\}\) and, when the denominator is positive, rescales
\(\hat{\bb}\leftarrow
B\hat{\bb}/
(\sum_{m\in\mathcal Z}\hat b_m\sum_{j\in\mathcal Z}x_{mj}^{\mathrm{no}})\).
It then clips the bonus vector at \(b^{\max}\), projects it onto \(\mathcal B\), and
re-solves DSMDP-M to obtain
\((\bb^{\mathbf H_1},\bx^{\mathbf H_1},Z^{\mathbf H_1})\).

\paragraph{Baseline 2: LP-Dual-Informed Heuristic.}
\label{app_h2}

The second heuristic, denoted by \(\mathbf H_2\), uses a selected optimal
branch of the no-bonus matching LP. Let
\(\bx^{\mathrm{no}}\), \(\boldsymbol{\eta}^{\mathrm{no}}\), and \(\lambda^{\mathrm{no}}\) denote its
primal solution, supply duals, and budget dual, respectively. Define
\begin{equation}
q_m^{\mathrm{no}}:=\sum_j x_{mj}^{\mathrm{no}},
\quad
g_m^{\mathrm{no}}:=
\sum_{i\in\mathcal Z}\eta_i^{\mathrm{no}}
\frac{\partial s_i(\boldsymbol{0})}{\partial b_m}
-(1+\lambda^{\mathrm{no}})q_m^{\mathrm{no}},
\quad
\bar g:=\sum_{m\in\mathcal Z}g_m^{\mathrm{no}}.
\label{ec:eq:h2_score}
\end{equation}
Using the quantities defined in Eq.~\eqref{ec:eq:h2_score}, \(\mathbf H_2\) returns
\(\bb^{\mathbf H_2}=\boldsymbol{0}\) if \(\bar g<0\) or \(\max_i\eta_i^{\mathrm{no}}=0\).
Otherwise, for tuning parameters \(c_\tau>0\), \(\tau^{\max}>0\),
and \(\omega\in[0,1]\), it sets
\(\tau=\min\{\bar g/c_\tau,\tau^{\max}\}\) and
\(\hat b_i=\tau[(1-\omega)+\omega\eta_i^{\mathrm{no}}/\max_\ell\eta_\ell^{\mathrm{no}}]\).
It then clips the bonus vector to \([0,b^{\max}]\),
adjusts neighboring bonus differences, and clips the vector again.
Finally, it recomputes driver responses and re-solves DSMDP-M
under the resulting \(\bb^{\mathbf H_2}\) for
\(\bx^{\mathbf H_2}\) and \(Z^{\mathbf H_2}\).

\paragraph{Baseline 3: NOMAD.}
NOMAD applies mesh adaptive direct search \citep{le2011nomad} to the bonus vector and solves the matching LP for each candidate.
The black-box procedure enforces the constraints defining \(\mathcal B\), while the matching LP enforces the realized bonus-spending constraint.
The separation allows NOMAD to search over implementable bonus maps while each black-box evaluation returns the conditional matching optimum.

\subsection{Sampling Ablation: ISA versus ISA-CG}
\label{app:isa_cg_ablation}

Table~\ref{tab:stochastic_ablation} reports the
original objectives and ISA's relative improvements
over ISA-CG under Original, MC, and MA-MC evaluations
for\(B=1500\) and \(B=1000\). At \(B=1500\), ISA-CG slightly outperforms ISA. At \(B=1000\), however, ISA-CG takes substantially smaller accepted steps near budget-related matching transitions, whereas ISA achieves higher profits in all nine instances by 0.19\%--1.43\%.

\begin{table}[htbp]
\centering
\caption{ISA relative to ISA-CG under slack and binding budgets in the New York City study.}
\label{tab:stochastic_ablation}
{\footnotesize
\setlength{\tabcolsep}{2pt}
\renewcommand{\arraystretch}{0.9}
\begin{tabular*}{0.98\linewidth}{@{\extracolsep{\fill}}c*{10}{c}@{}}
\toprule
& \multicolumn{5}{c}{\(B=1500\)}
& \multicolumn{5}{c}{\(B=1000\)} \\
\cmidrule(lr){2-6}
\cmidrule(lr){7-11}
& \multicolumn{2}{c}{Original objective}
& \multicolumn{3}{c}{ISA vs.\ ISA-CG (\%)}
& \multicolumn{2}{c}{Original objective}
& \multicolumn{3}{c}{ISA vs.\ ISA-CG (\%)} \\
\cmidrule(lr){2-3}
\cmidrule(lr){4-6}
\cmidrule(lr){7-8}
\cmidrule(lr){9-11}
ID & ISA-CG & ISA & Original & MC & MA-MC
& ISA-CG & ISA & Original & MC & MA-MC \\
\midrule
M1 & 8715 & 8703 & -0.14 & -0.13 & -0.13 & 8643 & 8659 & +0.19 & +0.21 & +0.21 \\
M2 & 8996 & 8974 & -0.24 & -0.24 & -0.24 & 8854 & 8904 & +0.57 & +0.79 & +0.79 \\
M3 & 8892 & 8881 & -0.12 & -0.13 & -0.13 & 8733 & 8814 & +0.92 & +1.21 & +1.21 \\
\midrule
O1 & 11886 & 11863 & -0.19 & -0.10 & -0.10 & 11410 & 11541 & +1.15 & +1.09 & +1.10 \\
O2 & 9188 & 9180 & -0.09 & -0.07 & -0.07 & 9039 & 9101 & +0.68 & +0.64 & +0.64 \\
O3 & 9801 & 9785 & -0.17 & -0.14 & -0.14 & 9638 & 9695 & +0.59 & +0.67 & +0.67 \\
\midrule
E1 & 10838 & 10829 & -0.09 & -0.09 & -0.09 & 10497 & 10608 & +1.06 & +1.07 & +1.07 \\
E2 & 10259 & 10253 & -0.05 & -0.04 & -0.04 & 10046 & 10116 & +0.70 & +0.67 & +0.67 \\
E3 & 11727 & 11715 & -0.11 & -0.08 & -0.08 & 11275 & 11436 & +1.43 & +1.42 & +1.42 \\
\bottomrule
\end{tabular*}
}
\compacttableend
\end{table}

\endgroup

\end{document}